\documentclass[oneside]{article}

\usepackage[english]{babel}
\usepackage[utf8x]{inputenc}
\usepackage[T1]{fontenc}

\usepackage{amsmath, amstext, amsthm, mathrsfs,dsfont,amssymb}
\usepackage{tikz}
\usepackage{graphicx}
\usepackage[all]{xy}

\usepackage{hyperref}
\hypersetup{colorlinks = true, urlcolor = blue, citecolor=blue, bookmarksopen = true}

\definecolor{dartmouthgreen}{rgb}{0.05, 0.5, 0.06}

\theoremstyle{plain}
\newtheorem{theorem}{Theorem}
\newtheorem{hypo}{Hypothesis}

\newtheorem{lemma}{Lemma}
\newtheorem*{lemma*}{Lemma}
\newtheorem{corollary}{Corollary}
\newtheorem*{corollary*}{Corollary}
\theoremstyle{definition}

\theoremstyle{remark}

\newtheorem{remark}{Remark}
\renewenvironment{proof}{\vspace{5mm}\noindent\textbf{Proof }}{\hspace*{\fill}$\Box$\medskip\vspace{5mm}}

\def \ind {\mathds{1}}
\def \ppp {\ldots}
\def \lc {\left\lbrace}
\def \rc {\right\rbrace}

\def \uk {\underline{\kappa}}
\def \ut {\underline{\tau}}
\def \uz {\underline{z}}
\def \ud {\underline{\delta}}
\def \sd {\overline{\delta}}

\def \A {\mathbb{A}}
\def \C {\mathbb{C}}
\def \D {\mathbb{D}}
\def \E {\mathbb{E}}
\def \M {\mathbb{M}}
\def \N {\mathbb{N}}
\def \P {\mathbb{P}}

\def \R {\mathbb{R}}
\def \S {\mathbb{S}}

\def \Z {\mathbb{Z}}

\def \Dc {\mathcal{D}}

\def \Ic {\mathcal{I}}
\def \Jc {\mathcal{J}}
\def \Kc {\mathcal{K}}
\def \Lc {\mathcal{L}}
\def \Mc {\mathcal{M}}
\def \Oc {\mathcal{O}}
\def \Pc {\mathcal{P}}
\def \Qc {\mathcal{Q}}

\def \Sc {\mathcal{S}}

\def \Uc {\mathcal{U}}

\def \Acc {\mathscr{A}}
\def \Bcc {\mathscr{B}}
\def \Ccc {\mathscr{C}}

\def \Gcc {\mathscr{G}}
\def \Hcc {\mathscr{H}}
\def \Lcc {\mathscr{L}}
\def \Pcc {\mathscr{P}}

\def \Gg {\textbf{G}}

\begin{document}
	\begin{center}
		{\large \bfseries Local Limit Theorem for Complex Valued Sequences}
	\end{center}

	\begin{center}
		Lucas \textsc{Coeuret}\footnote{Institut de Mathématiques de Toulouse ; UMR5219 ; Université de Toulouse ; CNRS ; UPS, 118 route de Narbonne, F-31062 Toulouse Cedex 9, France. Research of the author was supported by the Agence Nationale de la Recherche projects Nabuco (ANR-17-CE40-0025) and Indyana (ANR-21-CE40-0008-01), as well as by the Labex Centre International de Mathématiques et Informatique de Toulouse under grant agreement ANR-11-LABX-0040. E-mail: lucas.coeuret@math-univ.toulouse.fr}
	\end{center}
	
	\vspace{5mm}
	
	\begin{center}
		\textbf{Abstract}
	\end{center}
	
	In this article, we study the pointwise asymptotic behavior of iterated convolutions on the one dimensional lattice $\Z$. We generalize the so-called local limit theorem in probability theory to complex valued sequences. A sharp rate of convergence towards an explicitly computable attractor is proved together with a generalized Gaussian bound for the asymptotic expansion up to any order of the iterated convolution.
	
	\vspace{4mm}
	
	
	
	\vspace{4mm}
	
	For $1\leq q <+\infty$, we let $\ell^q(\Z)$ denote the Banach space of complex valued sequences indexed by $\Z$ and such that the norm:
	$$\left\|u\right\|_{\ell^q}:=\left(\sum_{j\in\Z}|u_j|^q\right)^\frac{1}{q}$$
	is finite. We also let $\ell^\infty(\Z)$ denote the Banach space of bounded complex valued sequences indexed by $\Z$ equipped with the norm
	$$\left\|u\right\|_{\ell^\infty}:=\sup_{j\in \Z}|u_j|.$$
	
	Throughout this article, we define the following sets:
	$$\Uc:=\lc z\in \C, |z|>1\rc,\quad \D:=\lc z\in \C, |z|<1\rc, \quad \S^1:=\lc z\in \C, |z|=1\rc,$$
	$$\overline{\Uc}:=\S^1\cup \Uc,\quad  \overline{\D}:=\S^1\cup \D.$$
	
	For $z\in \C$ and $r>0$, we let $B_r(z)$ denote the open ball in $\C$ centered at $z$ with radius $r$.
	
	For $E$ a Banach space, we denote $\Lc(E)$ the space of bounded operators acting on $E$ and $\left\|\cdot\right\|_{\Lc(E)}$ the operator norm. For $T$ in $\Lc(E)$, the notation $\sigma(T)$ stands for the spectrum of the operator $T$.
	
	Lastly, we let $\Mc_n(\C)$ denote the space of complex valued square matrices of size $n$ and for an element $M$ of $\Mc_n(\C)$, the notation $M^T$ stands for the transpose of $M$.
	
	We use the notation $\lesssim$ to express an inequality up to a multiplicative constant. Eventually, we let $C$ (resp. $c$) denote some large (resp. small) positive constants that may vary throughout the text (sometimes within the same line).
	
	\section{Introduction and main result}
	
	\subsection{Context}
	
	We define the convolution $a\ast b$ of two elements $a$ and $b$ of $\ell^1(\Z)$ by
	$$\forall j \in \Z, \quad (a\ast b)_j:=\sum_{l\in\Z}a_lb_{j-l}.$$
	
	When equipped with this product, $\ell^1(\Z)$ is a Banach algebra. For $a \in \ell^1(\Z)$, we define the Laurent operator $L_a$ associated with $a$ which acts on $\ell^q(\Z)$ for $q\in[1,+\infty]$ as
	$$\forall u\in \ell^q(\Z), \quad L_au:=a\ast u\in \ell^q(\Z).$$
	Young's inequality implies that those operators are well defined and are bounded for all $q\in [1,+\infty]$. Furthermore, we have that $L_{a\ast b}=L_a\circ L_b$ for $a,b\in \ell^1(\Z)$. Finally, Wiener's theorem \cite{Newman} characterizes the invertible elements of $\ell^1(\Z)$ and thus allows us to describe the spectrum of $L_a$ via the Fourier series $F$ associated with $a$:
	$$\sigma(L_a)=\lc F(t):=\sum_{k\in\Z}a_k e^{itk}, t\in\R\rc.$$
	We observe that the spectrum is independent of the index $q$ and that $F$ is continuous since $a$ belongs to $\ell^1(\Z)$.
	
	If we suppose that the sequence $a$ has real nonnegative coefficients and $\sum_{k\in\Z} a_k=1$, then the sequence $a^n:=a\ast \ppp\ast a$ is the probability distribution\footnote{We say that a sequence $a$ is the probability distribution of a random variable $Y$ with values in $\Z$ when $\P(Y=j)=a_j$ for all $j\in\Z$.} of the sum of $n$ independent random variables supported on $\Z$ each with the probability distribution $a$. A lot is known on the pointwise asymptotic behavior of the sequence $a^n$ in this case. In particular, the local limit theorem states, under suitable hypotheses on the sequence $a$, that there exists a family of functions $(q_\sigma:\R\rightarrow \R)_{\sigma\in\N\backslash\lc0,1\rc}$ such that for all $s\in \N^*$ we have the following asymptotic expansion for the elements $a_j^n$
	\begin{equation}\label{cas_proba}
		a_j^n - \frac{1}{\sqrt{2\pi V n}} \exp\left(-\frac{X_{n,j}^2}{2}\right) - \sum_{\sigma=2}^s\frac{q_\sigma(X_{n,j})}{n^\frac{\sigma}{2}}\underset{n\rightarrow+\infty}=o\left(\frac{1}{n^\frac{s}{2}}\right)
	\end{equation}
	with $X_{n,j}=\frac{j-n\alpha}{\sqrt{Vn}}$ where $\alpha=\sum_{k\in\Z}ka_k$ and $V=\sum_{k\in\Z}k^2a_k-\alpha^2$ are respectively the mean and the variance of a random variable with probability distribution $a$ and where the error term is uniform with respect to $j\in\Z$ (see \cite[Chapter VII, Theorem 13]{Petrov} for more details). Furthermore, the terms in the asymptotic expansion \eqref{cas_proba} can be explicitely computed using Hermite polynomials since the functions $q_\sigma$ are explicit linear combinations of derivatives of the Gaussian function $x\mapsto \exp\left(-\frac{x^2}{2}\right)$. The asymptotic expansion \eqref{cas_proba} gives a precise description of the asymptotic behavior of $a^n_j$ in the range $|j-n\alpha|\lesssim \sqrt{n}$ and implies that the convolution powers of $a$ are attracted towards the heat kernel.
	
	Following, among other works, \cite{D-S,R-S,Cou-Faye}, we are interested in generalizing the local limit theorem to the case where $a$ is complex valued. This problem is relevant for instance when one studies the large time behavior of finite difference approximation of evolution equations. Extending the works of Schoenberg \cite{Schoenberg}, Greville \cite{Greville} and Diaconis and Saloff-Coste \cite[Theorem 2.6]{D-S}, the article \cite{R-S} of Randles and Saloff-Coste already provides a generalization of the local limit theorem for a large class of complex valued finitely supported sequences. By doing so, the authors of \cite{R-S} describe an asymptotic expansion similar to \eqref{cas_proba} for $s=1$ and identify the leading asymptotic term (the so-called "attractors" in \cite{R-S}). Our goal in this paper is to generalize the result of \cite{R-S} by obtaining an asymptotic expansion similar to \eqref{cas_proba} for any $s\in\N$ with explicitely computable terms. We also prove a sharp rate of convergence together with a generalized Gaussian bound for the remainder of our new-found asymptotic expansion (see Theorem \ref{thPrinc}). In the case where $a$ is the probability distribution of a random variable, as above, the main theorem of this paper would translate in saying that, under suitable assumptions on $a$ (namely that $a$ is finitely supported with at least two nonzero elements), for all $s\in \N^*$, there exist two constants $C,c>0$ such that
	$$\forall n\in\N^*,\forall j\in\Z,\quad \left|a^n_j -\frac{1}{\sqrt{2\pi V n}} \exp\left(-\frac{X_{n,j}^2}{2}\right) - \sum_{\sigma=2}^s\frac{q_\sigma(X_{n,j})}{n^\frac{\sigma}{2}}\right| \leq \frac{C}{n^\frac{s+1}{2}} \exp\left(-c{X_{n,j}}^2\right)$$
	with $X_{n,j}=\frac{j-n\alpha}{\sqrt{Vn}}$. As an example of application, these improvements on the local limit theorem allow us in the probabilistic case to prove the well-known Berry-Esseen inequality (see \cite{Berry,Esseen}) which states that there exists a constant $C>0$ such that
	$$\forall n\in\N^*, \forall J\in\Z, \quad \left|\sum_{j\leq J}a^n_j-\sum_{j\leq J}\frac{1}{\sqrt{2\pi V n}}\exp\left(-\frac{|j-n\alpha|^2}{2nV}\right)\right|\leq \frac{C}{\sqrt{n}}.$$
	
	However, we will need stronger hypotheses on the elements of $\ell^1(\Z)$ than the conditions imposed in \cite{R-S}. We will consider here elements $a$ of $\ell^1(\Z)$ which are finitely supported and such that the sequence $(a^n)_{n\in\N}$ is bounded in $\ell^1(\Z)$. The fundamental contribution \cite{Thomee} by Thomée completely characterizes such elements and is an important starting point for our work.
	
	In the articles \cite{D-S} and \cite{R-S}, the proofs mainly rely on the use of Fourier analysis to express the elements $a^n_j$ via the Fourier series associated with $a$. In this paper, we will rather follow an approach usually referred to in partial differential equations as "spatial dynamics". It aims at using the functional calculus (see \cite[Chapter VII]{Conway}) to express the temporal Green's function (here the coefficients $a_j^n$) with the resolvent of the operator $L_a$ via the spatial Green's function which is the unique solution of 
	$$(zId- L_a)u= \delta, \quad z \in \C\backslash\sigma(L_a),$$
	where $\delta$ is the discrete Dirac mass $\delta := (\delta_{j,0})_{j\in \Z}$. This approach has already been used in \cite{Cou-Faye} to extend the result of \cite[Theorem 1.1]{D-S} and obtain a uniform generalized Gaussian bound for the elements $a^n_j$. It has also been used in \cite{Cou-Faye2} to prove similar results on finite rank perturbations of Toeplitz operators (convolution operators on $\ell^q(\N)$ rather than on $\ell^q(\Z)$). The present paper is very much inspired by \cite{Cou-Faye,Cou-Faye2} and we will use notations and methods similar to those articles. We will now present in more details the hypotheses we need on the elements $a\in\ell^1(\Z)$ that we shall consider and we shall then present our main theorem.
	
	\subsection{Hypotheses}
	
	We consider a given sequence $a\in \ell^1(\Z)$. We let $\Lcc_a$ be the bounded operator acting on $\ell^q(\Z)$ defined as 
	$$\forall u\in \ell^q(\Z), \quad \Lcc_a u := \left(\sum_{l\in\Z} a_lu_{j+l}\right)_{j\in\Z}.$$
	This operator is obviously linked to Laurent operators and could be written as one of them ($\Lcc_a=L_b$ for $b:=(a_{-j})_{j\in\Z}$). Our goal will be to study the powers $\Lcc_a^n$ for $n$ large. This problem arises for instance as the large time behavior of finite difference approximations of partial differential equations and is equivalent to studying the asymptotic behavior of the coefficients of $b^n:=b\ast\ppp\ast b$ as $n$ tends to infinity. We define the symbol $F$ associated with $a$ as 
	\begin{equation}\label{def_f}
		\forall \kappa \in \S^1, \quad F(\kappa):= \sum_{j\in \Z} a_j \kappa^j.
	\end{equation}
	The Wiener theorem \cite{Newman} allows us to conclude that the spectrum of $\Lcc_a$ is given, for any $q\in[1, +\infty]$, by:
	$$\sigma(\Lcc_a)=F(\S^1).$$
	
	We are now going to introduce some hypotheses that are necessary for the rest of the paper.
	
	\begin{hypo}\label{H1}
		The sequence $a$ is finitely supported and has at least two nonzero coefficients.
	\end{hypo}

	Looking at the definition of the operator $\Lcc_a$, in terms of applications for numerical analysis, this hypothesis translates the fact that we are only considering the case of explicit finite difference schemes. Hypothesis \ref{H1} implies that we can extend the definition \eqref{def_f} of $F$ to the pointed plane $\C\backslash\lc0\rc$ and $F$ becomes a holomorphic function on this domain. We introduce the two following elements 
	$$k_m:=\min \lc k\in\Z,\quad a_k\neq0\rc,\quad		k_M:=\max \lc k\in\Z,\quad a_k\neq0\rc.$$
	Observing that Hypothesis \ref{H1} implies $k_m<k_M$, we then distinguish three different possibilities:
	\begin{itemize}
		\item Case 1: $k_M\leq-1$. We then define $r:=-k_m$ and $p:=0$.
		\item Case 2: $k_m\leq 0 \leq k_M$. We then define $r:=-k_m$ and $p:=k_M$.
		\item Case 3: $1\leq k_m$. We then define $r:=0$ and $p:=k_M$.
	\end{itemize}
	In every case, we have $r,p\in\N$ and $-r<p$. Also, we have that 
	\begin{equation}\label{expLa}
		\forall u\in\ell^q(\Z),\forall j\in\Z, \quad (\Lcc_a u)_j=\sum_{l=-r}^pa_lu_{j+l}.
	\end{equation}
	The natural integers $r$ and $p$ we just introduced define the common stencil of the operators $\Lcc_a$ and the identity operator and they will be useful to study the so-called resolvent equation \eqref{def_G} below. We now introduce an assumption on the Laurent series $F$ which is based on \cite{Thomee}. Just like in \cite{D-S,R-S,Cou-Faye}, we normalize the sequence $a$ so that the maximum of $F$ on $\S^1$ is $1$. 
	
	\begin{hypo}\label{H2}
		There exists a finite set of distinct points $\lc\uk_1,\ppp,\uk_K\rc$, $K\geq 1$, in $\S^1$ such that for all $k\in\lc1,\ppp,K\rc$, $\uz_k:=F(\uk_k)$ belongs to $\S^1$ and
		$$\forall \kappa\in \S^1\backslash \lc\uk_1,\ppp,\uk_K\rc, \quad |F(\kappa)|<1.$$
		Moreover, we suppose that for each $k\in\lc1,\ppp,K\rc$, there exist a nonzero real number $\alpha_k$, an integer $\mu_k\geq1$ and a complex number $\beta_k$ with positive real part such that
		\begin{equation}\label{F}
			F(\uk_ke^{i\xi})\underset{\xi\rightarrow0}= \uz_k\exp(-i\alpha_k \xi - \beta_k \xi^{2\mu_k} + O(|\xi|^{2\mu_k+1})).
		\end{equation}
	\end{hypo}
	
	\begin{figure}
		\begin{center}
			\begin{tikzpicture}[scale=2]
				\fill[color=gray!20] (-1.5,-1.5) -- (-1.5,1.5) -- (1.5,1.5) -- (1.5,-1.5) -- cycle;
				\draw[color=black!60] (-1,1) node {$\Oc$};
				\fill[color=white] plot [samples = 100, domain=0:2*pi] ({cos(\x r)},{sin(\x r)/2}) -- cycle ;
				\fill[color=white] (0,0) circle (0.7) ;
				
				\draw (0,0) circle (1);
				\draw[dashed] (0,0) circle (0.7);
				\draw (0.2,-0.2) node {$\exp(-\underline{\eta})\S^1$};
				\draw (45:1.2) node {$\S^1$};
				\draw[thick,red] plot [samples = 100, domain=0:2*pi] ({cos(\x r)},{sin(\x r)/2});
				\draw[red] (0,0.3) node {$\sigma(\Lcc_a)$};
				\draw[blue] (1,0) node {$\bullet$} node[right] {$\uz_1$};
				\draw[blue] (-1,0) node {$\bullet$} node[left] {$\uz_2$};
			\end{tikzpicture}
			\caption{An example of spectrum $\sigma(\Lcc_a)$. The spectrum $\sigma(\Lcc_a)$ (in red) is inside the closed disk $\bar{\D}$ and touches the boundary $\S^1$ in finitely many points. In gray, we have $\Oc$ the intersection of the unbounded connected component of $\C\backslash\sigma(\Lcc_a)$ and $\lc z\in\C, |z|>\exp(-\underline{\eta})\rc$.}
			\label{spec}
		\end{center}
	\end{figure}
	
	Geometrically, this means that the spectrum $\sigma(\Lcc_a)$ is contained in the disk $\overline{\D}$ and it intersects $\S^1$ at finitely many points (see Figure \ref{spec} for an example with $K=2$, $\uz_1=1$, $\uz_2=-1$) and that the logarithm of $F$ has a specific asymptotic expansion at those intersection points. From a general point of view, it is proved in \cite[Theorem 1]{Thomee} that Hypothesis \ref{H2} is one of two conditions that characterize the elements $a$ of $\ell^1(\Z)$ such that the geometric sequence $(a^n)_{n\in\N}$ is bounded in $\ell^1(\Z)$. In the more specific field of numerical analysis, the condition \eqref{F} has been studied closely because of its link with the stability of finite difference approximations in the maximum norm (see \cite{Thomee}). We can observe that, under Hypotheses \ref{H1} and \ref{H2}, there holds
	$$\forall n\in\N^*,\quad \left\|\Lcc_a^n\right\|_{\Lc(\ell^2(\Z))}=\left\|F^n\right\|_{L^\infty(\S^1)}=1.$$
	It assures the $\ell^2$-stability, or strong stability (see \cite{Strang}, \cite{Tad}), of the numerical scheme defined as
	\begin{equation}\label{an_num}\lc
		\begin{array}{cc}
			u^{n+1}=\Lcc_a u^n,& n\geq0,\\ 
			u^0\in \ell^2(\Z). &
		\end{array}\right.
	\end{equation}
	However, it has further consequences, as the asymptotic expansion \eqref{F} assures the $\ell^q$-stability of the scheme \eqref{an_num} for every $q$ in $[1,+\infty]$ (see \cite[Theorem 1]{Thomee} which focuses on the $\ell^\infty$-stability but also studies the $\ell^q$-stability as a consequence). In terms of numerical scheme, the meaning of \eqref{F} is that the numerical scheme introduces an artificial numerical diffusion (like the Lax-Friedrichs scheme for example).   
	
	We now introduce yet another hypothesis.
	
	\begin{hypo}\label{H3}
		For all $k\in\lc1,\ppp,K\rc$, the set 
		$$\Ic_k:= \lc\nu \in\lc1,\ppp,K\rc, \quad \uz_\nu=\uz_k\rc$$
		has either one or two elements, where we recall that $\uz_\nu:=F(\uk_\nu)$. Moreover, if there are two distinct elements $\nu_{k,1}$ and $\nu_{k,2}$ in $\Ic_k$, then $\alpha_{\nu_{k,1}}\alpha_{\nu_{k,2}}<0$. 
	\end{hypo}
	
	Hypothesis \ref{H3} will simplify part of the analysis when we will study the spatial Green's function defined in \eqref{def_G} below. It will allow us to study precisely the spectrum of the matrix $\M(z)$ defined below as \eqref{def_M} near the tangency points $\uz_k$. Combining Hypothesis \ref{H3} with the fact that the $\alpha_k$'s are nonzero real numbers (see Hypothesis \ref{H2}) implies that, for $k\in\lc1,\ppp,K\rc$, we have three different possibilities:
	
	\begin{itemize}
		\item \textbf{Case I:} $\Ic_k$ is the singleton $\lc k\rc$ and $\alpha_k>0$,
		\item \textbf{Case II:} $\Ic_k$ is the singleton $\lc k\rc$ and $\alpha_k<0$,
		\item \textbf{Case III:} $\Ic_k$ has two distinct elements $\nu_{k,1}$ and $\nu_{k,2}$ such that $\alpha_{\nu_{k,1}}>0$ and $\alpha_{\nu_{k,2}}<0$.
	\end{itemize} 
	
	Distinguishing between those three cases will be useful later on. The three hypotheses we presented above will be crucial in the rest of the paper. Some hypotheses might be relaxable, but this would be considerations for future works.	
	
	Finally, by defining the discrete Dirac mass $\delta := (\delta_{j,0})_{j\in \Z}$, we introduce the so-called temporal Green's function defined by 
	\begin{equation}\label{defGreenTempo}
		\forall n\in\N, \forall j\in \Z, \quad \Gcc_j^n:=\left(\Lcc_a^n\delta\right)_j.
	\end{equation}
	It is interesting to observe that the equality between the operator $\Lcc_a$ and the Laurent operator $L_b$ with $b=\left(a_{-j}\right)_{j\in\Z}$ implies that
	$$\forall n\in\N, \forall j\in \Z, \quad \Gcc_j^n=b^n_j$$
	where $b^n=b\ast \ppp\ast b$.
	\subsection{Main results and comparison to previous results}
	
	Our main goal is to determine the asymptotic behavior of $\Gcc_j^n$ when $n$ becomes large. The identification of the leading asymptotic term was achieved in \cite[Theorem 1.2]{R-S}. We aim here at extending the result of \cite[Theorem 1.2]{R-S} into a complete asymptotic expansion up to any order and at proving sharp bounds for the remainder. To express the asymptotic expansion of $\Gcc_j^n$, we introduce the functions $H_{2\mu}^\beta: \R\rightarrow \C$, where $\mu\in\N^*$ and $\beta \in \C$ has positive real part, which are defined as
	$$\forall x \in \R,\quad H_{2\mu}^\beta(x) := \frac{1}{2\pi} \int_\R e^{ixu}e^{-\beta u^{2\mu}}du.$$
	We call those functions generalized Gaussians since for $\mu=1$, we have 
	$$\forall x \in \R,\quad H_{2}^\beta(x)=\frac{1}{\sqrt{4\pi\beta}}e^{-\frac{x^2}{4\beta}}.$$
	
	Let us state the main result of this paper.
	
	\begin{theorem}\label{thPrinc}
		Let $a\in \ell^1(\Z)$ which verifies Hypotheses \ref{H1}, \ref{H2} and \ref{H3}. Then, for all integers $s_1,\ppp,s_K\in\N$ there exist a family of polynomials $(\Pcc^k_\sigma)_{\sigma\in\lc1, \ppp,s_k\rc}$ in $\C[X,Y]$ for each $k\in\lc1, \ppp,K\rc$ and two positive constants $C,c$ such that for all $n\in\N^*$ and $j\in \Z$, there holds:
		\begin{equation}\label{devPrinc}
			\left|\Gcc_j^n-\sum_{k=1}^K\sum_{\sigma=1}^{s_k}\frac{\uz_k^n\uk_k^j}{n^\frac{\sigma}{2\mu_{k}}}\left(\Pcc^k_\sigma\left(X_{n,j,k},\frac{d}{dx}\right)H_{2\mu_k}^{\beta_k}\right)\left(X_{n,j,k}\right)\right|\leq\sum_{k=1}^K \frac{C}{n^\frac{s_k+1}{2\mu_k}}\exp\left(-c|X_{n,j,k}|^\frac{2\mu_k}{2\mu_k-1}\right)
		\end{equation}
		where $X_{n,j,k}=\frac{n\alpha_k-j}{n^\frac{1}{2\mu_k}}$.
	\end{theorem}
	
	Theorem \ref{thPrinc} gives the asymptotic behavior of the elements $\Gcc_j^n$ up to any order with a sharp generalized Gaussian estimate of the remainder. We would also like to point out that the proof of Theorem \ref{thPrinc} (mainly Lemmas \ref{lemEstInterm}, \ref{lemDevAsymp} and equality \eqref{egHjalp}) gives us an explicit expression of the polynomials $\Pcc^k_\sigma$ of Theorem \ref{thPrinc}. Examples are provided in Section \ref{secComput} where we compute these polynomials for $\sigma=1,2$ and numerically verify the claim of Theorem \ref{thPrinc} for some sequences $a$.
	
	The following lemma, which is proved using integration by parts, implies that we cannot prove the uniqueness of the polynomials $\Pcc^k_\sigma$ of Theorem \ref{thPrinc}.
	\begin{lemma}\label{lemH}
		For $\mu\in\N^*$, $\beta\in \C$ with positive real part and $m\in \N^*$, we have
		$$\forall x\in\R,\quad x{H_{2\mu}^\beta}^{(m)}(x)= (-1)^\mu2\mu\beta {H_{2\mu}^\beta}^{(m+2\mu-1)}(x)-m{H_{2\mu}^\beta}^{(m-1)}(x),$$
		and
		$$\forall x\in\R,\quad xH_{2\mu}^\beta(x)= (-1)^\mu2\mu\beta {H_{2\mu}^\beta}^{(2\mu-1)}(x).$$
	\end{lemma}
	In other words, one can either choose to multiply $H_{2\mu}^{\beta}$ by a polynomial or to differentiate it sufficiently many times. Hence, there may hold 
	$$P(\cdot, \frac{d}{dx})H_{2\mu}^\beta=0$$
	for a nonzero $P\in \C[X,Y]$.
	
	In our proof of Theorem \ref{thPrinc}, the polynomials $\Pcc^k_\sigma$ depend on the chosen integers $s_1, \ppp,s_k$. It might be possible to prove the existence of a family of polynomials $\left(\Pcc^k_{\sigma}\right)_{(k,\sigma)\in\lc 1,\ppp,K\rc \times \N^*}$ in $\C[X,Y]$ for which the estimates \eqref{devPrinc} are verified for all $s_1,\ppp,s_K\in\N$. However, we do not yet have a proof of this fact in full generality. We now compare Theorem \ref{thPrinc} with prior results:\newline
	
	$\bullet$ In the probabilistic case presented in the introduction, Theorem \ref{thPrinc} allows us to prove sharp bounds with Gaussian estimates on the remainder of the asymptotic expansion of $\Gcc_j^n$ that were not proved via the asymptotic expansion \eqref{cas_proba} of the local limit theorem. \newline
	
	$\bullet$ \cite[Theorem 3.1]{D-S} gives sharp generalized Gaussian estimates for the elements $\Gcc_j^n$ when the sequence $a$ satisfies Hypotheses \ref{H1}, \ref{H2} and \ref{H3} with a single tangency point (i.e. $K=1$), which in comparison to Theorem \ref{thPrinc} would match the case $s_k=0$.  \cite[Theorem 1]{Cou-Faye} generalizes those generalized Gaussian estimates for sequences $a$ with any number $K\in \N^*$ of tangency points and a relaxed Hypothesis \ref{H1}. Theorem \ref{thPrinc} thus improves those results by proving similar sharp generalized Gaussian estimates for the remainder of the asymptotic expansion of the elements $\Gcc_j^n$ up to any order $s_1,\ppp,s_K\in\N$.\newline
		
	$\bullet$ For a sequence $a\in\ell^1(\Z)$ which satisfies Hypotheses \ref{H1} and \ref{H2}, we introduce the the so-called "attractors":
	$$\forall k\in\lc1,\ppp, K\rc, \forall n\in\N^*, \forall j\in \Z, \quad \Hcc_{k,j}^n := \frac{{\uz_k}^n{\uk_k}^j}{n^\frac{1}{2\mu_k}}H_{2\mu_k}^{\beta_k}\left(\frac{j-n\alpha_k}{n^\frac{1}{2\mu_k}}\right).$$
	
	In \cite[Theorem 1.2]{R-S}, it is proved that if we introduce $\Kc_\mu=\lc k\in\lc1,\ppp,K\rc, \mu_k=\mu\rc$ where $\mu=\max_{k\in\lc1,\ppp,K\rc}\mu_k$, then
	\begin{equation}\label{eq_R-S}
		\Gcc_j^n-\sum_{k\in\Kc_\mu}\Hcc_{k,j}^n\underset{n\rightarrow +\infty}= o\left(\frac{1}{n^\frac{1}{2\mu}}\right)
	\end{equation}
	where the error term in \eqref{eq_R-S} is uniform on $\Z$. Compared to Theorem \ref{thPrinc}, this is equivalent to finding the asymptotic expansion up to order $s_1,\ppp,s_K=1$. The result of Randles and Saloff-Coste gives a precise description of the behavior of $\Gcc^n_j$ for $j$ such that
	\begin{equation}\label{dom}
		\left|j-n\alpha_k\right|\lesssim n^\frac{1}{2\mu},
	\end{equation}
	where $k\in \Kc_\mu$. Theorem \ref{thPrinc} allows us to extend the result of \cite{R-S} by going even farther in the asymptotic expansion of the elements $\Gcc_j^n$, and proving sharp generalized Gaussian bounds on the remainder with a more precise speed of convergence. However, \cite[Theorem 1.2]{R-S} also treats the case where the asymptotic expansion \eqref{F} has the form
	$$F(\uk_ke^{i\xi})\underset{\xi\rightarrow0}= \uz_k\exp(-i\alpha_k \xi +i \gamma_k \xi^{\nu_k} + O(|\xi|^{\nu_k+1})),$$
	where $\gamma_k$ is a real number and the integer $\nu_k\in \N\backslash\lc 0,1\rc$ can be even or odd. A generalization of Theorem \ref{thPrinc} in this difficult case has not yet been found, even though the result of \cite{JFCLW} indicates that such a result might be attainable.
	
	\subsection{Extending the result when the drift vanishes}
	
	As we have seen, Theorem \ref{thPrinc} allows us to have generalize the local limit theorem for complex valued sequences but it still has some limits. Relaxing some of the hypotheses we made could be interesting and theoretically doable in some cases. For example, Theorem \ref{thPrinc} is constrained by Hypothesis \ref{H2} which imposes that $\alpha_k$ is nonzero even though the result \cite[Theorem 1.2]{R-S} does not have this kind of restriction. The hypothesis $\alpha_k\neq0$ is essential in the proof of Theorem \ref{thPrinc} below but it seems to be a technical hypothesis that we would want to avoid. The following corollary will allow us to extend Theorem \ref{thPrinc} to some sequences $a$ for which we allow $\alpha_k$ to be equal to $0$. First, we introduce a relaxed version of Hypothesis \ref{H2}.
	
	\begin{hypo}[Hypothesis 2 bis]\label{H2_bis}
		The sequence $a$ verifies Hypothesis \ref{H2} but with the possibility that some $\alpha_k$ are equal to $0$.
	\end{hypo}
	
	We now consider a finitely supported sequence $a\in\ell^1(\Z)$ which verifies Hypothesis \ref{H2_bis} and let $J\in\Z$. Then, if we define the sequence $b=(a_{j+J})_{j\in\Z}$ and $\widetilde{F}$ the symbol associated with $b$, we have that $b$ satisfies Hypothesis \ref{H2_bis} since 
	$$\forall \kappa\in\S^1, \quad \widetilde{F}(\kappa)=\kappa^{-J}F(\kappa),$$
	and therefore
	$$\forall \kappa\in\S^1, \quad \left|\widetilde{F}(\kappa)\right|=\left|F(\kappa)\right|.$$
	Also, we have for $k\in\lc1,\ppp,K\rc$
	\begin{equation*}
		\widetilde{F}(\uk_ke^{i\xi})\underset{\xi\rightarrow0}= \uk_k^{-J}\uz_k\exp(-i(\alpha_k+J) \xi - \beta_k \xi^{2\mu_k} + o(|\xi|^{2\mu_k})).
	\end{equation*}
	Considering this new sequence $b$ allows us to "shift" the elements $\alpha_k$. In particular, if we choose $J$ large enough, then $b$ satisfies Hypothesis \ref{H2}. However, it is not clear that the sequence $b$ would satisfy Hypothesis \ref{H3}. We can then prove the following corollary of Theorem \ref{thPrinc} which generalizes Theorem \ref{thPrinc} in the case where $\alpha_k$ can be equal to $0$.
	
	\begin{corollary}\label{cor_prin}
		Let $a\in \ell^1(\Z)$ which verifies Hypotheses \ref{H1} and \ref{H2_bis}. If there exists some integer $J\in\Z$ such that the sequence $(a_{j+J})_{j\in\Z}$ verifies Hypotheses \ref{H2} and \ref{H3}, then for all $s_1,\ppp,s_K\in\N$ there exist a family of polynomials $(\Pcc^k_\sigma)_{\sigma\in\lc1, \ppp,s_k\rc}$ in $\C[X,Y]$ for each $k\in\lc1, \ppp,K\rc$ and two positive constants $C,c$ such that for all $n\in\N^*$ and $j\in \Z$
		$$\left|\Gcc_j^n-\sum_{k=1}^K\sum_{\sigma=1}^{s_k}\frac{\uz_k^n\uk_k^j}{n^\frac{\sigma}{2\mu_{k}}}\left(\Pcc^k_\sigma\left(X_{n,j,k},\frac{d}{dx}\right)H_{2\mu_k}^{\beta_k}\right)\left(X_{n,j,k}\right)\right|\leq\sum_{k=1}^K \frac{C}{n^\frac{s_k+1}{2\mu_k}}\exp\left(-c|X_{n,j,k}|^\frac{2\mu_k}{2\mu_k-1}\right)$$
		with $X_{n,j,k}=\frac{n\alpha_k-j}{n^\frac{1}{2\mu_k}}$.
	\end{corollary}

	We prove Corollary \ref{cor_prin} in Section \ref{secCor}.
	
	\subsection{Plan of the paper}
	
	The main goal of the paper is the proof of Theorem \ref{thPrinc}. As explained in the introduction, the proof of Theorem \ref{thPrinc} will rely on an approach referred to as spatial dynamics. In Section \ref{sec_GS}, we will introduce the spatial Green's function on which Coulombel and Faye proved holomorphic extension properties and sharp bounds in \cite[Section 2]{Cou-Faye}. Our goal in Section \ref{sec_GS} is to improve the analysis of \cite{Cou-Faye} and to obtain the precise behavior of the spatial Green's function for $z$ close to $\uz_k$ and to prove sharp bounds on the remainder. More precisely, the main novelty of this section is the introduction of the explicit function $f_k$ in Lemmas \ref{green_spatial_près_1} and \ref{green_spatial_près_2} which allows us to properly describe the spatial Green's function for $z$ close to $\uz_k$.
	
	In Section \ref{sec_GT}, we prove Theorem \ref{thPrinc} while assuming that the elements $\alpha_k$ are distinct. This assumption will allow us to separate the different Gaussian waves in the estimate \eqref{devPrinc}. Section \ref{subsec_est_far} will be dedicated to the easier part of the proof which is proving estimate \eqref{devPrinc} when $j$ is far from the axes $j=n\alpha_k$. The bulk of the proof resides in Sections \ref{subsecPlan}-\ref{subsecLemImpTh} which will be dedicated to proving estimate \eqref{devPrinc} when $j$ is close to the axes $j=n\alpha_k$. In Section \ref{sec_GS/GT}, we will express the elements $\Gcc_j^n$ with the spatial Green's function using functional calculus. We will then use the results of Section \ref{sec_GS} on the spatial Green's function to prove generalized Gaussian estimates on the difference of the elements $\Gcc_j^n$  and a linear combination of terms of the form
	\begin{equation}\label{termesY}
		\frac{1}{\left(\frac{j}{\alpha_k}\right)^\frac{l}{2\mu_k}}{H_{2\mu_k}^{\beta_k}}^{(m)}\left(Y_{n,j,k}\right)\quad \text{where }l\in\N^*, m\in\N,Y_{n,j,k}:=\frac{n\alpha_k-j}{\left(\frac{j}{\alpha_k}\right)^\frac{1}{2\mu_k}}.
	\end{equation}
	Keeping in mind that we are considering the case where $j$ is close to $n\alpha_k$, Section \ref{secAsymp} will deal with approaching the terms \eqref{termesY} with linear combinations of the following terms appearing in Theorem \ref{thPrinc}:
	$$\frac{1}{n^\frac{l}{2\mu_k}}\left(X_{n,j,k}\right)^{m_2}{H_{2\mu_k}^{\beta_k}}^{(m_1)}\left(X_{n,j,k}\right)\quad \text{where }l\in\N^*, m_1,m_2\in\N,X_{n,j,k}:=\frac{n\alpha_k-j}{n^\frac{1}{2\mu_k}}.$$ 
	Section \ref{subsecLemImpTh} will combine the results of the previous sections to conclude the proof of Theorem \ref{thPrinc} by constructing the polynomials $\Pc^k_\sigma$.
	
	In Section \ref{secConcluTh}, we prove Theorem \ref{thPrinc} when the elements $\alpha_k$ can be equal. We also prove Corollary \ref{cor_prin}.
	
	Finally, in Section \ref{secComput}, we will explicitly compute the polynomials $\Pcc^k_\sigma$ of Theorem \ref{thPrinc} for $\sigma=1,2$ for any $s_k\in \N\backslash\lc0,1\rc$ and numerically verify the estimate \eqref{devPrinc} of Theorem \ref{thPrinc} in two cases. The first one is the probabilistic case, i.e. a sequence $a$ with non negative coefficients. We will compare the result of Theorem \ref{thPrinc} with the local limit theorem. The second example will be the sequence $a$ associated with the so-called O3 scheme for the transport equation (see \cite{Despres}). This is an example of sequence $a$ where $\mu=2$ in the asymptotic expansion \eqref{F}.
	
	\section{Spatial Green's function}\label{sec_GS}
	
	From now on, we consider a sequence $a$ that satisfies Hypotheses \ref{H1}, \ref{H2} and \ref{H3}. In this section, we are going to introduce the spatial Green's function and prove some estimates for it. We will start by defining the necessary objects for our study. First, we can observe the following lemma for which the proof can be found in the Appendix (Section \ref{sec_appendix}).
	\begin{lemma}\label{arp}
		For $a\in\ell^1(\Z)$ which verifies Hypotheses \ref{H1} and \ref{H2}, we have that $a_{-r}$ and $a_p$ belong to $\D$.
	\end{lemma}
	
	We define for $z\in \C$ and $j\in \lc-r,\ppp,p\rc$
	\begin{equation}
		\A_j(z):= z \delta_{j,0}-a_j.\label{def_Aj}
	\end{equation}
	
	The definition of $r$ and $p$ implies that the functions $\A_{-r}$ and $\A_p$ can vanish at most on one point which are respectively $a_{-r}$ and $a_p$. Lemma \ref{arp} allows us to find $\underline{\eta}>0$ such that $\A_{-r}$ and $\A_p$ do not vanish on $\lc z \in \C, |z|>\exp(-\underline{\eta})\rc.$ We can therefore define for all $z \in \C$ such that $|z|>\exp(-\underline{\eta})$ the matrix
	\begin{equation}\label{def_M}\M(z):= \begin{pmatrix}
		-\frac{\A_{p-1}(z)}{\A_p(z)} & \ppp & \ppp & -\frac{\A_{-r}(z)}{\A_p(z)} \\
		1 & 0 & \ppp & 0\\
		0 & \ddots & \ddots & \vdots\\
		0& 0 & 1& 0
	\end{pmatrix} \in \Mc_{p+r}(\C).\end{equation}
	
	The application which associates $z$ with $\M(z)$ is holomorphic on the annulus $\lc z \in \C, |z|>\exp(-\underline{\eta})\rc.$ Moreover, since $\A_{-r}(z)\neq0$, the upper right coefficient of $\M(z)$ is always nonzero and $\M(z)$ is invertible. We define the open set $\Oc$ which corresponds to the intersection of the unbounded connected component of $\C\backslash F(\S^1)$ and $\lc z\in\C, |z|>\exp(-\underline{\eta})\rc$ (see Figure \ref{spec}).	Hypothesis \ref{H2} implies that $\overline{\Uc}\backslash\lc \uz_1,\ppp,\uz_K\rc$ is contained within $\Oc$.  By recalling that $\sigma(\Lcc_a)=F(\S^1)$, when we consider that $\Lcc_a$ acts on $\ell^2(\Z)$, we have the existence for every $z\in \Oc$ of a unique sequence $G(z):=(G_j(z))_{j\in\Z}\in \ell^2(\Z)$ such that
	\begin{equation}
		(zI-\Lcc_a)G(z)=\delta, \label{def_G}
	\end{equation}
	where $\delta$ still denotes the discrete Dirac mass. The sequence $G(z)$ is the so-called spatial Green's function which has already been studied in \cite{Cou-Faye}. In \cite[Lemma 2]{Cou-Faye}, we can find a proof of local sharp exponential bounds on $G_j(z)$ when $z\in \Oc$ is far from the tangency points $\uz_k$. This bound will be sufficient for our purpose. Furthermore, in \cite[Lemmas 3 and 4]{Cou-Faye}, the authors proved that the spatial Green's function $G_j(z)$ could be holomorphically extended near the points $\uz_k$ through the spectrum of the operator $\Lcc_a$ which is not immediate based on the definition \eqref{def_G} of the spatial Green's function and they proved sharp bounds on $G_j(z)$ in this case. To prove Theorem \ref{thPrinc}, we will need to get a more precise description of the behavior of the sequence $G_j(z)$ close to any tangency point $\uz_k$. This section will therefore follow \cite[Section 2]{Cou-Faye} and make it more precise by specifying where our study of the sequence $G(z)$ differs from \cite[Section 2]{Cou-Faye}.

	Using the functions $\A_l$ which are defined by \eqref{def_Aj}, the equation \eqref{def_G} can be rewritten as 
	$$\forall z\in\Oc,\forall j \in \Z, \quad \sum_{l=-r}^p\A_l(z)G_{j+l}(z)=\delta_{j,0}.$$
	
	We introduce the vectors
	$$\forall z\in\Oc,\forall j\in \Z, \quad W_j(z) := \begin{pmatrix}
		G_{j+p-1}(z) \\ \vdots \\ G_{j-r}(z)
	\end{pmatrix}\in\C^{p+r}, \quad  e := \begin{pmatrix}
		1\\0\\ \vdots \\ 0
	\end{pmatrix}\in \C^{p+r}.$$
	
	We then end up with the following dynamical system
	\begin{equation}
		\forall z\in \Oc, \forall j\in \Z, \quad W_{j+1}(z)-\M(z)W_j(z) = -\frac{\delta_{j,0}}{\A_p(z)}e. \label{Wj}
	\end{equation}
	
	The study of the recurrence relation \eqref{Wj} relies on the following lemma introduced in \cite{Kre} that studies the eigenvalues of $\M(z)$ for $z\in \Oc$ and $z\in\lc\uz_k, 1\leq k\leq K\rc$. We recall that we defined cases \textbf{I}, \textbf{II} and \textbf{III} according to the cardinality of $\Ic_k$ and the sign of $\alpha_k$ right after Hypothesis \ref{H3}. We also recall that we consider that the sequence $a$ verifies Hypotheses \ref{H1}, \ref{H2} and \ref{H3}.
	
	\begin{lemma}[Spectral Splitting]\label{spec_spl}
		For $z\in\C$ such that $|z|>\exp(-\underline{\eta})$, the eigenvalues $\kappa\in\C$ of the matrix $\M(z)$ are nonzero and satisfy the equality
		$$F(\kappa)=z.$$		
		Let $z\in \Oc$. Then the matrix $\M(z)$ has
		\begin{itemize}
			\item no eigenvalue on $\S^1$,
			\item $r$ eigenvalues in $\D\backslash\lc0\rc$ (that we call stable eigenvalues),
			\item  $p$ eigenvalues in $\Uc$ (that we call unstable eigenvalues).
		\end{itemize} 
		
		We now consider $k\in\lc1,\ppp,K\rc$. The eigenvalues of the matrix $\M(\uz_k)$ are described by the following possibilities depending on $k$.
		\begin{itemize}
			\item In case \textbf{I}, $\M(\uz_k)$ has $\uk_k\in\S^1$ as a simple eigenvalue, $r-1$ eigenvalues in $\D$ and $p$ eigenvalues in $\Uc$.
			\item In case \textbf{II}, $\M(\uz_k)$ has $\uk_k\in\S^1$ as a simple eigenvalue, $r$ eigenvalues in $\D$ and $p-1$ eigenvalues in $\Uc$.
			\item In case \textbf{III}, if we denote $\nu_{k,1}$ and $\nu_{k,2}$ the two distinct elements of $\Ic_k$, then $\M(\uz_k)$ has $\uk_{\nu_{k,1}}\in\S^1$ and $\uk_{\nu_{k,2}}\in\S^1$ as simple eigenvalues, $r-1$ eigenvalues in $\D$ and $p-1$ eigenvalues in $\Uc$.
		\end{itemize} 
	\end{lemma}
	
	Lemma \ref{spec_spl} is proved in \cite[Lemma 1]{Cou-Faye} and is the key to study the recurrence relation \eqref{Wj}. We now want to prove some estimates on the spatial Green's function $G(z)$. We recall that the set $\Oc$ is the intersection of the set $\lc z\in \C, |z|\geq\exp(-\underline{\eta})\rc$, where the matrix $\M(z)$ is defined, and the set $\sigma(\Lcc)=\C\backslash F(\S^1)$, where the spatial Green's function $G(z)$ is defined. We begin with the following lemma.
	
	\begin{lemma}[Bounds far from the tangency points \cite{Cou-Faye}]
		For all $\underline{z}\in \Oc$, there exist a radius $\delta>0$ and constants $C,c>0$ such that for all $j\in \Z$, $z\mapsto G_j(z)$ is holomorphic on $B_\delta(\underline{z})$ and satisfies
		$$\forall z\in B_\delta(\underline{z}), \forall j \in \Z, \quad |G_j(z)|\leq C\exp(-c|j|).$$
		\label{green_spatial_loin}
	\end{lemma}
	
	Lemma \ref{green_spatial_loin} is proved in \cite[Lemma 2]{Cou-Faye} and allows us to study the spatial Green's function far from the points $\uz_k$, where the spectrum of $\Lcc_a$ intersects the unit circle $\S^1$. We will now have to study the spatial Green's function $G(z)$ near those points $\uz_k$ while still remembering that $G_j(z)$ and the vector $W_j(z)$ are only defined on $\Oc$ in the neighborhood of $\uz_k$. We are going to extend holomorphically $G_j(z)$ in a whole neighborhood of $\uz_k$, and thus pass through the spectrum $\sigma(\Lcc_a)$.
	
	\begin{lemma}[Bounds close to the tangency points : cases\textbf{ I} and \textbf{II}]\label{green_spatial_près_1}
		Let $k\in\lc1,\ppp,K\rc$ so that we are either in case \textbf{I} or \textbf{II}. Then, there exist a radius $\varepsilon>0$, some constants $C,c>0$ and some holomorphic functions $\kappa_k,f_k:B_\varepsilon(\uz_k)\rightarrow \C$ such that for all $z \in B_\varepsilon(\uz_k)$, $\kappa_k(z)$ is a simple eigenvalue of $\M(z)$ with $\kappa_k(\uz_k)=\uk_k$, for all $j\in \Z$, the function $z\in B_\varepsilon(\uz_k)\cap \Oc \mapsto G_j(z)$  can be holomorphically extended on $B_\varepsilon(\uz_k)$ and
		
		\textbf{Case I: ($\alpha_k>0$)}
		\begin{align}
			\forall z\in B_\varepsilon(\uz_k), \forall j \geq 1, \quad |G_j(z)- f_k(z) \kappa_k(z)^j|\leq C\exp(-cj).\label{est_spa_B_cas1}\\
			\forall z\in B_\varepsilon(\uz_k), \forall j \leq0, \quad |G_j(z)|\leq C\exp(-c|j|).\label{est_spa_B_2_cas1}
		\end{align}
		
		\textbf{Case II: ($\alpha_k<0$)}
		\begin{align}
			\forall z\in B_\varepsilon(\uz_k), \forall j \geq 1, \quad |G_j(z)|\leq C\exp(-cj).\label{est_spa_B_cas2}\\
			\forall z\in B_\varepsilon(\uz_k), \forall j \leq0, \quad |G_j(z)- f_k(z) \kappa_k(z)^j|\leq C\exp(-c|j|).\label{est_spa_B_2_cas2}
		\end{align}
		Furthermore, we have
		\begin{equation}\label{fk_1}
			\forall z\in B_\varepsilon(\uz_k), \quad f_k(z)=-\mathrm{sgn}(\alpha_k)\frac{\kappa_k^\prime(z)}{\kappa_k(z)}.
		\end{equation}
	\end{lemma}
	
	\begin{lemma}[Bounds close to the tangency points : case\textbf{ III}]\label{green_spatial_près_2}
		Let $k\in\lc1,\ppp,K\rc$ so that we are in case\textbf{ III}. The set $\Ic_k$ has two elements $\nu_{k,1}$ and $\nu_{k,2}$ so that $\alpha_{\nu_{k,1}}>0$ and $\alpha_{\nu_{k,2}}<0$. Then, there exist a radius $\varepsilon>0$, some constants $C,c>0$ and some holomorphic functions $\kappa_{\nu_{k,1}},\kappa_{\nu_{k,2}},f_{\nu_{k,1}},f_{\nu_{k,2}}:B_\varepsilon(\uz_k)\rightarrow \C$ such that for all $z \in B_\varepsilon(\uz_k)$, $\kappa_{\nu_{k,1}}(z)$ and $\kappa_{\nu_{k,2}}(z)$ are simple eigenvalues of $\M(z)$ with $\kappa_{\nu_{k,1}}(\uz_k)=\uk_{\nu_{k,1}}$ and $\kappa_{\nu_{k,2}}(\uz_k)=\uk_{\nu_{k,1}}$, for all $j\in \Z$, the function $z\in B_\varepsilon(\uz_k)\cap \Oc \mapsto G_j(z)$  can be holomorphically extended on $B_\varepsilon(\uz_k)$ and
		\begin{align}
			\forall z\in B_\varepsilon(\uz_k), \forall j \geq 1, \quad |G_j(z)- f_{\nu_{k,1}}(z) \kappa_{\nu_{k,1}}(z)^j|\leq C\exp(-cj).\label{est_spa_B_cas3}\\
			\forall z\in B_\varepsilon(\uz_k), \forall j \leq0, \quad |G_j(z)- f_{\nu_{k,2}}(z) \kappa_{\nu_{k,2}}(z)^j|\leq C\exp(-c|j|).\label{est_spa_B_2_cas3}
		\end{align}
		Furthermore, knowing that $\uz_k=\uz_{\nu_{k,1}}=\uz_{\nu_{k,2}}$, we have that 
		\begin{equation}\label{fk_2}
			\forall z\in B_\varepsilon(\uz_k),\quad  f_{\nu_{k,1}}(z)=-\frac{\kappa_{\nu_{k,1}}^\prime(z)}{\kappa_{\nu_{k,1}}(z)},\quad f_{\nu_{k,2}}(z)=\frac{\kappa_{\nu_{k,2}}^\prime(z)}{\kappa_{\nu_{k,2}}(z)}.
		\end{equation}
	\end{lemma}
	
	Lemmas \ref{green_spatial_près_1} and \ref{green_spatial_près_2} are similar to \cite[Lemmas 3 and 4]{Cou-Faye} but instead of proving sharp bounds on the spatial Green's function, we express its precise behavior near the points $\uz_k$. This is the crucial improvement with respect to \cite{Cou-Faye} that will allow us to find their asymptotic behavior and prove a sharp bound for the remainder. 
	
	\begin{proof}\textbf{of Lemma \ref{green_spatial_près_1}}
		Our proof will follow that of \cite[Lemmas 3, 4]{Cou-Faye}. First, we observe that case\textbf{ II} would be dealt similarly as case\textbf{ I} and that case\textbf{ III} is a mixture of both cases\textbf{ I} and \textbf{II}. Therefore, we will only detail the proof of Lemma \ref{green_spatial_près_1} in case\textbf{ I} and leave the proof of Lemma \ref{green_spatial_près_2} to the interested reader. We therefore consider $k\in\lc1,\ppp,K\rc$ so that we are in case\textbf{ I}. Lemma \ref{spec_spl} implies that $\uk_k$ is a simple eigenvalue of $\M(\uz_k)$. Thus, we can find a holomorphic function $\kappa_k$ defined on a neighborhood $B_\varepsilon(\uz_k)$ of $\uz_k$ such that for all $z\in B_\varepsilon(\uz_k)$, $\kappa_k(z)$ is an algebraically simple eigenvalue of $\M(z)$ and $\kappa_k(\uz_k)=\uk_k$. We also know that for all $z\in B_\varepsilon(\uz_k)$, the vector
		$$R_k(z):= \begin{pmatrix}
			\kappa_k(z)^{p+r-1} \\
			\vdots\\
			\kappa_k(z)\\
			1
		\end{pmatrix}\in\C^{p+r}$$
		is an eigenvector of $\M(z)$ associated with $\kappa_k(z)$. Because of Lemma \ref{spec_spl}, even if we have to take a smaller radius $\varepsilon$, we can assume that for all $z\in B_\varepsilon(\uz_k)$, $\M(z)$ has $\kappa_k(z)$ as a simple eigenvalue, $r-1$ eigenvalues different from $\kappa_k(z)$ in $\D$ and $p$ eigenvalues different from $\kappa_k(z)$ in $\Uc$. We define $E^s(z)$ (resp. $E^u(z)$) the strictly stable (resp. strictly unstable) subspace of $\M(z)$ which corresponds to the subspace spanned by the generalized eigenvectors of $\M(z)$ associated with eigenvalues different from $\kappa_k(z)$ in $\D$ (resp. $\Uc$). We therefore know that $E^s(z)$ (resp. $E^u(z)$) has dimension $r-1$ (resp. $p$) thanks to Lemma \ref{spec_spl} and we have the decomposition
		$$\C^{p+r} = E^s(z) \oplus E^u(z) \oplus \text{ Span } R_k(z).$$ 
		The associated projectors are denoted $\pi^s(z)$, $\pi^u(z)$ and $\pi^k(z)$. Those linear maps commute with $\M(z)$ and depend holomorphically on $z\in B_\varepsilon(\uz_k)$ (see \cite[I. Problem 5.9]{Kato}).
		
		For all $z\in B_\varepsilon(\uz_k)\cap\Oc$ and $j\in \Z$, $G_j(z)$ and the vector $W_j(z)$ are well defined. Also, by Lemma \ref{spec_spl}, we have that $|\kappa_k(z)|<1$ for all $z\in B_\varepsilon(\uz_k)\cap\Oc$. By reasoning in the same manner as in the proof of \cite[Lemma 3]{Cou-Faye}, we have for all $z\in B_\varepsilon(\uz_k)\cap\Oc$ and $j\in \Z$
		\begin{align}
			\pi^u(z)W_j(z) &=  -\frac{\ind_{j\in]-\infty,0]}}{\A_p(z)}\M(z)^{j-1}\pi^u(z)e, \label{wj_u_2}\\
			\pi^s(z)W_j(z) &= \frac{\ind_{j\in[1,+\infty [}}{\A_p(z)}\M(z)^{j-1}\pi^s(z)e, \label{wj_s_2}\\
			\pi^k(z)W_j(z) &= \frac{\ind_{j\in[1,+\infty [}}{\A_p(z)}\M(z)^{j-1}\pi^k(z)e = \frac{\ind_{j\in[1,+\infty [}}{\A_p(z)}\kappa_k(z)^{j-1}\pi^k(z)e. \label{wj_k}
		\end{align}
		
		We observe that the right hand side in the equations \eqref{wj_u_2}, \eqref{wj_s_2} and \eqref{wj_k} can be holomorphically extended on $B_\varepsilon(\uz_k)$. Therefore, we can extend holomorphically the applications which associates $z$ to $\pi^s(z)W_j(z)$, $\pi^u(z)W_j(z)$ and $\pi^k(z)W_j(z)$ on the whole open ball $B_\varepsilon(\uz_k)$ and this allows us to extend $W_j(z)$ on $B_\varepsilon(\uz_k)$. Since $G_j(z)$ is a coordinate of the vector $W_j(z)$, the holomorphic extension property is proved.
		
		By reasoning in the same manner as in the proof of the inequality \cite[(23)]{Cou-Faye}, we prove that there exist two constants $C,c>0$ such that 
		$$\forall z \in B_\varepsilon(\uz_k), \forall j \in \Z, \quad \|\pi^s(z)W_j(z) + \pi^u(z)W_j(z)\|\leq C\exp(-c|j|).$$
		This implies that 
		$$\forall z \in B_\varepsilon(\uz_k), \forall j \in \Z, \quad \|W_j(z) - \pi^k(z)W_j(z)\|\leq C\exp(-c|j|).$$
		
		This is now where our proof differs from the proof of \cite[Lemmas 3, 4]{Cou-Faye}. In \cite{Cou-Faye}, the authors find bounds on $\pi^k(z)W_j(z)$ and thus obtain estimates on $G_j(z)$. In our case, we have a stronger hypothesis (Hypothesis \ref{H1}) that allows us to have a much simpler expression \eqref{wj_k} of $\pi^k(z)W_j(z)$ and this will enable us to find the precise behavior of $G_j(z)$.
		
		For $j\leq 0$, we observe from \eqref{wj_k} that $\pi^k(z)W_j(z)=0$ and that $G_j(z)$ is a component of $W_j(z)$. We therefore get the inequality \eqref{est_spa_B_2_cas1}.
		
		We now consider the case $j\geq 1$. We have that $G_j(z)=(W_j(z))_p$ for all $z\in B_\varepsilon(\uz_k)$ where $(X)_p$ refers to the $p$-th coordinate of a vector $X\in\C^{p+r}$. Then, 
		$$\forall z \in B_\varepsilon(\uz_k), \quad |G_j(z) - (\pi^k(z)W_j(z))_p|\leq C\exp(-c|j|).$$
		
		We then define the holomorphic function 
		$$\begin{array}{cccc}
			f_k: & B_\varepsilon(\uz_k) & \rightarrow & \C \\
			& z & \mapsto & \frac{1}{\A_p(z)\kappa_k(z)}(\pi^k(z)e)_p
		\end{array}.$$
		By observing that $(\pi^k(z)W_j(z))_p = f_k(z)\kappa_k(z)^j$, we get the inequality \eqref{est_spa_B_cas1} and it now remains to obtain the expression \eqref{fk_1}. We first need to determine the spectral projector $\pi^k(z)$. We recall that $\kappa_k(z)\in\S^1$ is a simple eigenvalue of $\M(z)$ and the vector 
		$$R_k(z)= \begin{pmatrix}
			\kappa_k(z)^{p+r-1} \\
			\vdots\\
			\kappa_k(z)\\
			1
		\end{pmatrix}\in\C^{p+r}$$
		is an eigenvector of $\M(z)$ associated with $\kappa_k(z)$. We also know that there exists a unique eigenvector $L_k(z)=(l_j(z))_{j\in\lc1,\ppp,p+r\rc}\in\C^{p+r}$ of $\M(z)^T$ associated with the eigenvalue $\kappa_k(z)$ such that
		$$L_k(z)\cdot R_k(z)=1$$
		where the symmetric bilinear form $\cdot$ on $\C^{p+r}$ is defined as\footnote{Observe that this symetric bilinear form is not the Hermitian product on $\C^{p+r}$.}
		$$\forall X,Y\in \C^{p+r},\quad X\cdot Y := \sum_{l=1}^{p+r} X_iY_i.$$
		Then, we have that 
		$$\forall Y\in \C^{p+r},\quad \pi^k(z)Y = (L_k(z)\cdot Y)R_k(z).$$
		Thus, applying to the vector $e$ implies that
		\begin{equation}\label{fInt}
			 f_k(z)= \frac{l_1(z)\kappa_k(z)^{r-1}}{\A_p(z)}.
		\end{equation}
		
		We thus need to find the value of the coefficient $l_1(z)$. Since $L_k(z)$ is an eigenvalue of $\M(z)^T$ for the eigenvalue $\kappa_k(z)$, we get
		$$\forall j \in\lc1,\ppp,p+r\rc,\quad l_j(z) = -\left(\sum_{l=-r}^{p-j} \frac{\A_l(z)}{\kappa_k(z)^{p-j+1-l}}\right)\frac{l_1(z)}{\A_p(z)}.$$
		We now have an expression of each $l_j(z)$ depending on $l_1(z)$. To determine the value of $l_1(z)$, we have to use the normalization condition that we have made between $L_k(z)$ and $R_k(z)$. We have
		$$1= L_k(z)\cdot R_k(z) = \sum_{j=1}^{p+r} \kappa_k(z)^{p+r-j}l_j(z) = -\left(\sum_{j=1}^{p+r} \sum_{l=-r}^{p-j} \A_l(z)\kappa_k(z)^{l+r-1}\right)\frac{l_1(z)}{\A_p(z)}.$$
		By the expression of $\A_l(z)$, this implies that
		\begin{align*}
			1=  -\left( \sum_{l=-r}^{p} (p-l)\A_l(z)\kappa_k(z)^{l+r-1}\right)\frac{l_1(z)}{\A_p(z)} &= -\left(p\kappa_k(z)^{r-1}z - \sum_{l=-r}^{p} (p-l)a_l\kappa_k(z)^{l+r-1}\right)\frac{l_1(z)}{\A_p(z)}\\
			&=-\left(p\kappa_k(z)^{r-1}(z - F(\kappa_k(z))) + \kappa_k(z)^rF^\prime(\kappa_k(z))\right)\frac{l_1(z)}{\A_p(z)}.
		\end{align*}
		Since $\kappa_k(z)$ is an eigenvalue of $\M(z)$, Lemma \ref{spec_spl} implies that 
		$$F(\kappa_k(z))=z\quad \text{and} \quad \kappa_k^\prime(z)F^\prime(\kappa_k(z))=1.$$
		Thus, 
		$$1=-\frac{\kappa_k(z)^rl_1(z)}{\kappa_k^\prime(z)\A_p(z)}.$$
		Combining this equality with \eqref{fInt} implies the equality \eqref{fk_1}.
	\end{proof}

	\section{Temporal Green's function}\label{sec_GT}
	
	We are now ready to start proving Theorem \ref{thPrinc}. In Section \ref{subsec_est_far}, we will prove the result of the theorem far from the axes $j=n\alpha_k$. In this regime, the estimates proved in \cite[Theorem 1]{Cou-Faye} on $\Gcc_j^n$ and estimates on the derivatives of the function $H_{2\mu}^\beta$ will allow us to prove bounds that are even stronger than those claimed in Theorem \ref{thPrinc}. The bulk of the proof will happen in the case where $j-n\alpha_k$ is close to $0$ as the limiting estimates of Theorem \ref{thPrinc} occur in this case. Section \ref{subsecPlan} will summarize the idea of the proof in the case where $j$ is close to $n\alpha_k$ and Sections \ref{sec_GS/GT}-\ref{subsecLemImpTh} give the details. The main tools are the use of functional calculus (see \cite[Chapter VII]{Conway}) to express the elements $\Gcc_j^n$ with the spatial Green's function $G_j(z)$ and the estimates on the spatial Green's function proved in Section \ref{sec_GS}.
	
	Before we start, we are going to make two hypotheses to simplify the proof. The first one is that $-1\notin\lc\uz_1,\ppp,\uz_K\rc$. This hypothesis is actually not restrictive. If it were not verified, we would just have to multiply the sequence $a$ by some well chosen element of $\S^1$ to find a new sequence $b$ that will verify this hypothesis and prove the theorem for this new sequence. The theorem for our previous sequence $a$ would directly follow.
	
	The second hypothesis we make is that all $\alpha_k$ are distinct from one another. This hypothesis has a real impact on the proof, symplifying greatly some parts of the calculations. We will come back in Section \ref{sec_alpha_eg} to the case where the elements $\alpha_k$ can be equal and explain which elements of the proof should be modified.
	
	\subsection{Estimates far from the axes $j=n\alpha_k$}\label{subsec_est_far}
	
	As explained at the beginning of the section, we suppose that all $\alpha_k$ are distinct from one another. Without loss of generality, we suppose that we arranged them so that there holds:
	$$\alpha_1<\ppp<\alpha_k<\ppp<\alpha_K.$$
	
	For all $k\in\lc1,\ppp,K\rc$, we define two elements $\ud_k,\sd_k\in\R^*$ such that $\ud_k,\sd_k$ and $\alpha_k$ have the same sign and
	$$\ud_1<\alpha_1<\sd_1<\ppp<\ud_k<\alpha_k<\sd_k<\ppp<\ud_K<\alpha_K<\sd_K.$$
	We now define for every $k\in\lc1,\ppp,K\rc$ the sector
	$$\Dc_k:= \lc(n,j)\in \N^*\times \Z, \quad n\ud_k\leq j \leq n\sd_k\rc$$
	that do not intersect each other. We also introduce
	$$\Dc:=\bigcup_{k=1}^K \Dc_k.$$

	\begin{figure}
		\begin{center}
			\begin{tikzpicture}[scale=4]
				
				\draw (-0.5,0.5) node[fill=white] {$\Dc_1$};
				\draw (0.71,0.75) node[fill=white] {$\Dc_2$};
				\draw (0.5,0.3) node[fill=white] {$\Dc_3$};
				
				\fill[color=gray!20] (0,0) -- (-1,0.7) -- (-1,0.3) -- cycle;
				
				\fill[color=gray!20] (0,0) -- (1/3,1) -- (0.7,1) -- cycle;

				\fill[color=gray!20] (0,0) -- (1,0.4) -- (1,0.1) -- cycle;
				
				\draw[color=gray!60] (-1,0) grid[step=0.2] (1,1);
				\draw[->] (-1,0) -- (1,0) node[right] {$j$};
				\draw[->] (0,-0.1) -- (0,1) node[above right] {$n$};
				
				\draw[thick, red] (0,0) -- (-1,0.3) node[left] {$\ud_1$};
				\draw[thick, blue] (0,0) -- (-1,0.5) node[left] {$\alpha_1$};
				\draw[thick, red] (0,0) -- (-1,0.7) node[left] {$\sd_1$};
				
				\draw[thick, red] (0,0) -- (1/3,1) node[above] {$\ud_2$};
				\draw[thick, blue] (0,0) -- (1/2,1) node[above] {$\alpha_2$};
				\draw[thick, red] (0,0) -- (0.7,1) node[above] {$\sd_2$};
				
				\draw[thick, red] (0,0) -- (1,0.4) node[right] {$\ud_3$};
				\draw[thick, blue] (0,0) -- (1,0.25) node[right] {$\alpha_3$};
				\draw[thick, red] (0,0) -- (1,0.1) node[right] {$\sd_3$};
			\end{tikzpicture}
			\caption{An illustration of the sectors $\Dc_k$. Here, we have $\alpha_1=-2$, $\alpha_2=0.5$ and $\alpha_3=4$. The rays labeled $\alpha_k$ (resp. $\ud_k$, $\sd_k$) correspond to the ray $j=n\alpha_k$ (resp. $j=n\ud_k$, $j=n\sd_k$). We observe that, because $\ud_k$, $\alpha_k$ and $\sd_k$ have the same sign, $j$ and $\alpha_k$ have the same sign for $(n,j)\in\Dc_k$. Also, the sectors $\Dc_k$ do not intersect each other.}
			\label{sec}
		\end{center}
	\end{figure}
	
	We represent the sectors $\Dc_k$ on the Figure \ref{sec}. In this section, we are going to prove the following two lemmas, which give estimates on the Green's function $\Gcc_j^n$ and on the elements in its asymptotic expansion \eqref{devPrinc} outside of the sectors $\Dc_k$.
	
	\begin{lemma}\label{prop_loin_G}
		We have that
		$$\forall (n,j)\in\N\times\Z, \quad j<-np \text{ or } j>nr \Rightarrow \Gcc_j^n=0.$$
		Furthermore, there exist two constants $C,c>0$ such that
		\begin{equation}\label{inLoinG}
			\forall (n,j)\in(\N^*\times\Z)\backslash \Dc, \quad -np\leq j\leq nr \Rightarrow \left|\Gcc_j^n\right|\leq C\exp(-c(n+|j|)).
		\end{equation}
	\end{lemma}
	
	\begin{lemma}\label{prop_loin_H}
		We consider $k\in\lc1,\ppp,K\rc$ and $\Pcc\in\C[X,Y]$. For all $s\in \N$, there exist two constants $C,c>0$ such that
		\begin{equation}\label{inLoinH}
			\forall (n,j)\in(\N^*\times\Z)\backslash \Dc_k, \quad \left|\left(\Pcc\left(X_{n,j,k},\frac{d}{dx}\right)H_{2\mu_k}^{\beta_k}\right)\left(X_{n,j,k}\right) \right|\leq  \frac{C}{n^\frac{s+1}{2\mu_k}}\exp\left(-c\left|X_{n,j,k}\right|^\frac{2\mu_k}{2\mu_k-1}\right)
		\end{equation}
		where $X_{n,j,k}:=\frac{n\alpha_k-j}{n^\frac{1}{2\mu_k}}$.
	\end{lemma}
	
	Both lemmas are proved in a similar way.

	\begin{proof}\textbf{of Lemma \ref{prop_loin_G}}
		The first part of Lemma \ref{prop_loin_G} is directly proved recursively using the definition \eqref{defGreenTempo} of the elements $\Gcc_j^n$ and the equality \eqref{expLa} on the operator $\Lcc_a$. We now focus our attention on the inequality \eqref{inLoinG} of Lemma \ref{prop_loin_G}. The result \cite[Theorem 1]{Cou-Faye} gives us the existence of two constants $C,c>0$ such that
		$$\forall n \in\N^*, \forall j \in \Z, \quad |\Gcc_j^n|\leq \sum_{k=1}^K\frac{C}{n^\frac{1}{2\mu_k}} \exp\left(-c\left(\frac{|j-n\alpha_k|}{n^\frac{1}{2\mu_k}}\right)^\frac{2\mu_k}{2\mu_k-1}\right).$$
		
		For a sufficiently small $\tilde{c}>0$, we have that
		\begin{equation}\label{inegLemPropLoinG}
			\forall k\in\lc1,\ppp,K\rc, \forall (n,j)\in(\N^*\times \Z)\backslash \Dc, \quad -np\leq j \leq nr \Rightarrow  c\left(\frac{|j-n\alpha_k|}{n^\frac{1}{2\mu_k}}\right)^\frac{2\mu_k}{2\mu_k-1}\geq \tilde{c}(n+|j|).
		\end{equation}
		
		Therefore, we prove that there exist two positive constants $C,c$ such that
		$$\forall (n,j)\in(\N^*\times\Z)\backslash \Dc, \quad -np\leq j \leq nr \Rightarrow \left|\Gcc_j^n\right|\leq C\exp(-c(n+|j|)).$$
	\end{proof}
	
	To prove Lemma \ref{prop_loin_H}, we use the following lemma which gives sharp estimates on the derivatives of the function $H_{2\mu}^\beta$.
	\begin{lemma}\label{ineg_H}
		For $\mu\in\N^*$, $\beta\in \C$ with positive real part and $m\in \N$, there exist two constants $C,c>0$ such that
		$$\forall x\in \R, \quad  \left|{H_{2\mu}^\beta}^{(m)}(x)\right|\leq C\exp\left(-c|x|^\frac{2\mu}{2\mu-1}\right).$$
	\end{lemma}
	This lemma is proved in \cite[Proposition 5.3]{Rob}. For the sake of completeness, we give a complete proof in the appendix (Section \ref{sec_appendix}). 
	
	\begin{proof}\textbf{of Lemma \ref{prop_loin_H}}
		We fix a $k\in\lc1,\ppp,K\rc$ and we verify the estimate of Lemma \ref{prop_loin_H} for the monomial $\Pcc=X^{l_X}Y^{l_Y}$ where $l_X,l_Y\in\N$. We use Lemma \ref{ineg_H} which implies the existence of two constants $C,c>0$ such that 
		$$\forall (n,j)\in\N^*\times \Z,\quad \left|\left(\frac{n\alpha_k-j}{n^\frac{1}{2\mu_k}}\right)^{l_X}\left(H_{2\mu_k}^{\beta_k}\right)^{(l_Y)}\left(\frac{n\alpha_k-j}{n^\frac{1}{2\mu_k}}\right)\right|\leq C\left(\frac{|n\alpha_k-j|}{n^\frac{1}{2\mu_k}}\right)^{l_X}\exp\left(-c\left(\frac{|n\alpha_k-j|}{n^\frac{1}{2\mu_k}}\right)^\frac{2\mu_k}{2\mu_k-1}\right).$$
		This implies that there exists $\tilde{C}>0$ such that
		$$\forall (n,j)\in\N^*\times \Z,\quad \left|\left(\frac{n\alpha_k-j}{n^\frac{1}{2\mu_k}}\right)^{l_X}\left(H_{2\mu_k}^{\beta_k}\right)^{(l_Y)}\left(\frac{n\alpha_k-j}{n^\frac{1}{2\mu_k}}\right)\right|\leq \tilde{C}\exp\left(-\frac{c}{2}\left(\frac{|n\alpha_k-j|}{n^\frac{1}{2\mu_k}}\right)^\frac{2\mu_k}{2\mu_k-1}\right).$$
		Using the definition of the set $\Dc_k$, we prove the existence of a constant $\tilde{c}>0$ such that
		$$\forall (n,j)\in (\N^*\times\Z)\backslash\Dc_k, \quad \frac{c}{4}\left(\frac{|n\alpha_k-j|}{n^\frac{1}{2\mu_k}}\right)^\frac{2\mu_k}{2\mu_k-1}\geq \tilde{c} n.$$
		Therefore, we easily conclude that there exist two positive constants $C,c$ such that the inequality \eqref{inLoinH} of Lemma \ref{prop_loin_H} is verified for $\Pcc = X^{l_X}Y^{l_Y}$.
	\end{proof}
	
	Now that the two Lemmas \ref{prop_loin_G} and \ref{prop_loin_H} are proved, we observe that for any family of polynomials $(\Pcc^k_\sigma)_{k\in\lc1, \ppp,K\rc, \sigma\in\N^*}$ which belong to $\C[X,Y]$, for all $s_1, \ppp,s_K\in\N$, there exist two positive constants $C,c$ such that for all $(n,j)\in\N^*\times\Z\backslash\Dc$
	\begin{equation}\label{inHorsD}
		\left|\Gcc_j^n -\sum_{k=1}^K\sum_{\sigma=1}^{s_k}\frac{\uz_k^n\uk_k^j}{n^{\frac{\sigma}{2\mu_k}}} \left(\Pcc^k_\sigma\left(X_{n,j,k},\frac{d}{dx}\right)H_{2\mu_k}^{\beta_k}\right)\left(X_{n,j,k}\right)\right|  \leq \sum_{k=1}^K\frac{C}{n^\frac{s_k+1}{2\mu_k}}\exp\left(-c\left|X_{n,j,k}\right|^\frac{2\mu_k}{2\mu_k-1}\right),
	\end{equation}
	and for any $k_0\in\lc1,\ppp,K\rc$ and for all $(n,j)\in\Dc_{k_0}$, since the sets $\Dc_k$ do not intersect each other
	\begin{multline}\label{inHorsDk}
		\left|\Gcc_j^n -\sum_{k=1}^K\sum_{\sigma=1}^{s_k}\frac{\uz_k^n\uk_k^j}{n^{\frac{\sigma}{2\mu_k}}} \left(\Pcc^k_\sigma\left(X_{n,j,k},\frac{d}{dx}\right)H_{2\mu_k}^{\beta_k}\right)\left(X_{n,j,k}\right)\right|  \leq \sum_{\underset{k\neq k_0}{k=1}}^K\frac{C}{n^\frac{s_k+1}{2\mu_k}}\exp\left(-c\left|X_{n,j,k}\right|^\frac{2\mu_k}{2\mu_k-1}\right) \\ +\left|\Gcc_j^n -\sum_{\sigma=1}^{s_{k_0}}\frac{\uz_{k_0}^n\uk_{k_0}^j}{n^{\frac{\sigma}{2\mu_{k_0}}}} \left(\Pcc^{k_0}_{\sigma}\left(X_{n,j,k_0},\frac{d}{dx}\right)H_{2\mu_{k_0}}^{\beta_{k_0}}\right)\left(X_{n,j,k_0}\right)\right|
	\end{multline}	
	with $X_{n,j,k}=\frac{n\alpha_k-j}{n^\frac{1}{2\mu_k}}$. There just remains to find a family of polynomials $(\Pcc^k_\sigma)_{k,\sigma}$ to bound the last term in \eqref{inHorsDk} when $(n,j)\in\Dc_{k_0}$.
	
	\subsection{Plan of the proof of Theorem \ref{thPrinc} close to the axes $j=n\alpha_k$}\label{subsecPlan}
	
	We claim that to conclude the proof of Theorem \ref{thPrinc}, there only remains to prove the following lemma:
	\begin{lemma}\label{lemPrinc}
		For all $k\in\lc1,\ppp,K\rc$ and $s_k\in\N$, there exist a family of polynomials $(\Pcc^k_\sigma)_{\sigma\in\lc1,\ppp,s_k\rc}$ in $\C[X,Y]$ and two positive constants $C,c$ such that
		\begin{equation}\label{inPrinc}
				\forall (n,j)\in\Dc_k,  \quad  \left|\Gcc_j^n -\sum_{\sigma=1}^{s_k}\frac{\uz_k^n\uk_k^j}{n^{\frac{\sigma}{2\mu_k}}} \left(\Pcc^k_\sigma\left(X_{n,j,k},\frac{d}{dx}\right)H_{2\mu_k}^{\beta_k}\right)\left(X_{n,j,k}\right)\right|\leq \frac{C}{n^\frac{s_k+1}{2\mu_k}}\exp\left(-c|X_{n,j,k}|^\frac{2\mu_k}{2\mu_k-1}\right)
		\end{equation}
		with $X_{n,j,k}=\frac{n\alpha_k-j}{n^\frac{1}{2\mu_k}}$.
	\end{lemma}
	
	Once the existence of families of polynomials $(\Pcc^k_\sigma)_{k,\sigma}$ satisfying Lemma \ref{lemPrinc} is proved, the inequalities \eqref{inHorsD} and \eqref{inHorsDk} we deduced from Lemmas \ref{prop_loin_G} and \ref{prop_loin_H} imply that Theorem \ref{thPrinc} is also verified for the same family of polynomials. It is important to observe that we use intensively the fact that the sectors $\Dc_k$ do not intersect each other. In Section \ref{sec_alpha_eg}, we will see that when the elements $\alpha_k$ are not supposed to be different, we will need to adapt Lemma \ref{lemPrinc} to take into account that for each sector there could be multiple generalized Gaussian waves that are superposed in the estimate \eqref{devPrinc}.
	
	 We now focus our attention on proving Lemma \ref{lemPrinc}. We fix $k\in\lc1, \ppp,K\rc$ and $s\in \N$. For $s=0$, the result has been proved in \cite[Theorem 1]{Cou-Faye}. Therefore, we will focus on the case where $s\geq 1$. The proof of Lemma \ref{lemPrinc} in this case will be separated in three steps:
	
	$\bullet$ \textbf{Step 1:} In Section \ref{sec_GS/GT}, we will express the elements $\Gcc_j^n$ using the spatial Green's function $G_j(z)$ via the inverse Laplace transform and use the results of Section \ref{sec_GS} to prove the following lemma:
	
	\begin{lemma}\label{lemEstInterm}
		For all $k\in\lc1,\ppp,K\rc$ and for all $s\in \N^*$, there exist two positive constants $C,c$ such that for all $(n,j)\in \Dc_k$
		\begin{multline*}
			 \left|\Gcc_j^n -\frac{\uz_k^n\uk_k^j}{2\pi} \int_{-\infty}^{+\infty}P_{s,k}(it+\ut_k)\left(\sum_{l=0}^{s-1}\frac{(j R_{s,k}(it+\ut_k))^l}{l!}\right)\exp\left(it\left(n-\frac{j}{\alpha_k}\right) -\frac{j}{\alpha_k}\frac{\beta_k}{\alpha_k^{2\mu_k}}t^{2\mu_k}\right) dt\right|\\ \leq \frac{C}{n^\frac{s+1}{2\mu_k}}\exp\left(-c\left(\frac{|n\alpha_k-j|}{n^\frac{1}{2\mu_k}}\right)^\frac{2\mu_k}{2\mu_k-1}\right)
		\end{multline*} 
		where $\ut_k:=i\theta_k$ is the only element of $i ]-\pi,\pi[$ such that
		$$\uz_k =\exp(\ut_k) =\exp(i\theta_k)$$
		and the polynomial functions $P_{s,k}$ and $R_{s,k}$ have explicit expressions defined in Lemma \ref{lem_varpi}.
	\end{lemma}
	
	$\bullet$ \textbf{Step 2:} We observe that in Lemma \ref{lemEstInterm}, we approach the elements $\Gcc_j^n$ for $(n,j)\in \Dc_k$ by an explicit linear combination of the following terms where $l,m\in \N$ and $m\geq (2\mu_k+1)l$
	\begin{equation}\label{egHjalp}
		\frac{j^l}{2\pi}\int_{-\infty}^{+\infty}(it)^m\exp\left(it\left(n-\frac{j}{\alpha_k}\right) -\frac{j}{\alpha_k}\frac{\beta_k}{\alpha_k^{2\mu_k}}t^{2\mu_k}\right) dt= \frac{\alpha_k^{m+l}|\alpha_k|}{\left(\frac{j}{\alpha_k}\right)^\frac{m-2\mu_kl+1}{2\mu_k}}{H_{2\mu_k}^{\beta_k}}^{(m)}\left(\frac{n\alpha_k-j}{\left(\frac{j}{\alpha_k}\right)^\frac{1}{2\mu_k}}\right).
	\end{equation}
	
	If we compare the terms in \eqref{egHjalp} with the terms appearing in the estimate \eqref{devPrinc} of Theorem \ref{thPrinc}, since we are considering $(n,j)\in \Dc_k$, we see that $\frac{j}{\alpha_k}$ is close to $n$. Therefore, once Lemma \ref{lemEstInterm} is proved, we will only need some standard analysis in Section \ref{secAsymp} to prove the following lemma.
	
	\begin{lemma}\label{lemDevAsymp}
		For all $s \in \N$, $m\in \N$, $l\in\N\backslash\lc0\rc$ and $k\in \lc 1, \ppp,K\rc$, if we consider $d\in \N$ such that
		$$d\geq \frac{s+1}{2\mu_k-1}$$
		then there exist two constants $C,c>0$ such that for all $(n,j)\in \Dc_k$,
		$$\left|\frac{{H_{2\mu_k}^{\beta_k}}^{(m)}\left(Y_{n,j,k}\right)}{\left(\frac{j}{\alpha_k}\right)^\frac{l}{2\mu_k}}- \sum_{k_1=0}^{d-1}\sum_{k_3=0}^{d-1} \frac{\Bcc^k_{l,k_1,k_3}}{n^\frac{l+(2\mu_k-1)k_3}{2\mu_k}}\left(X_{n,j,k}\right)^{k_1+k_3}{H_{2\mu_k}^{\beta_k}}^{(m+k_1)}\left(X_{n,j,k}\right) \right|  \leq \frac{C}{n^\frac{s+1}{2\mu_k}}\exp\left(-c\left|X_{n,j,k}\right|^\frac{2\mu_k}{2\mu_k-1}\right)$$
		where $Y_{n,j,k}:= \frac{n\alpha_k-j}{\left(\frac{j}{\alpha_k}\right)^\frac{1}{2\mu_k}}$, $X_{n,j,k}:=\frac{n\alpha_k-j}{n^\frac{1}{2\mu_k}}$ and 
		$$\Bcc^k_{l,k_1,k_3}:= \sum_{k_2=0}^{k_1}\frac{\binom{k_1}{k_2}(-1)^{k_1-k_2}}{k_1!k_3!\alpha_k^{k_3}}\left(\prod_{k_4=0}^{k_3-1} \frac{l+k_2}{2\mu_k}+k_4\right)  .$$
	\end{lemma}
	
	$\bullet$ \textbf{Step 3:} In Section \ref{subsecLemImpTh}, we will explicitly construct the polynomials $\Pcc^k_\sigma$ satisfying Lemma \ref{lemPrinc} using Lemmas \ref{lemEstInterm} and \ref{lemDevAsymp}. This will conclude the proof of Lemma \ref{lemPrinc} and Theorem \ref{thPrinc} in the case where the elements $\alpha_k$ are distinct.
	
	\subsection{Step 1: Link between the spatial and temporal Green's functions and proof of Lemma \ref{lemEstInterm}}\label{sec_GS/GT}
	
	As explained at the end of the previous section, we start by proving Lemma \ref{lemEstInterm}. The first step will be to express the elements $\Gcc_j^n$ via the spatial Green's function $G_j(z)$. The equation \eqref{def_G} implies by using the inverse Laplace transform that if we define a path which surrounds $\sigma(\Lcc_a)=F(\S^1)$, like for example $\widetilde{\Gamma}_\rho=\exp(\rho)\S^1$ for $0<\rho\leq\pi$, then 
	$$\forall n\in \N^*, \forall j\in \Z, \quad \Gcc_j^n=\frac{1}{2i\pi} \int_{\widetilde{\Gamma}_\rho} z^nG_j(z)dz.$$
	We fix this choice of path for now but we are going to modify it in what follows. The idea will be to deform the path on which we integrate so that we can best use the estimates on $G_j(z)$ proved in Section \ref{sec_GS}. We start with a change of variable $z=\exp(\tau)$ in the previous equality. Therefore, if we define $\Gamma_\rho:= \lc \rho+il, l\in[-\pi,\pi]\rc$ and $\Gg_j(\tau)= e^\tau G_j(e^\tau)$, then 
	\begin{equation}\label{exp_Gcc}
		\forall n\in\N^*, \forall j\in\Z, \quad \Gcc_j^n=\frac{1}{2i\pi} \int_{\Gamma_\rho} e^{n\tau}\Gg_j(\tau)d\tau.
	\end{equation}
	We will therefore need a lemma that allows us to get from estimates on $G_j(z)$ to estimates on $\Gg_j(\tau)$. First, recalling that $\uz_k\neq -1$, we define for all $k\in\lc1,\ppp,K\rc$ the unique element $\ut_k:=i\theta_k$ of $i ]-\pi,\pi[$ such that
	$$\uz_k =\exp(\ut_k) =\exp(i\theta_k) .$$
	We also introduce for all $k\in\lc1,\ppp,K\rc$ the unique $\tilde{\theta}_k\in]-\pi,\pi]$ such that
	$$\uk_k=e^{i\tilde{\theta}_k}.$$
	
	We now introduce a lemma to pass from estimates on $G_j(z)$ to estimates on $\Gg_j(\tau)$.

	\begin{lemma}\label{lem_varpi}
		There exist a radius $\varepsilon_\star>0$ and for all $k\in \lc1,\ppp,K\rc$ two holomorphic functions $\varpi_k:B_{\varepsilon_\star}(\ut_k)\rightarrow \C$ and $g_k:B_{\varepsilon_\star}(\ut_k)\rightarrow \C$ such that for all $\varepsilon\in]0,\varepsilon_\star[$, there exist a width $\eta_\varepsilon\in]0,\varepsilon[$ and two constants $C,c>0$ such that if we define
		$$U_\varepsilon:=\lc\tau\in \C, \Re(\tau)\in]-\eta_\varepsilon, \pi], \Im(\tau)\in [-\pi,\pi]\rc \quad \text{and} \quad \Omega_\varepsilon := U_\varepsilon \backslash\bigcup_{k=1}^KB_\varepsilon(\ut_k),$$
		then for all $j\in\Z$, the application $\tau\mapsto \Gg_j(\tau)$ can be holomorphically extended on $U_\varepsilon\cup\bigcup_{k=1}^KB_\varepsilon(\ut_k)$ and we have that
		\begin{equation}\label{varpi_loin}
			\forall \tau \in \Omega_\varepsilon, \forall j\in\Z, \quad |\Gg_j(\tau)|\leq Ce^{-c|j|}.
		\end{equation}
		Also, for all $k\in\lc1,\ppp,K\rc$, depending on the case, we have that 
		
		\textbf{Case I:}
		\begin{align}
			\forall \tau \in B_\varepsilon(\ut_k), \forall j\geq 1, \quad |\Gg_j(\tau)-e^\tau g_k(\tau)e^{j\varpi_k(\tau)}|\leq Ce^{-c|j|}, \label{varpi_près_pos_I}\\
			\forall \tau \in B_\varepsilon(\ut_k), \forall j\leq 0, \quad |\Gg_j(\tau)|\leq Ce^{-c|j|}, \label{varpi_près_neg_I}
		\end{align}
		
		\textbf{Case II:}
		\begin{align}
			\forall \tau \in B_\varepsilon(\ut_k), \forall j\geq 1, \quad |\Gg_j(\tau)|\leq Ce^{-c|j|}, \label{varpi_près_pos_II}\\
			\forall \tau \in B_\varepsilon(\ut_k), \forall j\leq 0, \quad |\Gg_j(\tau)-e^\tau g_k(\tau)e^{j\varpi_k(\tau)}|\leq Ce^{-c|j|}, \label{varpi_près_neg_II}
		\end{align}
		
		\textbf{Case III:}
		\begin{align}
			\forall \tau \in B_\varepsilon(\ut_k), \forall j\geq 1, \quad |\Gg_j(\tau)-e^\tau g_{\nu_{k,1}}(\tau)e^{j\varpi_{\nu_{k,1}}(\tau)}|\leq Ce^{-c|j|}, \label{varpi_près_pos_III}\\
			\forall \tau \in B_\varepsilon(\ut_k), \forall j\leq 0, \quad |\Gg_j(\tau)-e^\tau g_{\nu_{k,2}}(\tau)e^{j\varpi_{\nu_{k,2}}(\tau)}|\leq Ce^{-c|j|}, \label{varpi_près_neg_III}
		\end{align}
		where we have $\Ic_k=\lc\nu_{k,1}, \nu_{k,2}\rc$, $\alpha_{\nu_{k,1}}>0$ and $\alpha_{\nu_{k,2}}<0$.
		
		For all $k\in\lc1,\ppp,K\rc$, we have
		\begin{equation}\label{eqVarpi}
			\varpi_k(\tau)\underset{\tau\rightarrow \ut_k}=i\tilde{\theta}_k -\frac{(\tau-\ut_k)}{\alpha_k}+(-1)^{\mu_k+1}\frac{\beta_k}{\alpha_k^{2\mu_k+1}}(\tau-\ut_k)^{2\mu_k}+o(|\tau-\ut_k|^{2\mu_k}).
		\end{equation}
		and
		\begin{equation}\label{eqG}
			\forall \tau \in B_{\varepsilon_\star}(\ut_k), \quad e^\tau g_k(\tau)=-\mathrm{sgn}(\alpha_k)\varpi_k^\prime(\tau)
		\end{equation}
		
		For $s\in \N^*$, we define the functions 
		$$\begin{array}{cccc}
			P_{s,k}: & \tau\in \C & \mapsto & -\mathrm{sgn}(\alpha_k) \displaystyle\sum_{l=0}^{s-1}\frac{\varpi_k^{(l+1)}(\ut_k)}{l!}(\tau-\ut_k)^l, \\
			\varphi_k : &\tau\in \C & \mapsto & i\tilde{\theta}_k-\displaystyle\frac{(\tau-\ut_k)}{\alpha_k}+(-1)^{\mu_k+1}\frac{\beta_k}{\alpha_k^{2\mu_k+1}}(\tau-\ut_k)^{2\mu_k},\\
			Q_{s,k}: & \tau\in \C & \mapsto &  \displaystyle\sum_{l=0}^{2\mu_k+s-1}\frac{\varpi_k^{(l)}(\ut_k)}{l!}(\tau-\ut_k)^l, \\
			R_{s,k}: & \tau\in \C & \mapsto & Q_{s,k}(\tau)-\varphi_k(\tau).
		\end{array}$$
		The functions $P_{s,k}$, $Q_{s,k}$ and $\varphi_k$ are asymptotic expansions of the function $e^\tau g_k$ and $\varpi_k$  at $\ut_k$ up to different orders. We can then define a bounded holomorphic function $\xi_{s,k}:B_{\varepsilon_\star}(\ut_k)\mapsto \C$ such that
		$$\forall \tau \in B_{\varepsilon_\star}(\ut_k), \quad \varpi_k(\tau)=Q_{s,k}(\tau)+\xi_{s,k}(\tau)(\tau-\ut_k)^{2\mu_k+s}.$$ 
		We then can prove that there exist two positive constants $A_R, A_I$ such that for all $\tau \in B_{\varepsilon_\star}(\ut_k)$
		\begin{align}
			\alpha_k\Re(\varphi_k(\tau))&\leq  -\Re(\tau-\ut_k)+A_R\Re(\tau-\ut_k)^{2\mu_k}-A_I\Im(\tau-\ut_k)^{2\mu_k},\label{estVarphi}\\
			\alpha_k\Re(\varpi_k(\tau))+|\alpha_k||\xi_{s,k}(\tau)(\tau-\ut_k)^{2\mu_k+s}|&\leq  -\Re(\tau-\ut_k)+A_R\Re(\tau-\ut_k)^{2\mu_k}-A_I\Im(\tau-\ut_k)^{2\mu_k},\label{estVarpi}\\
			\alpha_k\Re(\varphi_k(\tau))+|\alpha_k||R_{s,k}(\tau)|&\leq  -\Re(\tau-\ut_k)+A_R\Re(\tau-\ut_k)^{2\mu_k}-A_I\Im(\tau-\ut_k)^{2\mu_k}.\label{estR}
		\end{align}
	\end{lemma}

	\begin{proof}
		Using the Lemmas \ref{green_spatial_près_1} and \ref{green_spatial_près_2} and writing $\kappa_k(z)=\exp(\omega_k(z))$ for $z$ near $\uz_k$ with $\omega_k(\uz_k)=i\tilde{\theta}_k$, we can define for a choice of $\varepsilon_\star$ small enough two holomorphic functions $\varpi_k$ and $g_k$ such that
		$$\forall \tau\in B_{\varepsilon_\star}(\ut_k),\quad \varpi_k(\tau)= \omega_k(e^\tau), g_k(\tau)= f_k(e^\tau).$$  
		Lemmas \ref{green_spatial_près_1} and \ref{green_spatial_près_2} directly imply the inequalities \eqref{varpi_près_pos_I}, \eqref{varpi_près_neg_I}, \eqref{varpi_près_pos_II}, \eqref{varpi_près_neg_II}, \eqref{varpi_près_pos_III} and \eqref{varpi_près_neg_III} on the open balls $B_{\varepsilon_\star}(\ut_k)$ and the fact that the functions $\tau\mapsto \Gg_j(\tau)$ are holomorphic on $B_{\varepsilon_\star}(\ut_k)$. We now consider $\varepsilon\in]0,\varepsilon_\star[$. The inequalities we just proved remain true on $B_\varepsilon(\ut_k)$.  Using a compactness argument and Lemma \ref{green_spatial_loin}, we also get the existence of $\eta_\varepsilon$ and the inequality \eqref{varpi_loin}.
		
		We observe that the asymptotic expansion \eqref{F} implies that 
		$$\tau-\ut_k\underset{\tau\rightarrow \ut_k}{=} -\alpha_k(\varpi_k(\tau)-i\tilde{\theta}_k) +(-1)^{\mu_k+1}\beta_k  (\varpi_k(\tau)-i\tilde{\theta}_k)^{2\mu_k} + O\left(\left|\varpi_k(\tau)-i\tilde{\theta}_k\right|^{2\mu_k+1}\right).$$
		We then deduce the equation \eqref{eqVarpi}.
		
		For $\tau \in B_{\varepsilon_\star}(\ut_k)$, the equations \eqref{fk_1} and \eqref{fk_2} imply the equality \eqref{eqG}.
	
	 There only remains to prove the existence of $A_R$ and $A_I$ to verify the inequalities \eqref{estVarphi} - \eqref{estR}.
		
	We are going to prove \eqref{estVarphi} first. Because of Young's inequality, we have that for $l\in\lc1,\ppp,2\mu_k-1\rc$, for all $\delta>0$, there exists $C_\delta>0$ such that for all $\tau \in \C$ 
	$$|\Re(\tau)|^l|\Im(\tau)|^{2\mu_k-l}\leq \delta \Im(\tau)^{2\mu_k} + C_\delta \Re(\tau)^{2\mu_k}.$$
	
	Furthermore, we have that 
	$$\alpha_k\Re(\varphi_k(\tau))  =-\Re(\tau-\ut_k) + (-1)^{\mu_k+1}\left(\frac{\Re(\beta_k)}{\alpha_k^{2\mu_k}}\Re((\tau-\ut_k)^{2\mu_k}) - \frac{\Im(\beta_k)}{\alpha_k^{2\mu_k}}\Im((\tau-\ut_k)^{2\mu_k})\right).$$
	
	Then, for $\delta>0$, there exists $C_\delta>0$ such that 
	$$\alpha_k\Re(\varphi_k(\tau))  \leq -\Re(\tau-\ut_k) +\Re(\tau-\ut_k)^{2\mu_k} \left(\frac{\Re(\beta_k)}{\alpha_k^{2\mu_k}}+C_\delta\right)+ \Im(\tau-\ut_k)^{2\mu_k} \left(-\frac{\Re(\beta_k)}{\alpha_k^{2\mu_k}}+\delta\right).$$
	
	Therefore, by taking $\delta$ small enough, we can end the proof of inequality \eqref{estVarphi}. The proof of inequality \eqref{estVarpi} is similar. We have for $\tau\in B_{\varepsilon_\star}(\ut_k)$
	\begin{multline*}
		\alpha_k\Re(\varpi_k(\tau)) +|\alpha_k||\xi_{s,k}(\tau)(\tau-\ut_k)^{2\mu_k+s}|
		\leq -\Re(\tau-\ut_k)+  |\alpha_k|\left (2|\xi_{s,k}(\tau)| |\tau-\ut_k|^{2\mu_k+s} +|R_{s,k}(\tau)|\right)\\+ (-1)^{\mu_k+1}\left(\frac{\Re(\beta_k)}{\alpha_k^{2\mu_k}}\Re((\tau-\ut_k)^{2\mu_k}) - \frac{\Im(\beta_k)}{\alpha_k^{2\mu_k}}\Im((\tau-\ut_k)^{2\mu_k})\right) .
	\end{multline*}
	
	We know there exists $c_1,c_2>0$ such that 
	$$\forall k\in\lc1,\ppp,K\rc, \forall \tau \in \C,\quad |\tau|^{2\mu_k}\leq c_1\Re(\tau)^{2\mu_k} + c_2 \Im(\tau)^{2\mu_k}.$$
	
	Since $\xi_{s,k}$ and $\frac{R_{s,k}}{X^{2\mu_k+1}}$ can be bounded by some constant $\tilde{C}>0$ on $B_{\varepsilon_\star}(\ut_k)$, using the same reasoning as previously gives us 
	\begin{multline*}
		\alpha_k\Re(\varpi_k(\tau)) +|\alpha_k||\xi_{s,k}(\tau)(\tau-\ut_k)^{2\mu_k+s}|
		\leq -\Re(\tau-\ut_k) + |\alpha_k| \tilde{C}\left(2 \varepsilon_\star^{s}+\varepsilon_\star\right) (c_1\Re(\tau-\ut_k)^{2\mu_k} + c_2 \Im(\tau-\ut_k)^{2\mu_k})  \\ +\Re(\tau-\ut_k)^{2\mu_k} \left(\frac{\Re(\beta_k)}{\alpha_k^{2\mu_k}}+C_\delta\right)+ \Im(\tau-\ut_k)^{2\mu_k} \left(-\frac{\Re(\beta_k)}{\alpha_k^{2\mu_k}}+\delta\right).
	\end{multline*}
	
	Taking $\delta$ and ${\varepsilon_\star}$ small enough allows us to prove \eqref{estVarpi}. We prove the inequality \eqref{estR} the same way.	
	\end{proof}
	
	\begin{remark}
		We observe that the constants in the inequalities \eqref{varpi_près_pos_I}, \eqref{varpi_près_neg_I}, \eqref{varpi_près_pos_II}, \eqref{varpi_près_neg_II}, \eqref{varpi_près_pos_III} and \eqref{varpi_près_neg_III} can (and will) be chosen uniformly with respect to $\varepsilon\in]0,\varepsilon_\star[$. However, it is not the case for the constants in inequality \eqref{varpi_loin}.
	\end{remark}
	
	\subsubsection{ Choice of integration paths for the proof of Lemma \ref{lemEstInterm}}
	
	From now on, we fix a $k\in\lc1,\ppp,K\rc$ and an integer $s\in\N\backslash\lc0\rc$ and our goal is to prove the claim of Lemma \ref{lemEstInterm} for this $k$ and $s$, i.e. we want to prove the existence of two positive constants $C,c$ such that for all $(n,j)\in\Dc_k$ we have
	\begin{multline}\label{inlemEstINterm}
		\left|\Gcc_j^n -\frac{\uz_k^n\uk_k^j}{2\pi} \int_{-\infty}^{+\infty}P_{s,k}(it+\ut_k)\left(\sum_{l=0}^{s-1}\frac{(j R_{s,k}(it+\ut_k))^l}{l!}\right)\exp\left(it\left(n-\frac{j}{\alpha_k}\right) -\frac{j}{\alpha_k}\frac{\beta_k}{\alpha_k^{2\mu_k}}t^{2\mu_k}\right) dt\right|\\ \leq \frac{C}{n^\frac{s+1}{2\mu_k}}\exp\left(-c\left(\frac{|n\alpha_k-j|}{n^\frac{1}{2\mu_k}}\right)^\frac{2\mu_k}{2\mu_k-1}\right).
	\end{multline}
	We will suppose that $\alpha_k>0$. The major consequence is that for $(n,j)\in\Dc_k$, we have $j\geq 1$. This implies that we will use the inequalities \eqref{varpi_près_pos_I}, \eqref{varpi_près_pos_II} and \eqref{varpi_près_pos_III}. The case where $\alpha_k<0$ would need some little modifications, in particular we will have that $j\leq 0$ for $(n,j)\in\Dc_k$ and we would rather use the inequalities \eqref{varpi_près_neg_I}, \eqref{varpi_près_neg_II} and \eqref{varpi_près_neg_III}.
	
	Before we begin with the proof, we will need to introduce some lemmas and define some elements. First, we can easily prove the following lemma which allows us to pass from bounds that are exponentially decaying in $n$ to the generalized Gaussian bounds expected in \eqref{inlemEstINterm}.
	
	\begin{lemma}\label{est_expn}
		We consider $C,c>0$. Then, for all $s\in \N^*$, there exist $\tilde{C},\tilde{c}>0$ such that 
		$$\forall (n,j)\in\Dc_k, \quad  C\exp(-cn)\leq \frac{\tilde{C}}{n^\frac{s+1}{2\mu_k}}\exp\left(-\tilde{c}\left|X_{n,j,k}\right|^\frac{2\mu_k}{2\mu_k-1}\right)$$
		with $X_{n,j,k}:=\frac{n\alpha_k-j}{n^\frac{1}{2\mu_k}}$.
	\end{lemma}	
	
	We now apply Lemma \ref{lem_varpi} and consider $\varepsilon\in]0,\varepsilon_\star[$ small enough so that
	$$\forall i,j\in\lc1,\ppp,K\rc, \quad \uz_i\neq \uz_j \Rightarrow B_\varepsilon(\ut_i)\cap B_\varepsilon(\ut_j)=\emptyset $$
	and 
	$$\forall l\in\lc1,\ppp,K\rc,\quad B_\varepsilon(\ut_l)\subset \lc\tau\in \C, \quad\Im(\tau)\in [-\pi,\pi]\rc.$$	
	This can be done because we supposed that $\uz_l\neq -1$ which implies $\ut_l\notin\lc-i\pi,i\pi\rc$ for all $l$. We also introduce some conditions on the values $\eta_\varepsilon$ we defined in Lemma \ref{lem_varpi} which will be useful later on in the proof, especially for Lemma \ref{lemDécal}. We define the function
	\begin{equation}\label{defR}
		\begin{array}{cccc}
			r_\varepsilon: & ]0,\varepsilon[& \rightarrow & \R \\
			& \eta & \mapsto & \sqrt{\varepsilon^2-\eta^2}
		\end{array}
	\end{equation}
	which serves to define the extremities of $-\eta+i\R\cap B_\varepsilon(\ut_k)$ for any $k\in\lc 1,\ppp,K\rc$. We impose that $\eta_\varepsilon$ is small enough so that 
	\begin{equation}\label{inEta}
		\eta_\varepsilon<\sqrt{\frac{3}{4}}\varepsilon.
	\end{equation}
	This condition implies that
	$$r_\varepsilon(\eta_\varepsilon)>\frac{\varepsilon}{2}.$$
	Finally, we also impose that
	\begin{equation}\label{inEta2}
		\forall k\in\lc 1, \ppp,K\rc, \quad \eta_\varepsilon+A_R\eta_\varepsilon^{2\mu_k} -A_I \left(\frac{\varepsilon}{2}\right)^{2\mu_k}<0.
	\end{equation}
	We now fix a constant $\eta\in]0,\eta_\varepsilon[$ which we will use to express the modified path on which we will integrate the right-hand term of equality \eqref{exp_Gcc}. 	 
	
	We will now follow a strategy developed in \cite{ZH}, which has also been used in \cite{Godillon}, \cite{Cou-Faye} and \cite{Cou-Faye2}, and introduce a family of parameterized curves. For $\tau_p\in \R$, we introduce
	$$\Psi_k(\tau_p)=\tau_p-A_R{\tau_p}^{2\mu_k}.$$
	The function $\Psi_k$ is continuous and strictly increasing on $\left]-\infty,\left(\frac{1}{2\mu_kA_R}\right)^\frac{1}{2\mu_k-1}\right[$. We choose $\varepsilon$ small enough so that it is strictly increasing on $]-\infty,\varepsilon]$. We can therefore introduce for $\tau_p\in[-\eta,\varepsilon]$
	$$\Gamma_{k,p} = \lc \tau\in \C, -\eta\leq \Re(\tau)\leq \tau_p, \quad \Re(\tau-\ut_k) - A_R \Re(\tau-\ut_k)^{2\mu_k} +  A_I \Im(\tau-\ut_k)^{2\mu_k}= \Psi_k(\tau_p)\rc.$$
	It is a symmetric curve with respect to the axis $\R+\ut_k=\R+i\theta_k$ which intersects this axis on the point $\tau_p+\ut_k$. If we introduce $\ell_{k,p}= \left(\frac{\Psi_k(\tau_p)-\Psi_k(-\eta)}{A_I}\right)^\frac{1}{2\mu_k}$, then $-\eta +i(\theta_k+\ell_{k,p})$ and $-\eta +i(\theta_k-\ell_{k,p})$ are the end points of $\Gamma_{k,p}$. We can also introduce a parametrization of this curve by defining $\gamma_{k,p}:[-\ell_{k,p}, \ell_{k,p}]\rightarrow \C$ such that 
	\begin{equation}\label{param}
		\forall \tau_p\in\left[-\eta,\varepsilon\right], \forall t\in[-\ell_{k,p},\ell_{k,p}],\quad \Im(\gamma_{k,p}(t))=t+\theta_k, \quad \Re(\gamma_{k,p}(t))=h_{k,p}(t):=\Psi_k^{-1}\left(\Psi_k(\tau_p)-A_It^{2\mu_k}\right).
	\end{equation}
	
	The above parametrization immediately yields that there exists a constant $M>0$ such that 
	\begin{equation}\label{hp}
		\forall \tau_p \in[-\eta,\varepsilon], \forall t \in[-\ell_{k,p},\ell_{k,p}], \quad |h_{k,p}^\prime(t)|\leq M.
	\end{equation}
	Also, there exists a constant $c_\star>0$ such that 
	\begin{equation}
		\forall \tau_p\in[-\eta,\varepsilon], 	\forall \tau \in \Gamma_{k,p}, \quad \Re(\tau-\ut_k)-\tau_p\leq -c_\star \Im(\tau-\ut_k)^{2\mu_k}. \label{ine_Re}
	\end{equation}
	
	We introduce those integration paths $\Gamma_{k,p}$ because they allow us to use optimally the inequalities \eqref{estVarphi}-\eqref{estR}. For example, if we seek to bound $e^{n\tau+j\varpi_k(\tau)}$ when $(n,j)\in \Dc_k$ and $\tau\in \Gamma_{k,p}$, it follows from the equality $\mathrm{sgn}(j)=\mathrm{sgn}(\alpha_k)$ and the inequalities \eqref{estVarpi} and \eqref{ine_Re} that
	\begin{align}
		\begin{split}
		n\Re(\tau-\ut_k)+j\Re(\varpi_k(\tau))& \leq n\Re(\tau-\ut_k)-\frac{j}{\alpha_k} \left(\Re(\tau-\ut_k)-A_R\Re(\tau-\ut_k)^{2\mu_k}+A_I\Im(\tau-\ut_k)^{2\mu_k}\right)\\
		& \leq -nc_\star \Im(\tau-\ut_k)^{2\mu_k}- \left(\frac{j}{\alpha_k}-n\right)\tau_p +\frac{j}{\alpha_k}A_R\tau_p^{2\mu_k}.
		\end{split}\label{estClas}
	\end{align}
	
	Such calculations will happen regularly in the following proof (see Lemmas \ref{ine_taup} and \ref{ine_res}). There remains to make an appropriate choice of $\tau_p$ depending on $n$ and $j$ that minimizes the right-hand side of the inequality \eqref{estClas} whilst the paths $\Gamma_{k,p}$ remain within the ball $B_{\varepsilon}(\ut_k)$. Even if we have to consider a smaller $\eta$, we can define a real number $0<\varepsilon_{k,0}<\varepsilon$ such that the curve $\Gamma_{k,p}$ associated to $\tau_p=\varepsilon_{k,0}$ intersects the axis $-\eta+i\R$ within $B_\varepsilon(\ut_k)$. Then, we let 
	$$\zeta_k=\frac{j-n\alpha_k}{2\mu_k n}, \quad \gamma_k=\frac{A_Rj}{n}, \quad \rho_k\left(\frac{\zeta_k}{\gamma_k}\right)=\mathrm{sgn}(\zeta_k)\left(\frac{|\zeta_k|}{\gamma_k}\right)^\frac{1}{2\mu_k-1}.$$
	The inequality \eqref{estClas} thus becomes
	\begin{equation}\label{estClas2}
		n\Re(\tau-\ut_k)+j\Re(\varpi_k(\tau))\leq -nc_\star \Im(\tau-\ut_k)^{2\mu_k}+\frac{n}{\alpha_k}(\gamma_k \tau_p^{2\mu_k}-2\mu_k\zeta_k\tau_p).
	\end{equation}
	Our limiting estimates will come from the case where $\zeta_k$ is close to $0$. We observe that the condition $(n,j)\in\Dc_k$ implies 
	\begin{equation}
		A_R\ud_k\leq \gamma_k \leq A_R\sd_k.\label{ineg_gamma}
	\end{equation}
	
	Moreover, we have that $\rho_k\left(\frac{\zeta_k}{\gamma_k}\right)$ is the unique real root of the polynomial
	$$\gamma_k x^{2\mu_k-1}=\zeta_k.$$
	
	Then, we take 
	$$\tau_p:=\lc\begin{array}{ccc}
		\rho_k\left(\frac{\zeta_k}{\gamma_k}\right), & \text{ if }\rho_k\left(\frac{\zeta_k}{\gamma_k}\right)\in[-\frac{\eta}{2},\varepsilon_{k,0}],& \text{(Case \textbf{A})}\\
		\varepsilon_{k,0}, & \text{ if }\rho_k\left(\frac{\zeta_k}{\gamma_k}\right)>\varepsilon_{k,0}, &\text{(Case \textbf{B})}\\
		-\frac{\eta}{2}, & \text{ if }\rho_k\left(\frac{\zeta_k}{\gamma_k}\right)<-\frac{\eta}{2}.&\text{(Case \textbf{C})}\\
	\end{array}\right.$$
	The case \textbf{A} corresponds to the choice to minimize the right-hand side of \eqref{estClas2}. The cases \textbf{B} and \textbf{C} allow the path $\Gamma_{k,p}$ to stay within $B_\varepsilon(\ut_k)$.
	
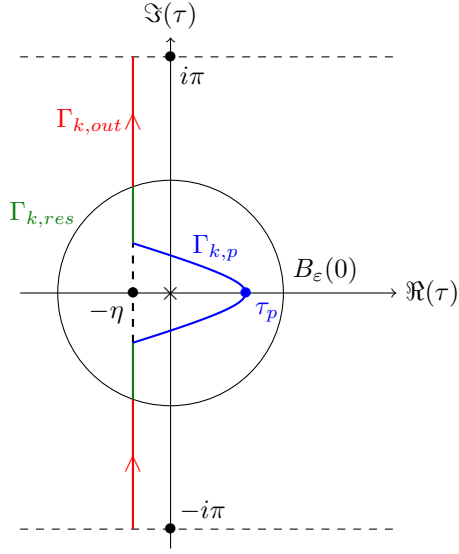
\begin{figure}
		\begin{center}
			\begin{tikzpicture}[scale=1]
				\draw[->] (-2,0) -- (3,0) node[right] {$\Re(\tau)$};
				\draw[->] (0,-3.4) -- (0,3.4) node[above] {$\Im(\tau)$};
				\draw[dashed] (-2,pi) -- (3,pi);
				\draw[dashed] (-2,-pi) -- (3,-pi);
				\draw (0,pi) node {$\bullet$} node[below right] {$i\pi$};
				\draw (0,-pi) node {$\bullet$} node[above right] {$-i\pi$};
				\draw (-0.5,0) node {$\bullet$} node[below left] {$-\eta$};
				\draw (0,0) node {$\times$} circle (1.5);
				\draw (1.5,0) node[above right] {$B_\varepsilon(0)$};
				\draw[blue] (1,0) node {$\bullet$} node[below right] {$\tau_p$};
				\draw[red,thick] (-0.5,-pi) -- (-0.5,{-sqrt(1.5^2-0.5^2)}) node[midway, sloped] {$>$};
				\draw[red,thick] (-0.5,{sqrt(1.5^2-0.5^2)}) -- (-0.5,pi) node[midway, left] {$\Gamma_{k,out}$} node[midway, sloped] {$>$};
				\draw[dartmouthgreen,thick] (-0.5,{-sqrt(1.5^2-0.5^2)}) -- (-0.5,{-sqrt((3/4+0.5+1/16)/3)});
				\draw[dartmouthgreen,thick] (-0.5,{sqrt((3/4+0.5+1/16)/3)}) -- (-0.5,{sqrt(1.5^2-0.5^2)}) ;
				\draw[dartmouthgreen] (-1.7,{(sqrt(1.5^2-0.5^2)+sqrt((3/4+0.5+1/16)/3))/2}) node {$\Gamma_{k,res}$};
				\draw[dashed,thick] (-0.5,{-sqrt((3/4+0.5+1/16)/3)}) -- (-0.5,{sqrt((3/4+0.5+1/16)/3)});
				\draw[thick,blue] plot [samples = 100, domain={-sqrt((3/4+0.5+1/16)/3)}:{sqrt((3/4+0.5+1/16)/3)}] ({2*(1-sqrt(1/4+3*abs(\x)^2))},\x) ;
				\draw[blue] (0.6,0.6) node {$\Gamma_{k,p}$};
			\end{tikzpicture}
			\caption{A representation of the path $\Gamma_k$ for $\ut_k=0$. It is composed of $\Gamma_{k,out}$ (in red), $\Gamma_{k,res}$ (in green) and $\Gamma_{k,p}$ (in blue). The section of $\Gamma_k$ which lies inside the ball $B_\varepsilon(\ut_k)$ (i.e. the reunion of $\Gamma_{k,res}$ and $\Gamma_{k,p}$) is notated $\Gamma_{k,in}$.}
			\label{chem}
		\end{center}
	\end{figure}

	There just remains to define the path $\Gamma_k$ defined on the Figure \ref{chem}. As we can see, it follows the ray $-\eta+i[-\pi,\pi]$ and is deformed inside $B_\varepsilon(\ut_k)$ into the path $\Gamma_{k,p}$. We define
	\begin{align*}
		\Gamma_{k,res}:=&\lc-\eta +it, t\in[-\pi,\pi]\backslash[\theta_k-\ell_{k,p},\theta_k+ \ell_{k,p}]\rc\cap B_\varepsilon(\ut_k),\\ \Gamma_{k,out}:=&\lc-\eta +it, t\in[-\pi,\pi]\rc\cap B_\varepsilon(\ut_k)^c,\\
		\Gamma_{k,in}:=&\Gamma_{k,p}\cup\Gamma_{k,res},\\
		\Gamma_k:= & \Gamma_{k,in} \cup \Gamma_{k,out}.
	\end{align*}
	
	Using Cauchy's formula and taking into account the "$2i\pi$-periodicity" of $\Gg_j(\tau)$, we have that for all $n\in\N^*$ and $j\in\Z$
	\begin{equation}
		\Gcc_j^n = \frac{1}{2i\pi}\int_{\Gamma_\rho} e^{n\tau}\Gg_j(\tau)d\tau= \frac{1}{2i\pi}\int_{\Gamma_k} e^{n\tau}\Gg_j(\tau)d\tau. \label{new_exp_G}
	\end{equation}

	In order to prove Lemma \ref{lemEstInterm}, we will start by proving the following lemma.
	
	\begin{lemma}\label{lemEstIntermMVC}
		For all $k\in\lc1,\ppp,K\rc$ and for all $s\in \N^*$, there exist two positive constants $C,c$ such that for all $(n,j)\in \Dc_k$
		$$\left|\Gcc_j^n -\frac{1}{2i\pi} \int_{\Gamma_{k,in}}P_{s,k}(\tau)\left(\sum_{l=0}^{s-1}\frac{(j R_{s,k}(\tau))^l}{l!}\right)e^{n\tau}e^{j\varphi_k(\tau)}d\tau\right|\leq \frac{C}{n^\frac{s+1}{2\mu_k}}\exp\left(-c\left(\frac{|n\alpha_k-j|}{n^\frac{1}{2\mu_k}}\right)^\frac{2\mu_k}{2\mu_k-1}\right).$$
	\end{lemma}
	
	Our main focus now will be to prove Lemma \ref{lemEstIntermMVC}. We observe that the triangular inequality implies
	\begin{equation}\label{egDecoup}
		\left|\Gcc_j^n -\frac{1}{2i\pi} \int_{\Gamma_{k,in}}P_{s,k}(\tau)\left(\sum_{l=0}^{s-1}\frac{(j R_{s,k}(\tau))^l}{l!}\right)e^{n\tau}e^{j\varphi_k(\tau)}d\tau\right| \leq \frac{1}{2\pi}\sum_{l=1}^8 E_l
	\end{equation}
	where
	\begin{align*}
		E_1& =  \left|\int_{\Gamma_{k,out}}e^{n\tau}\Gg_j(\tau)d\tau\right|, & E_2& = \left| \int_{\Gamma_{k,in}}e^{n\tau}\left(\Gg_j(\tau)-e^\tau g_k(\tau) \exp(j\varpi_k(\tau))\right)d\tau \right|  , \\
		E_3& =  \left| \int_{\Gamma_{k,p}}e^{n\tau+j\varpi_k(\tau)} \left( e^\tau g_k(\tau) - P_{s,k}(\tau)\right)d\tau \right| , & E_4& =  \left| \int_{\Gamma_{k,res}}e^{n\tau+j\varpi_k(\tau)} \left( e^\tau g_k(\tau) - P_{s,k}(\tau)\right)d\tau \right|, \\
		E_5& =  \left| \int_{\Gamma_{k,p}}P_{s,k}(\tau)e^{n\tau} \left( e^{j\varpi_k(\tau)} - e^{jQ_{s,k}(\tau)}\right)d\tau \right|, & E_6& =  \left| \int_{\Gamma_{k,res}}P_{s,k}(\tau)e^{n\tau} \left( e^{j\varpi_k(\tau)} - e^{jQ_{s,k}(\tau)}\right)d\tau \right|,\\
	\end{align*}
	$$E_7 =  \left|\int_{\Gamma_{k,p}} P_{s,k}(\tau) e^{n\tau+j\varphi_k(\tau)} \left(e^{jR_{s,k}(\tau)}-\sum_{l=0}^{s-1}\frac{(jR_{s,k}(\tau))^l}{l!}\right) d\tau\right|,$$
	$$E_8 =  \left|\int_{\Gamma_{k,res}} P_{s,k}(\tau) e^{n\tau+j\varphi_k(\tau)} \left(e^{jR_{s,k}(\tau)}-\sum_{l=0}^{s-1}\frac{(jR_{s,k}(\tau))^l}{l!}\right) d\tau\right|.$$
	
	We will now have to determine estimates on all these terms depending on $k$ (case\textbf{ I}, \textbf{II} and \textbf{III}) and also on $\tau_p$ and $\Gamma_{k,p}$:
	\begin{itemize}
		\item \textbf{Case A:} $\rho_k\left(\frac{\zeta_k}{\gamma_k}\right) \in \left[-\frac{\eta}{2},\varepsilon_{k,0}\right]$,
		\item \textbf{Case B:} $\rho_k\left(\frac{\zeta_k}{\gamma_k}\right) >\varepsilon_{k,0}$,
		\item \textbf{Case C:} $\rho_k\left(\frac{\zeta_k}{\gamma_k}\right) <-\frac{\eta}{2}$.
	\end{itemize}
	The main contribution will come from the terms $E_3$, $E_5$ and $E_7$ in the case \textbf{A}. We will prove much sharper estimates for the other terms.
	
	\subsubsection{Preliminary lemmas}
	
	Before we start to determine the estimates on the different terms, we are going to introduce some lemmas to simplify the redaction. Those lemmas assemble inequalities in the different cases (\textbf{A}, \textbf{B} and \textbf{C}) for which the proofs are similar with variations depending on the case we are in. They mainly rely on the inequalities \eqref{estVarphi}, \eqref{estVarpi} and \eqref{estR}. The proofs of those lemmas can be found in the appendix.
	
	We start with a lemma which will be useful to study the terms $E_5$, $E_6$, $E_7$ and $E_8$.
	
	\begin{lemma}[Inequalities in $B_{\varepsilon_\star}(\ut_k)$]\label{est_Beps}
		There exists $C>0$ such that for all $\tau \in B_{\varepsilon_\star}(\ut_k)$ and $(n,j)\in\Dc_k$, we have 
		$$ \left|e^{n\tau} \left( e^{j\varpi_k(\tau)} -e^{jQ_{s,k}(\tau)}\right)\right| \leq C n|\tau-\ut_k|^{2\mu_k+s} \exp(n \Re(\tau-\ut_k)+j (\Re(\varpi_k(\tau)) + \mathrm{sgn}(\alpha_k)|\xi_{s,k}(\tau)(\tau-\ut_k)^{2\mu_k+s+1}|)) $$
		and
		$$\left|e^{n\tau+j\varphi_k(\tau)} \left( e^{jR_{s,k}(\tau)} -\sum_{l=0}^{s-1}\frac{(jR_{s,k}(\tau))^l}{l!}\right)\right| \leq C \left(n|\tau-\ut_k|^{2\mu_k+1}\right)^s \exp(n \Re(\tau-\ut_k)+j (\Re(\varphi_k(\tau)) + \mathrm{sgn}(\alpha_k)|R_{s,k}(\tau)|)). $$
	\end{lemma}
	
	This next lemma will be useful for terms where the integral is defined along the path $\Gamma_{k,p}$ (terms $E_3$, $E_5$ and $E_7$).
	
	\begin{lemma}[Inequalities on $\Gamma_{k,p}$]\label{ine_taup}
		For $(n,j)\in\N^*\times \Z$ such that $\mathrm{sgn}(j)=\mathrm{sgn}(\alpha_k)$ and $\tau\in \Gamma_{k,p}$, we have
		
		$\bullet$ Case \textbf{A}: $\rho_k\left(\frac{\zeta_k}{\gamma_k}\right) \in \left[-\frac{\eta}{2},\varepsilon_{k,0}\right]$
		\begin{align}
			n\Re(\tau-\ut_k)+j(\Re(\varpi_k(\tau)) + \mathrm{sgn}(\alpha_k)|\xi_{s,k}(\tau)(\tau-\ut_k)^{2\mu_k+s}|) &\leq -nc_\star\Im(\tau-\ut_k)^{2\mu_k} - \frac{n}{\alpha_k} (2\mu_k-1)\gamma_k\left(\frac{|\zeta_k|}{\gamma_k}\right)^\frac{2\mu_k}{2\mu_k-1},\label{ine_varpi}\\
			n\Re(\tau-\ut_k)+j\Re(\varpi_k(\tau)) &\leq -nc_\star\Im(\tau-\ut_k)^{2\mu_k} - \frac{n}{\alpha_k} 	(2\mu_k-1)\gamma_k\left(\frac{|\zeta_k|}{\gamma_k}\right)^\frac{2\mu_k}{2\mu_k-1},\label{ine_varpi2}\\
			n\Re(\tau-\ut_k)+j(\Re(\varphi_k(\tau)) + \mathrm{sgn}(\alpha_k)|R_{s,k}(\tau)|) &\leq -nc_\star\Im(\tau-\ut_k)^{2\mu_k} - \frac{n}{\alpha_k} 	(2\mu_k-1)\gamma_k\left(\frac{|\zeta_k|}{\gamma_k}\right)^\frac{2\mu_k}{2\mu_k-1}.\label{ine_varpi3}
		\end{align}	
		
		$\bullet$ Case \textbf{B}: $\rho_k\left(\frac{\zeta_k}{\gamma_k}\right) >\varepsilon_{k,0}$
		\begin{align}
			n\Re(\tau-\ut_k)+j(\Re(\varpi_k(\tau)) + \mathrm{sgn}(\alpha_k)|\xi_{s,k}(\tau)(\tau-\ut_k)^{2\mu_k+s}|)&\leq -\frac{n}{\alpha_k}(2\mu_k-1)A_R\ud_k\varepsilon_{k,0}^{2\mu_k},\label{ine_varpi_cas2}\\
			n\Re(\tau-\ut_k)+j\Re(\varpi_k(\tau)) &\leq -\frac{n}{\alpha_k}(2\mu_k-1)A_R\ud_k\varepsilon_{k,0}^{2\mu_k},\label{ine_varpi2_cas2}\\
			n\Re(\tau-\ut_k)+j(\Re(\varphi_k(\tau)) + \mathrm{sgn}(\alpha_k)|R_{s,k}(\tau)|) & \leq -\frac{n}{\alpha_k}(2\mu_k-1)A_R\ud_k\varepsilon_{k,0}^{2\mu_k}.\label{ine_varpi3_cas2}
		\end{align}
		
		$\bullet$ Case \textbf{C}: $\rho_k\left(\frac{\zeta_k}{\gamma_k}\right) <-\frac{\eta}{2}$
		\begin{align}
			n\Re(\tau-\ut_k)+j(\Re(\varpi_k(\tau)) + \mathrm{sgn}(\alpha_k)|\xi_{s,k}(\tau)(\tau-\ut_k)^{2\mu_k+s}|) &\leq -\frac{n}{\alpha_k}(2\mu_k-1)A_R\ud_k\left(\frac{\eta}{2}\right)^{2\mu_k},\label{ine_varpi_cas3}\\
			n\Re(\tau-\ut_k)+j\Re(\varpi_k(\tau)) &\leq -\frac{n}{\alpha_k}(2\mu_k-1)A_R\ud_k\left(\frac{\eta}{2}\right)^{2\mu_k},\label{ine_varpi2_cas3}\\
			n\Re(\tau-\ut_k)+j(\Re(\varphi_k(\tau)) + \mathrm{sgn}(\alpha_k)|R_{s,k}(\tau)|) & \leq-\frac{n}{\alpha_k}(2\mu_k-1)A_R\ud_k\left(\frac{\eta}{2}\right)^{2\mu_k}.\label{ine_varpi3_cas3}
		\end{align}
	\end{lemma}
	
		Finally, we introduce in the next lemma some inequalities that will help us for the terms with integrals defined on $\Gamma_{k,res}$ (terms $E_4$, $E_6$ and $E_8$).
	
	\begin{lemma}[Inequalities on $\Gamma_{k,res}$]\label{ine_res}
		For $(n,j)\in\N^*\times \Z$ such that $\mathrm{sgn}(j)=\mathrm{sgn}(\alpha_k)$ and $\tau\in \Gamma_{k,res}$, we have in all cases
		\begin{align}
			n\Re(\tau-\ut_k)+j(\Re(\varpi_k(\tau))+ \mathrm{sgn}(\alpha_k) \left|\xi_{s,k}(\tau)(\tau-\ut_k)^{2\mu_k+s}\right|) &\leq -n\frac{\eta}{2}, \label{ine_varpi_res}\\
			n\Re(\tau-\ut_k)+j(\Re(\varpi_k(\tau)) &\leq -n\frac{\eta}{2}, \label{ine_varpi_res2}\\
			n\Re(\tau-\ut_k)+j(\Re(\varphi_k(\tau)) + \mathrm{sgn}(\alpha_k)|R_{s,k}(\tau)|)&\leq -n\frac{\eta}{2}. \label{ine_varpi_res3}
		\end{align}
	\end{lemma}
	
	\subsubsection{Estimates of part of the terms}
	
	We are going to first prove estimates for the terms where the proof will not depend on the case\textbf{ A}, \textbf{B} or \textbf{C} in which we are.
	
	\underline{$\bullet$ Estimate for $E_2$:}
	
	We introduce the path $\Gamma_{\eta,k}$ defined as
	$$\Gamma_{\eta,k}:=\lc -\eta+it, t\in[-\pi,\pi]\rc \cap B_\varepsilon(\ut_k).$$
	Using Cauchy's formula, we have that
	$$\int_{\Gamma_{k,in}}e^{n\tau}\left(\Gg_j(\tau)-e^\tau g_k(\tau) \exp(j\varpi_k(\tau))\right)d\tau= \int_{\Gamma_{\eta,k}}e^{n\tau}\left(\Gg_j(\tau)-e^\tau g_k(\tau) \exp(j\varpi_k(\tau))\right)d\tau.$$
	Because we supposed that $\alpha_k>0$, depending on whether we are in case\textbf{ I} or \textbf{III}, the previous equality and the inequalities \eqref{varpi_près_pos_I} and \eqref{varpi_près_pos_III} imply
	$$\left| \int_{\Gamma_{k,in}}e^{n\tau}\left(\Gg_j(\tau)-e^\tau g_k(\tau) \exp(j\varpi_k(\tau))\right)d\tau \right|  \lesssim e^{-n\eta - c j}.$$
	
	\underline{$\bullet$ Estimate for $E_4$:}

	The inequality \eqref{ine_varpi_res2} implies
	$$\left| \int_{\Gamma_{k,res}}e^{n\tau+j\varpi_k(\tau)} \left( e^\tau g_k(\tau) - P_{s,k}(\tau)\right)d\tau \right| \lesssim \int_{\Gamma_{k,res}} \exp\left(n\Re(\tau)+ j\Re(\varpi_k(\tau))\right)|d\tau| \lesssim e^{-n\frac{\eta}{2}} .$$
	
	\underline{$\bullet$ Estimate for $E_6$:}
	
	If we use Lemma \ref{est_Beps}, we have
	\begin{multline*}
		\left| \int_{\Gamma_{k,res}}P_{s,k}(\tau)e^{n\tau} \left( e^{j\varpi_k(\tau)} - e^{jQ_{s,k}(\tau)}\right)d\tau \right|
		\\ \lesssim \int_{\Gamma_{k,res}} \exp(n \Re(\tau-\ut_k)+j (\Re(\varpi_k(\tau)) + \mathrm{sgn}(\alpha_k)|\xi_{s,k}(\tau)(\tau-\ut_k)^{2\mu_k+s}|))n|\tau-\ut_k|^{2\mu_k+s} |d\tau|.
	\end{multline*}
	
	Therefore, the inequality \eqref{ine_varpi_res} implies
	
	$$ 	\left| \int_{\Gamma_{k,res}}P_{s,k}(\tau)e^{n\tau} \left( e^{j\varpi_k(\tau)} - e^{jQ_{s,k}(\tau)}\right)d\tau\right|	\lesssim   ne^{-n\frac{\eta}{2}}
	\lesssim  e^{-n\frac{\eta}{4}}.$$
	
	\underline{$\bullet$ Estimate for $E_8$:}
	
	If we use Lemma \ref{est_Beps}, we have
	\begin{multline*}
		\left| \int_{\Gamma_{k,res}}P_{s,k}(\tau)e^{n\tau+j\varphi_k(\tau)} \left( e^{jR_{s,k}(\tau)} - \sum_{l=0}^{s-1}\frac{(jR_{s,k}(\tau))^l}{l!}\right)d\tau \right|
		\\ \lesssim \int_{\Gamma_{k,res}} \exp(n \Re(\tau-\ut_k)+j (\Re(\varpi_k(\tau)) + \mathrm{sgn}(\alpha_k)|R_{s,k}(\tau)|))(n|\tau-\ut_k|^{2\mu_k+1})^s |d\tau|.
	\end{multline*}
	
	Therefore, the inequality \eqref{ine_varpi_res3} implies
	
	$$ 	\left| \int_{\Gamma_{k,res}}P_{s,k}(\tau)e^{n\tau+j\varphi_k(\tau)} \left( e^{jR_{s,k}(\tau)} - \sum_{l=0}^{s-1}\frac{(jR_{s,k}(\tau))^l}{l!}\right)d\tau \right|	\lesssim   n^{s}e^{-n\frac{\eta}{2}}
	\lesssim  e^{-n\frac{\eta}{4}}.$$
	
	It remains to study the terms $E_1$, $E_3$, $E_5$ and $E_7$.
	
	\subsubsection{The terms $E_3$, $E_5$ and $E_7$, Case A : $\rho_k\left(\frac{\zeta_k}{\gamma_k}\right) \in \left[-\frac{\eta}{2},\varepsilon_{k,0}\right]$}
	
	This part of the proof is the most important because those terms will create the limiting estimates.
	
	\underline{$\bullet$ Estimate for $E_3$:}
	
	Because of Taylor's theorem, we have
	$$E_3  =  \left| \int_{\Gamma_{k,p}}\left( e^\tau g_k(\tau) - P_{s,k}(\tau)\right)e^{n\tau+j\varpi_k(\tau)} d\tau \right| \lesssim \int_{\Gamma_{k,p}} |\tau-\ut_k|^s\exp\left(n\Re(\tau)+ j\Re(\varpi_k(\tau))\right)|d\tau|.$$
	
	The inequality \eqref{ine_varpi2} implies 
	$$ E_3 \lesssim  \int_{\Gamma_{k,p}} |\tau-\ut_k|^s e^{-nc_\star\Im(\tau-\ut_k)^{2\mu_k}}|d\tau| \exp\left(-\frac{n}{\alpha_k} (2\mu_k-1)\gamma_k\left(\frac{|\zeta_k|}{\gamma_k}\right)^\frac{2\mu_k}{2\mu_k-1}\right).$$
	
	But, the inequality \eqref{ineg_gamma} and the fact that $\rho_k\left(\frac{\zeta_k}{\gamma_k}\right)=\tau_p$ imply
	$$\frac{n}{\alpha_k} (2\mu_k-1)\gamma_k\left(\frac{|\zeta_k|}{\gamma_k}\right)^\frac{2\mu_k}{2\mu_k-1}\geq \frac{2\mu_k-1}{\alpha_k}A_R\ud_k n|\tau_p|^{2\mu_k}.$$
	
	If we introduce $c>0$ small enough, then
	$$ E_3 \lesssim  \int_{\Gamma_{k,p}} |\tau-\ut_k|^s e^{-nc_\star\Im(\tau-\ut_k)^{2\mu_k}}|d\tau| \exp\left(-cn|\tau_p|^{2\mu_k}\right).$$
	
	Using the parametrization \eqref{param} and the inequality \eqref{hp}, we have
	$$\int_{\Gamma_{k,p}} |\tau-\ut_k|^s e^{-nc_\star\Im(\tau-\ut_k)^{2\mu_k}}|d\tau| \lesssim \int_{-\ell_{k,p}}^{\ell_{k,p}}(|\tau_p|^s+|t|^s)e^{-nc_\ast t^{2\mu_k}}dt. $$
	
	The change of variables $u=n^\frac{1}{2\mu_k}t$ and the fact that the function $\displaystyle x\geq0\mapsto x^s\exp\left(-\frac{c}{2}x^{2\mu_k}\right)$ is bounded imply
	$$\lc \begin{array}{c}
		\displaystyle \int_{-\ell_{k,p}}^{\ell_{k,p}}|t|^se^{-nc_\ast t^{2\mu_k}}dt\lesssim \frac{1}{n^\frac{s+1}{2\mu_k}},\\
		\displaystyle \int_{-\ell_{k,p}}^{\ell_{k,p}}|\tau_p|^se^{-nc_\ast t^{2\mu_k}}dt\lesssim \frac{1}{n^\frac{s+1}{2\mu_k}}\exp\left(\frac{c}{2}n|\tau_p|^{2\mu_k}\right).
	\end{array}\right.$$
	
	Thus,
	$$E_3 \lesssim \frac{1}{n^\frac{s+1}{2 \mu_k}} \exp\left(-\frac{c}{2}n|\tau_p|^{2\mu_k}\right). $$
	
	Lastly, the inequality \eqref{ineg_gamma} implies that we have a constant $\tilde{c}>0$ independent from $j$ and $n$ such that 
	$$\frac{c}{2}n|\tau_p|^{2\mu_k} \geq \tilde{c}\left(\frac{|j-n\alpha_k|}{n^\frac{1}{2\mu_k}}\right)^\frac{2\mu_k}{2\mu_k-1}$$
	so,
	$$E_3 \lesssim \frac{1}{n^\frac{s+1}{2\mu_k}} \exp\left(-\tilde{c} \left(\frac{|j-n\alpha_k|}{n^\frac{1}{2\mu_k}}\right)^\frac{2\mu_k}{2\mu_k-1}\right).$$	
	
	\underline{$\bullet$ Estimate for $E_5$:}
	
	Using Lemma \ref{est_Beps} and the inequality \eqref{ine_varpi}, we have
	\begin{align*}
		E_5 &=  \left| \int_{\Gamma_{k,p}}P_{s,k}(\tau)e^{n\tau}	 \left( e^{j\varpi_k(\tau)} -  e^{jQ_{s,k}(\tau)}\right)d\tau \right|\\
		&\lesssim  \int_{\Gamma_{k,p}} n|\tau-\ut_k|^{2\mu_k+s} \exp(n \Re(\tau-\ut_k)+j (\Re(\varpi_k(\tau)) + \mathrm{sgn}(\alpha_k)|\xi_{s,k}(\tau)(\tau-\ut_k)^{2\mu_k+s}|)) |d\tau|\\
		& \lesssim  \exp\left(- \frac{n}{\alpha_k} (2\mu_k-1)\gamma_k\left(\frac{|\zeta_k|}{\gamma_k}\right)^\frac{2\mu_k}{2\mu_k-1}\right)n\int_{\Gamma_{k,p}} |\tau-\ut_k|^{2\mu_k+s} \exp(-nc_\star\Im(\tau-\ut_k)^{2\mu_k}) |d\tau|.
	\end{align*}
	
	Just like in the estimate for the previous term, because of the inequality \eqref{ineg_gamma}, if we introduce $c>0$ small enough, we have
	
	$$E_5 \lesssim n \exp\left(-cn|\tau_p|^{2\mu_k}\right) \int_{\Gamma_{k,p}}|\tau-\ut_k|^{2\mu_k+s}\exp(-nc_\star\Im(\tau-\ut_k)^{2\mu_k})|d\tau|.$$
	
	The same reasoning as for the estimate of $E_3$ implies that
	$$n\int_{\Gamma_{k,p}}|\tau-\ut_k|^{2\mu_k+s}\exp(-nc_\star\Im(\tau-\ut_k)^{2\mu_k})|d\tau| \lesssim n \int_{-\ell_{k,p}}^{\ell_{k,p}}|t|^{2\mu_k+s}e^{-nc_\ast t^{2\mu_k}}dt +n \int_{-\ell_{k,p}}^{\ell_{k,p}}|\tau_p|^{2\mu_k+s}e^{-nc_\ast t^{2\mu_k}}dt.$$
	
	The change of variables $u=n^\frac{1}{2\mu_k}t$ and the fact that the function $\displaystyle x\geq0\mapsto x^{2\mu_k+s}\exp\left(-\frac{c}{2}x^{2\mu_k}\right)$ is bounded imply
	$$\lc \begin{array}{c}
		\displaystyle n\int_{-\ell_{k,p}}^{\ell_{k,p}}|t|^{2\mu_k+s}e^{-nc_\ast t^{2\mu_k}}dt\lesssim \frac{1}{n^\frac{s+1}{2\mu_k}},\\
		\displaystyle n\int_{-\ell_{k,p}}^{\ell_{k,p}}|\tau_p|^{2\mu_k+s}e^{-nc_\ast t^{2\mu_k}}dt\lesssim \frac{1}{n^\frac{s+1}{2\mu_k}}\exp\left(\frac{c}{2}n|\tau_p|^{2\mu_k}\right).
	\end{array}\right.$$
	
	Thus,
	$$E_5 \lesssim \frac{1}{n^\frac{s+1}{2\mu_k}} \exp\left(-\frac{c}{2}n|\tau_p|^{2\mu_k}\right). $$
	
	Lastly, the inequality on $\gamma_k$ \eqref{ineg_gamma} implies that we have a constant $\tilde{c}>0$ independent from $j$ and $n$ such that 
	$$\frac{c}{2}n|\tau_p|^{2\mu_k} \geq \tilde{c}\left(\frac{|j-n\alpha_k|}{n^\frac{1}{2\mu_k}}\right)^\frac{2\mu_k}{2\mu_k-1}$$
	so,
	$$E_5 \lesssim \frac{1}{n^\frac{s+1}{2\mu_k}} \exp\left(-\tilde{c} \left(\frac{|j-n\alpha_k|}{n^\frac{1}{2\mu_k}}\right)^\frac{2\mu_k}{2\mu_k-1}\right).$$

	\underline{$\bullet$ Estimate for $E_7$:}
	
	Using Lemma \ref{est_Beps} and the inequality \eqref{ine_varpi3}, we have
	\begin{align*}
		E_7 &=  \left| \int_{\Gamma_{k,p}}P_{s,k}(\tau)e^{n\tau+j\varphi_k(\tau)}	 \left( e^{jR_{s,k}(\tau)} - \sum_{l=0}^{s-1}\frac{(jR_{s,k}(\tau))^l}{l!}\right) d\tau\right|\\
		&\lesssim  \int_{\Gamma_{k,p}} (n|\tau-\ut_k|^{2\mu_k+1})^s \exp(n \Re(\tau-\ut_k)+j (\Re(\varphi_k(\tau)) + \mathrm{sgn}(\alpha_k)|R_{s,k}(\tau)|)) |d\tau|\\
		& \lesssim  \exp\left(- \frac{n}{\alpha_k} (2\mu_k-1)\gamma_k\left(\frac{|\zeta_k|}{\gamma_k}\right)^\frac{2\mu_k}{2\mu_k-1}\right)n^s\int_{\Gamma_{k,p}} |\tau-\ut_k|^{s(2\mu_k+1)} \exp(-nc_\star\Im(\tau-\ut_k)^{2\mu_k}) |d\tau|.
	\end{align*}
	
	Just like in the estimate for the previous term, because of the inequality \eqref{ineg_gamma}, if we introduce $c>0$ small enough, we have
	
	$$E_7 \lesssim n^s \exp\left(-cn|\tau_p|^{2\mu_k}\right) \int_{\Gamma_{k,p}}|\tau-\ut_k|^{s(2\mu_k+1)}\exp(-nc_\star\Im(\tau-\ut_k)^{2\mu_k})|d\tau|.$$
	
	The same reasoning as for the estimate of $E_3$ implies that
	\begin{multline*}
		n^s\int_{\Gamma_{k,p}}|\tau-\ut_k|^{s(2\mu_k+1)}\exp(-nc_\star\Im(\tau-\ut_k)^{2\mu_k})|d\tau| \\ \lesssim n^s \int_{-\ell_{k,p}}^{\ell_{k,p}}|t|^{s(2\mu_k+1)}e^{-nc_\ast t^{2\mu_k}}dt +n^s \int_{-\ell_{k,p}}^{\ell_{k,p}}|\tau_p|^{s(2\mu_k+1)}e^{-nc_\ast t^{2\mu_k}}dt.
	\end{multline*}
	
	The change of variables $u=n^\frac{1}{2\mu_k}t$ and the fact that the function $\displaystyle x\geq0\mapsto x^{s(2\mu_k+1)}\exp\left(-\frac{c}{2}x^{2\mu_k}\right)$ is bounded imply
	$$\lc \begin{array}{c}
		\displaystyle n^s\int_{-\ell_{k,p}}^{\ell_{k,p}}|t|^{s(2\mu_k+1)}e^{-nc_\ast t^{2\mu_k}}dt\lesssim \frac{1}{n^\frac{s+1}{2\mu_k}},\\
		\displaystyle n^s\int_{-\ell_{k,p}}^{\ell_{k,p}}|\tau_p|^{s(2\mu_k+1)}e^{-nc_\ast t^{2\mu_k}}dt\lesssim \frac{1}{n^\frac{s+1}{2\mu_k}}\exp\left(\frac{c}{2}n|\tau_p|^{2\mu_k}\right).
	\end{array}\right.$$
	
	Thus,
	$$E_7 \lesssim \frac{1}{n^\frac{s+1}{2\mu_k}} \exp\left(-\frac{c}{2}n|\tau_p|^{2\mu_k}\right). $$
	
	Lastly, the inequality on $\gamma_k$ \eqref{ineg_gamma} implies that we have a constant $\tilde{c}>0$ independent from $j$ and $n$ such that 
	$$\frac{c}{2}n|\tau_p|^{2\mu_k} \geq \tilde{c}\left(\frac{|j-n\alpha_k|}{n^\frac{1}{2\mu_k}}\right)^\frac{2\mu_k}{2\mu_k-1}$$
	so,
	$$E_7 \lesssim \frac{1}{n^\frac{s+1}{2\mu_k}} \exp\left(-\tilde{c} \left(\frac{|j-n\alpha_k|}{n^\frac{1}{2\mu_k}}\right)^\frac{2\mu_k}{2\mu_k-1}\right).$$	
	
	\subsubsection{The terms $E_3$, $E_5$ and $E_7$, Case B and C:}
	
	We now consider that we are either in case\textbf{ B} or case\textbf{ C} (i.e. $\rho_k(\frac{\zeta_k}{\gamma_k})\notin\left[-\frac{\eta}{2},\varepsilon_{k,0}\right]$).
	
	\underline{$\bullet$ Estimate for $E_3$:}
	
	Because of Taylor's theorem, we have
	$$E_3  =  \left| \int_{\Gamma_{k,p}}e^{n\tau+j\varpi_k(\tau)} \left( e^\tau g_k(\tau) - P_{s,k}(\tau)\right)d\tau \right| \lesssim \int_{\Gamma_{k,p}} |\tau-\ut_k|^s\exp\left(n\Re(\tau)+ j\Re(\varpi_k(\tau))\right)|d\tau|.$$
	
	Using the inequality \eqref{ine_varpi2_cas2} or \eqref{ine_varpi2_cas3} whether we are in case\textbf{ B} or \textbf{C}, they imply that there exists $c>0$ independent from $j$ and $n$ such that 
	$$E_3\lesssim e^{-cn}. $$
	
	\underline{$\bullet$ Estimate for $E_5$:}
	
	Using Lemma \ref{est_Beps}, we have
	\begin{align*}
		E_5 &=  \left| \int_{\Gamma_{k,p}}e^{n\tau}P_{s,k}(\tau) \left( e^{j\varpi_k(\tau)} - e^{jQ_{s,k}(\tau)}\right)d\tau \right|\\
		&\lesssim  \int_{\Gamma_{k,p}} n|\tau-\ut_k|^{2\mu_k+s} \exp(n \Re(\tau-\ut_k)+j (\Re(\varpi_k(\tau)) +\mathrm{sgn}(\alpha_k)|\xi_{s,k}(\tau)(\tau-\ut_k)^{2\mu_k+s}|)) |d\tau|.
	\end{align*}
	
	Using the inequality \eqref{ine_varpi_cas2} or \eqref{ine_varpi_cas3} whether we are in case\textbf{ B} or \textbf{C}, they imply that there exists $c>0$ independent from $j$ and $n$ such that 
	$$E_5\lesssim ne^{-cn} \lesssim e^{-\frac{c}{2}n}. $$
	
	\underline{$\bullet$ Estimate for $E_7$:}
	
	Using Lemma \ref{est_Beps}, we have
	\begin{align*}
		E_7 &=  \left| \int_{\Gamma_{k,p}}P_{s,k}(\tau)e^{n\tau+j\varphi_k(\tau)}	 \left( e^{jR_{s,k}(\tau)} - \sum_{l=0}^{s-1}\frac{(jR_{s,k}(\tau))^l}{l!}\right) d\tau\right|\\
		&\lesssim  \int_{\Gamma_{k,p}} (n|\tau-\ut_k|^{2\mu_k+1})^s \exp(n \Re(\tau-\ut_k)+j (\Re(\varphi_k(\tau)) + \mathrm{sgn}(\alpha_k)|R_{s,k}(\tau)|)) |d\tau|.
	\end{align*}
		
	Using the inequality \eqref{ine_varpi3_cas2} or \eqref{ine_varpi3_cas3} whether we are in case\textbf{ B} or \textbf{C}, they imply that there exists $c>0$ independent from $j$ and $n$ such that 
	$$E_7\lesssim n^se^{-cn} \lesssim e^{-\frac{c}{2}n}. $$
	
	\subsubsection{Estimate for the term $E_1$}
	
	\underline{$\bullet$ Estimate for $E_1$:}
	
	We recall that
	$$E_1 =  \left|\int_{\Gamma_{k,out}} e^{n\tau} \Gg_j(\tau) d\tau\right|.$$
	
	For $\tau \in \Gamma_{k,out}$, we have different estimates depending on whether we are inside a ball $B_\varepsilon(\ut_l)$ or not. Therefore, we introduce the set of distinct points $$\lc\hat{\tau}_1,\ppp,\hat{\tau}_R\rc = \lc\ut_l, \quad l\in\lc1,\ppp,K\rc\rc\backslash \lc\ut_k\rc.$$
	It allows us to decompose the path $\Gamma_{k,out}$ as 
	$$\Gamma_{k,out}:= \bigcup_{l=0}^R \widehat{\Gamma}_l,$$
	where for all $l\in\lc1,\ppp,R\rc$
	$$\widehat{\Gamma}_l := \Gamma_{k,out}\cap B_\varepsilon(\hat{\tau}_l)$$
	and
	$$\widehat{\Gamma}_0 := \Gamma_{k,out}\backslash \bigcup_{l=1}^R \widehat{\Gamma}_l.$$

	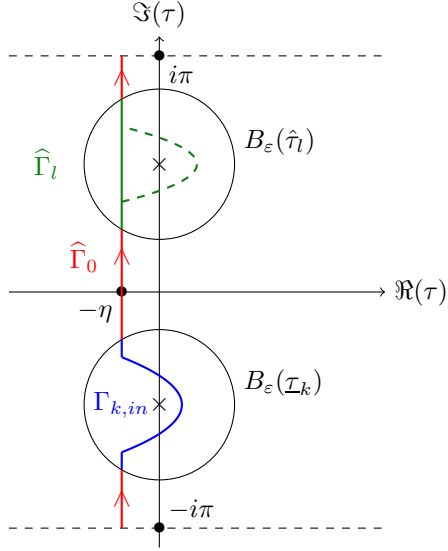
\begin{figure}
		\begin{center}
			\begin{tikzpicture}[scale=1]
				\draw[->] (-2,0) -- (3,0) node[right] {$\Re(\tau)$};
				\draw[->] (0,-3.4) -- (0,3.4) node[above] {$\Im(\tau)$};
				\draw[dashed] (-2,pi) -- (3,pi);
				\draw[dashed] (-2,-pi) -- (3,-pi);
				\draw (0,pi) node {$\bullet$} node[below right] {$i\pi$};
				\draw (0,-pi) node {$\bullet$} node[above right] {$-i\pi$};
				\draw (-0.5,0) node {$\bullet$} node[below left] {$-\eta$};

				\draw (0,-1.5) node {$\times$} circle (1);
				\draw (1,-1.5) node[above right] {$B_\varepsilon(\ut_k)$};
				
				\draw (0,1.7) node {$\times$} circle (1);
				\draw (1,1.7) node[above right] {$B_\varepsilon(\hat{\tau}_l)$};
				
				\draw[red,thick] (-0.5,-pi) -- (-0.5,{-1.5-sqrt(1-0.5^2)}) node[midway, sloped] {$>$};
				\draw[red,thick] (-0.5,{-1.5+sqrt(1-0.5^2)}) -- (-0.5,{1.7-sqrt(1-0.5^2)}) node[near end, sloped] {$>$};
				\draw[red] (-1,{(-1.5+sqrt(1-0.5^2)+3*(1.7-sqrt(1-0.5^2)))/4}) node {$\widehat{\Gamma}_0$};
				\draw[red,thick] (-0.5,{1.7+sqrt(1-0.5^2)}) -- (-0.5,pi)  node[midway, sloped] {$>$};
				
				\draw[thick,blue] plot [samples = 100, domain={-sqrt(0.4)}:{sqrt(0.4)}] ({0.3-2*abs(\x)^2},\x-1.5) ;
				\draw[blue,thick] (-0.5,{-1.5-sqrt(0.75)}) -- (-0.5,{-1.5-sqrt(0.4)});
				\draw[blue,thick] (-0.5,{-1.5+sqrt(0.4)}) -- (-0.5,{-1.5+sqrt(0.75)});
				\draw[blue] (-0.5,-1.5) node {$\Gamma_{k,in}$};

				\draw[dartmouthgreen,thick] (-0.5,{1.7-sqrt(0.75)}) -- (-0.5,{1.7+sqrt(0.75)});
				\draw[dartmouthgreen] (-1.5,1.7) node {$\widehat{\Gamma}_l$};
				\draw[thick,dartmouthgreen,dashed] plot [samples = 100, domain={-1/2}:{1/2}] ({0.5-4*abs(\x)^2},\x+1.7) ;
			\end{tikzpicture}
			\caption{This is a representation of $\Gamma_k$ where we decompose $\Gamma_{k,out}$. The red path corresponds to $\widehat{\Gamma}_0$ the part of $\Gamma_{k,out}$ which lies outside the balls $B_\varepsilon(\hat{\tau}_l)$. The green path corresponds to $\widehat{\Gamma}_l$ the part of $\Gamma_{k,out}$ which lies inside the ball $B_\varepsilon(\hat{\tau}_l)$. The dashed green path corresponds to the deformation we use in the proof of the estimate for $E_1$. }
			\label{chem_2}
		\end{center}
	\end{figure}
	
	This decomposition of $\Gamma_{k,out}$ is represented on Figure \ref{chem_2}. The inequality \eqref{varpi_loin} gives us that
	
	$$\left|\int_{\widehat{\Gamma}_0} e^{n\tau} \Gg_j(\tau) d\tau\right|\lesssim e^{-n\eta-c|j|}.$$
	
	We now consider $l\in\lc1,\ppp,R\rc$. There are two possibilities because of Hypothesis \ref{H3}:
	
	$\bullet$ The set $\lc i\in\lc 1,\ppp,R\rc,  \quad \ut_i=\hat{\tau}_l\rc$ is the singleton $\lc i\rc$ with $\alpha_i<0$ (i.e. we are in case\textbf{ II}). Then, knowing that for $(n,j)\in\Dc_k$ we have $j\geq 1$, because of the inequality \eqref{varpi_près_pos_II}, we have 
	$$ \left|\int_{\widehat{\Gamma}_l} e^{n\tau} \Gg_j(\tau) d\tau\right|\lesssim e^{-n\eta-c|j|}.$$
	
	$\bullet$ The set $\lc i\in\lc 1,\ppp,R\rc,  \quad \ut_i=\hat{\tau}_l\rc$ is the singleton $\lc i\rc$ with $\alpha_i>0$ (i.e. we are in case\textbf{ I}) or it has two distinct elements $\lc i,j\rc$ with $\alpha_i>0$ and $\alpha_j<0$ (i.e. we are in case\textbf{ III}). Either way, the inequalities \eqref{varpi_près_pos_I} and \eqref{varpi_près_pos_III} imply that 
	$$ \left|\int_{\widehat{\Gamma}_l} e^{n\tau} \Gg_j(\tau) d\tau\right|\leq 2\pi Ce^{-n\eta-c|j|} +\left|\int_{\widehat{\Gamma}_l} \exp(n\tau+j\varpi_i(\tau))e^\tau g_i(\tau) d\tau\right|  .$$
	
	Just like we defined the path $\Gamma_{k,p}$, $\Gamma_{k,res}$ and $\Gamma_{k,in}:= \Gamma_{k,p}\sqcup\Gamma_{k,res}$, we can define a path $\Gamma_{i,p}$, $\Gamma_{i,res}$ and $\Gamma_{i,in}:= \Gamma_{i,p}\sqcup\Gamma_{i,res}$. The path $\Gamma_{i,in}$ is represented with a dashed green line on the Figure \ref{chem_2}. Using Cauchy's formula, we then have 
	$$\int_{\widehat{\Gamma}_l} \exp(n\tau+j\varpi_i(\tau))e^\tau g_i(\tau) d\tau =\int_{\Gamma_{i,in}} \exp(n\tau+j\varpi_i(\tau))e^\tau g_i(\tau) d\tau$$
	
	The function $\tau\mapsto e^\tau g_i(\tau)$ can be bounded so we just have to bound $\displaystyle\int_{\Gamma_{i,in}} \exp(n\Re(\tau-\ut_i)+j\Re(\varpi_i(\tau)))d|\tau|$. We observe that the proofs of the Lemmas \ref{ine_taup} and \ref{ine_res} are also true for $\Gamma_{i,p}$ and $\Gamma_{i,res}$. Using the inequality \eqref{ine_varpi_res2} for the integral along the path $\Gamma_{i,res}$, we prove that there exists a constant $c>0$ independent from $n$  and $j$ so that
	$$\int_{\Gamma_{i,res}} \exp(n\Re(\tau-\ut_i)+j\Re(\varpi_i(\tau)))d|\tau|\lesssim e^{-cn}.$$
	
	It remains to bound the integral along the path $\Gamma_{i,p}$. In the case\textbf{ A} (i.e. $\rho_i(\frac{\zeta_i}{\gamma_i})\in[-\frac{\eta}{2},\varepsilon_{i,0}]$), we observe that for $(n,j)\in\Dc_k$, $\gamma_i$ is bounded between two positive constants and
	$$|\zeta_i|\geq \frac{1}{2\mu_i}\min(|\alpha_i-\ud_k|,|\alpha_i-\sd_k|).$$
	Therefore, using the inequality \eqref{ine_varpi2} and the previous observation in case\textbf{ A} and using the inequalities \eqref{ine_varpi2_cas2} and \eqref{ine_varpi2_cas3} in cases\textbf{ B} and \textbf{C}, we prove that there exists a constant $c>0$ independent from $n$  and $j$ so that
	$$\int_{\Gamma_{i,p}} \exp(n\Re(\tau-\ut_i)+j\Re(\varpi_i(\tau)))d|\tau|\lesssim e^{-cn}.$$
	
	Therefore, there exists a constant $c>0$ such that
	$$\forall (n,j)\in \Dc_k, \quad \left|\int_{\widehat{\Gamma}_l} e^{n\tau} \Gg_j(\tau) d\tau\right| \lesssim e^{-cn}.$$
	This gives a sharp estimate of $E_1$.
	
	If we recapitulate the estimates we found, we can define two constants $C,c>0$ such that 
	$$\forall (n,j)\in\Dc_k, \forall l\in\lc1,2,4,6,8\rc, \quad E_l\leq Ce^{-cn},$$
	and 
	$$\forall (n,j)\in\Dc_k, \forall l\in\lc3,5,7\rc, \quad E_l\leq \frac{C}{n^\frac{s+1}{2\mu_k}} \exp\left(-c \left(\frac{|j-n\alpha_k|}{n^\frac{1}{2\mu_k}}\right)^\frac{2\mu_k}{2\mu_k-1}\right).$$
	
	The estimates we proved on all the terms and Lemma \ref{est_expn} allow us to conclude the proof of Lemma \ref{lemEstIntermMVC}.
	
	\subsubsection{From Lemma \ref{lemEstIntermMVC} to Lemma \ref{lemEstInterm}}
	
	Now that Lemma \ref{lemEstIntermMVC} is proved, we know that there exist two positive constants $C,c$ such that for all $(n,j)\in\Dc_k$,
	\begin{equation}\label{estInterm}
		\left|\Gcc_j^n -\frac{1}{2i\pi} \int_{\Gamma_{k,in}}P_{s,k}(\tau)\left(\sum_{l=0}^{s-1}\frac{(j R_{s,k}(\tau))^l}{l!}\right)e^{n\tau}e^{j\varphi_k(\tau)}d\tau\right|\leq \frac{C}{n^\frac{s+1}{2\mu_k}}\exp\left(-c\left(\frac{|n\alpha_k-j|}{n^\frac{1}{2\mu_k}}\right)^\frac{2\mu_k}{2\mu_k-1}\right).
	\end{equation}

	Proving Lemma \ref{lemEstInterm} amounts to proving a similar estimate as \eqref{estInterm} where the integration path would be $\lc it+\ut_k, t\in\R\rc$. This is the goal of this subsection. We prove the following lemma, which will use the conditions \eqref{inEta} and \eqref{inEta2} we introduced on $\eta_\varepsilon$.
	
	\begin{lemma}\label{lemDécal}
		We define the path
		$$\Gamma^0_{k, in} := \lc it, \quad t \in [\theta_k-r_\varepsilon(\eta),\theta_k+r_\varepsilon(\eta)]\rc$$
		where the function $r_\varepsilon$ is defined in \eqref{defR}. Then, for all $m\in \N^*$, there exist two positive constants $C,c$ such that
		$$\forall (n,j)\in\Dc_k, \quad \left|\int_{\Gamma^0_{k, in}}(\tau-\ut_k)^me^{n\tau}e^{j\varphi_k(\tau)}d\tau -\int_{\Gamma_{k, in}}(\tau-\ut_k)^me^{n\tau}e^{j\varphi_k(\tau)}d\tau \right|\leq Ce^{-cn}.$$
	\end{lemma}
	
	\begin{proof}
		As in Figure \ref{figDecal}, we define the paths
		$$\Gamma^+_{comp} := \lc t + i (\theta_k+r_\varepsilon(\eta)), \quad t \in[-\eta,0]\rc, \quad \Gamma^-_{comp} := \lc t + i (\theta_k-r_\varepsilon(\eta)), \quad t \in[-\eta,0]\rc.$$
		
		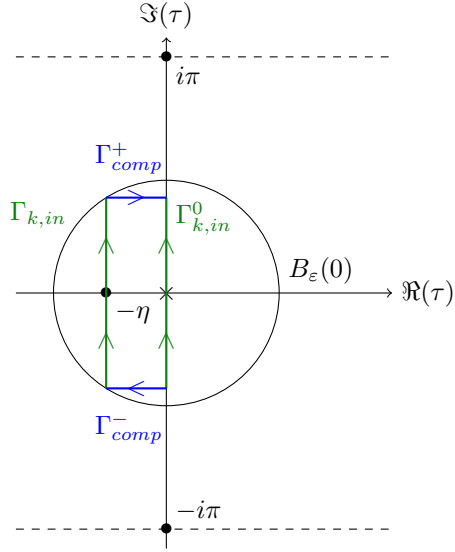
\begin{figure}
			\begin{center}
				\begin{tikzpicture}[scale=1]
					\draw[->] (-2,0) -- (3,0) node[right] {$\Re(\tau)$};
					\draw[->] (0,-3.4) -- (0,3.4) node[above] {$\Im(\tau)$};
					\draw[dashed] (-2,pi) -- (3,pi);
					\draw[dashed] (-2,-pi) -- (3,-pi);
					\draw (0,pi) node {$\bullet$} node[below right] {$i\pi$};
					\draw (0,-pi) node {$\bullet$} node[above right] {$-i\pi$};
					\draw (-0.8,0) node {$\bullet$} node[below right] {$-\eta$};
					\draw (0,0) node {$\times$} circle (1.5);
					\draw (1.5,0) node[above right] {$B_\varepsilon(0)$};
					
					\draw[dartmouthgreen,thick] (-0.8,{-sqrt(1.5^2-0.8^2)}) -- (-0.8,{sqrt(1.5^2-0.8^2)}) node[near end, sloped] {$>$}  node[near start, sloped] {$>$};
					\draw[dartmouthgreen] (-1.7,{(sqrt(1.5^2-0.5^2)+sqrt((3/4+0.5+1/16)/3))/2}) node {$\Gamma_{k,in}$};
					
					\draw[dartmouthgreen,thick] (0,{-sqrt(1.5^2-0.8^2)}) -- (0,{sqrt(1.5^2-0.8^2)}) node[near end, sloped] {$>$}  node[near start, sloped] {$>$};
					\draw[dartmouthgreen] (0.5,{1}) node {$\Gamma_{k,in}^0$};
					
					\draw[blue,thick] (0,{-sqrt(1.5^2-0.8^2)}) -- (-0.8,{-sqrt(1.5^2-0.8^2)}) node[midway, sloped, rotate=180] {$>$};
					\draw[blue] (-0.5,-{1.8}) node {$\Gamma_{comp}^-$};
					
					\draw[blue,thick] (-0.8,{sqrt(1.5^2-0.8^2)}) -- (0,{sqrt(1.5^2-0.8^2)}) node[midway, sloped] {$>$};
					\draw[blue] (-0.5,{1.8}) node {$\Gamma_{comp}^+$};
				\end{tikzpicture}
				\caption{A representation of the path $\Gamma_{k,in}$,  $\Gamma_{k,in}^0$ and $\Gamma^\pm_{comp}$ for $\ut_k=0$ used in Lemma \ref{lemDécal}.}
				\label{figDecal}
			\end{center}
		\end{figure}
	
		Cauchy's formula then implies that
		\begin{multline*}
			\left|\int_{\Gamma^0_{k, in}}(\tau-\ut_k)^me^{n\tau}e^{j\varphi_k(\tau)}d\tau -\int_{\Gamma_{k, in}}(\tau-\ut_k)^me^{n\tau}e^{j\varphi_k(\tau)}d\tau \right|\\\leq\left|\int_{\Gamma^+_{comp}}(\tau-\ut_k)^me^{n\tau}e^{j\varphi_k(\tau)}d\tau\right|+\left|\int_{\Gamma^-_{comp}}(\tau-\ut_k)^me^{n\tau}e^{j\varphi_k(\tau)}d\tau\right|.
		\end{multline*}
		We need to find estimates for the two terms on the right-hand side. Both terms will be bounded similarly so we will focus on the first one. Since $\Gamma^+_{comp}\subset B_\varepsilon(\ut_k)$, we have
		$$\left|\int_{\Gamma^+_{comp}}(\tau-\ut_k)^me^{n\tau}e^{j\varphi_k(\tau)}d\tau\right|\lesssim \int_{-\eta}^0\exp\left(nt+j\Re\left(\varphi_k\left(t+i(\theta_k+r_\varepsilon(\eta))\right)\right)\right)dt.$$
		For $t\in ]-\eta, 0[$, since $t+i(\theta_k+r_\varepsilon(\eta))\in B_{\varepsilon_\star}(\ut_k)$ and $\frac{j}{\alpha_k}>0$, using the inequality \eqref{estVarphi}, we prove
		$$nt+j\Re(\varphi_k(t+i(\theta_k+r_\varepsilon(\eta))))\leq \frac{j}{\alpha_k}\left(\eta+A_R\eta^{2\mu_k}-A_Ir_\varepsilon(\eta)^{2\mu_k}\right).$$
		Using the inequality \eqref{inEta}, we have that $r_\varepsilon(\eta)\geq r_\varepsilon(\eta_\varepsilon)>\frac{\varepsilon}{2}$. Inequality \eqref{inEta2} then implies that
		$$\eta+A_R\eta^{2\mu_k}-A_Ir_\varepsilon(\eta)^{2\mu_k}\leq \eta_\varepsilon+A_R\eta_\varepsilon^{2\mu_k}-A_I\left(\frac{\varepsilon}{2}\right)^{2\mu_k}<0.$$
		Since $(n,j)\in \Dc_k$, we have that $\frac{j}{\alpha_k}\geq \frac{\tilde{\delta}_k}{\alpha_k} n$ so there must exist $c>0$ such that
		$$\forall (n,j)\in \Dc_k, \forall t \in ]-\eta, 0[, \quad nt+j\Re(\varphi_k(t+i(\theta_k+r_\varepsilon(\eta))))\leq -cn.$$
		This concludes the proof of Lemma \ref{lemDécal}.
	\end{proof}
	
	Using Lemma \ref{lemEstIntermMVC} and the estimate \eqref{estInterm}, we have thus proved that for all $s\in \N^*$, there exist two positive constants $C,c$ such that for all $(n,j)\in\Dc_k$
	\begin{multline}\label{inGinterm}
		 \left|\Gcc_j^n -\frac{\uz_k^n\uk_k^j}{2\pi} \int_{-r_\varepsilon(\eta)}^{r_\varepsilon(\eta)}P_{s,k}(it+\ut_k)\left(\sum_{l=0}^{s-1}\frac{(j R_{s,k}(it+\ut_k))^l}{l!}\right)\exp\left(it\left(n-\frac{j}{\alpha_k}\right) -\frac{j}{\alpha_k}\frac{\beta_k}{\alpha_k^{2\mu_k}}t^{2\mu_k}\right) dt\right|\\ \leq \frac{C}{n^\frac{s+1}{2\mu_k}}\exp\left(-c\left(\frac{|n\alpha_k-j|}{n^\frac{1}{2\mu_k}}\right)^\frac{2\mu_k}{2\mu_k-1}\right).
	\end{multline}
	
	There just remains to prove the following lemma to conclude the proof of Lemma \ref{lemEstInterm}.
	
	\begin{lemma}\label{lemExt}
		For all $m\in\N$ and $c_0>0$, there exist two positive constants $C,c>0$ such that
		$$\forall (n,j)\in \Dc_k, \quad \int_{r_\varepsilon(\eta)}^{+\infty} t^m\exp\left(-\frac{j}{\alpha_k}c_0t^{2\mu_k}\right)dt\leq Ce^{-cn}. $$
	\end{lemma}
	
	\begin{proof}
		The proof is done recursively and using the following equality proved by integrating by parts
		\begin{multline}\label{lemExt_egRecu}
			\int_{r_\varepsilon(\eta)}^{+\infty} t^m\exp\left(-\frac{j}{\alpha_k}c_0t^{2\mu_k}\right)dt = \frac{r_\varepsilon(\eta)^{m+1-2\mu_k}}{2\mu_kc_0\frac{j}{\alpha_k}}\exp\left(-c_0r_\varepsilon(\eta)^{2\mu_k}\frac{j}{\alpha_k}\right)\\+\frac{m+1-2\mu_k}{2\mu_kc_0\frac{j}{\alpha_k}}\int_{r_\varepsilon(\eta)}^{+\infty} t^{m-2\mu_k}\exp\left(-\frac{j}{\alpha_k}c_0t^{2\mu_k}\right)dt.
		\end{multline}
		
		\begin{itemize}
			\item For $m\in\lc 0, \ppp, 2\mu_k-1\rc$, since the second term of the sum on the right hand side of \eqref{lemExt_egRecu} is non-positive, using the fact that $(n,j)\in\Dc_k$, we directly prove the result.
			\item If we consider $\tilde{m}\geq 2\mu_k$ such that the result of lemma has been proved for all $m\in \lc0, \ppp, \tilde{m}-1\rc$, then the equality \eqref{lemExt_egRecu} implies the result for $m=\tilde{m}$.
		\end{itemize}
	\end{proof}
	
	Combining Lemmas \ref{lemExt}, \ref{est_expn} and the inequality \eqref{inGinterm}, we easily conclude the proof of Lemma \ref{lemEstInterm}.

	\subsection{Step 2 : Proof of Lemma \ref{lemDevAsymp}}\label{secAsymp}
	
	As we explained in Section \ref{subsecPlan}, Lemma \ref{lemEstInterm} and the equality \eqref{egHjalp} imply that we proved generalized Gaussian estimates on the difference between the elements $\Gcc_j^n$ and a linear combination of 
	$$\frac{1}{\left(\frac{j}{\alpha_k}\right)^\frac{l}{2\mu_k}}{H_{2\mu_k}^{\beta_k}}^{(m)}\left(\frac{n\alpha_k-j}{\left(\frac{j}{\alpha_k}\right)^\frac{1}{2\mu_k}}\right)\quad \text{where }l\in\N^*, m\in\N.$$
	We now need to approach the above terms by the elements appearing in Theorem \ref{thPrinc}, i.e. a linear combination of 
	$$\frac{1}{n^\frac{l}{2\mu_k}}\left(\frac{n\alpha_k-j}{n^\frac{1}{2\mu_k}}\right)^{m_2}{H_{2\mu_k}^{\beta_k}}^{(m_1)}\left(\frac{n\alpha_k-j}{n^\frac{1}{2\mu_k}}\right)\quad \text{where }l\in\N^*, m_1,m_2\in\N.$$
	
	This is the goal of Lemma \ref{lemDevAsymp} that we recall here:
	\begin{lemma*}[Lemma \ref{lemDevAsymp}]
		For all $s \in \N$, $m\in \N$, $l\in\N\backslash\lc0\rc$ and $k\in \lc 1, \ppp,K\rc$, if we consider $d\in \N$ such that
		$$d\geq \frac{s+1}{2\mu_k-1}$$
		then there exist two constants $C,c>0$ such that for all $(n,j)\in \Dc_k$,
		$$\left|\frac{{H_{2\mu_k}^{\beta_k}}^{(m)}\left(Y_{n,j,k}\right)}{\left(\frac{j}{\alpha_k}\right)^\frac{l}{2\mu_k}}- \sum_{k_1=0}^{d-1}\sum_{k_3=0}^{d-1} \frac{\Bcc^k_{l,k_1,k_3}}{n^\frac{l+(2\mu_k-1)k_3}{2\mu_k}}\left(X_{n,j,k}\right)^{k_1+k_3}{H_{2\mu_k}^{\beta_k}}^{(m+k_1)}\left(X_{n,j,k}\right) \right|  \leq \frac{C}{n^\frac{s+1}{2\mu_k}}\exp\left(-c\left|X_{n,j,k}\right|^\frac{2\mu_k}{2\mu_k-1}\right)$$
		where $Y_{n,j,k}:= \frac{n\alpha_k-j}{\left(\frac{j}{\alpha_k}\right)^\frac{1}{2\mu_k}}$, $X_{n,j,k}:=\frac{n\alpha_k-j}{n^\frac{1}{2\mu_k}}$ and 
		$$\Bcc^k_{l,k_1,k_3}:= \sum_{k_2=0}^{k_1}\frac{\binom{k_1}{k_2}(-1)^{k_1-k_2}}{k_1!k_3!\alpha_k^{k_3}}\left(\prod_{k_4=0}^{k_3-1} \frac{l+k_2}{2\mu_k}+k_4\right)  .$$
	\end{lemma*}

	First, we prove the following lemma.
	
	\begin{lemma}\label{lemDevH}
		For all $s \in \N$, $m\in \N$ and $k\in \lc 1, \ppp,K\rc$, if we consider $d\in \N$ such that
		$$d\geq \frac{s+1}{2\mu_k-1}$$
		then there exist two constants $C,c>0$ such that for all $(n,j)\in \Dc_k$,
		\begin{multline*}
		\left|{H_{2\mu_k}^{\beta_k}}^{(m)}\left(Y_{n,j,k}\right)- \sum_{k_1=0}^{d-1}\frac{{H_{2\mu_k}^{\beta_k}}^{(m+k_1)}\left(X_{n,j,k}\right)}{k_1!}(n\alpha_k-j)^{k_1}\left(\left(\frac{\alpha_k}{j}\right)^\frac{1}{2\mu_k}-\left(\frac{1}{n}\right)^\frac{1}{2\mu_k}\right)^{k_1}\right| \\ \leq \frac{C}{n^\frac{s+1}{2\mu_k}}\exp\left(-c\left|X_{n,j,k}\right|^\frac{2\mu_k}{2\mu_k-1}\right).
		\end{multline*}
	\end{lemma}
	
	\begin{proof}
		We will apply Taylor's Theorem to bound the term on the left hand side of the inequality. We observe using the bounds of Lemma \ref{ineg_H}  on the derivatives of $H_{2\mu_k}^{\beta_k}$ that there exist two positive constants $C,c$ such that
		\begin{equation}\label{lemDevHin1}
			\forall (n,j)\in \Dc_k, \forall x\in \left[X_{n,j,k},Y_{n,j,k}\right], \quad \left|{H_{2\mu_k}^{\beta_k}}^{(m+d)}(x)\right|\leq C\exp\left(-c\left|X_{n,j,k}\right|^\frac{2\mu_k}{2\mu_k-1}\right).
		\end{equation}
		We also observe that the mean value inequality implies that there exists a constant $C>0$ such that
		\begin{equation}\label{lemDevHin2}
			\forall (n,j)\in\Dc_k, \quad \left|\left(\frac{\alpha_k}{j}\right)^\frac{1}{2\mu_k}-\left(\frac{1}{n}\right)^\frac{1}{2\mu_k}\right|\leq \frac{C}{n^{1+\frac{1}{2\mu_k}}}|n\alpha_k-j|.
		\end{equation}
		Combining Taylor's Theorem and both inequalities \eqref{lemDevHin1} and \eqref{lemDevHin2}, we can prove the existence of two positive constants $C,c$ such that for all $(n,j)\in \Dc_k$
		\begin{multline*}
			\left|{H_{2\mu_k}^{\beta_k}}^{(m)}\left(Y_{n,j,k}\right)- \sum_{k_1=0}^{d-1}\frac{{H_{2\mu_k}^{\beta_k}}^{(m+k_1)}\left(X_{n,j,k}\right)}{k_1!}(n\alpha_k-j)^{k_1}\left(\left(\frac{\alpha_k}{j}\right)^\frac{1}{2\mu_k}-\left(\frac{1}{n}\right)^\frac{1}{2\mu_k}\right)^{k_1}\right|  \\ \leq\frac{C}{n^{d\left(1-\frac{1}{2\mu_k}\right)}}\left|X_{n,j,k}\right|^{2d}\exp\left(-c\left|X_{n,j,k}\right|^\frac{2\mu_k}{2\mu_k-1}\right).
		\end{multline*}
		Since the function $x\mapsto x^{2d}\exp\left(-\frac{c}{2}x^\frac{2\mu_k}{2\mu_k-1}\right)$ is bounded, our choice for $d$ allows us to conclude.
	\end{proof}
	
	Using Lemma \ref{lemDevH}, we have now approached the elements $\Gcc_j^n$ via a linear combination of 
	\begin{equation}\label{termes}
		\frac{(n\alpha_k-j)^{k_1}}{n^\frac{k_1-k_2}{2\mu_k}\left(\frac{j}{\alpha_k}\right)^\frac{l+k_2}{2\mu_k}}{H_{2\mu_k}^{\beta_k}}^{(m+k_1)}\left(\frac{n\alpha_k-j}{n^\frac{1}{2\mu_k}}\right)\quad \text{where }l\in\N^*, m\in\N, k_1\in\N, k_2\in\lc0,\ppp,k_1\rc.
	\end{equation}
	We approach the terms in \eqref{termes} using the following lemma.
	
	\begin{lemma}\label{lemDevJAlp}
		We consider $s \in \N$, $m\in \N$, $l\in \N\backslash\lc0\rc$, $k_1\in \N$, $k_2\in\lc0, \ppp,k_1\rc$ and $k\in \lc 1, \ppp,K\rc$. We define the function
		$$\Psi_q:x\in \R_+^*\rightarrow \frac{1}{x^q}.$$
		If we consider $d\in \N$ such that
		$$d\geq \frac{s+1}{2\mu_k-1}$$
		then there exist two constants $C,c>0$ such that for all $(n,j)\in \Dc_k$,
		\begin{multline*}
			\left|{H_{2\mu_k}^{\beta_k}}^{(m+k_1)}\left(X_{n,j,k}\right) \frac{(n\alpha_k-j)^{k_1}}{n^\frac{k_1-k_2}{2\mu_k}}\left(\Psi_{\frac{l+k_2}{2\mu_k}}\left(\frac{j}{\alpha_k}\right)-\sum_{k_3=0}^{d-1}\frac{\Psi_{\frac{l+k_2}{2\mu_k}}^{(k_3)}\left(n\right)}{k_3!}\left(\frac{j}{\alpha_k}-n\right)^{k_3}\right)\right| \\ \leq \frac{C}{n^\frac{s+1}{2\mu_k}}\exp\left(-c\left|X_{n,j,k}\right|^\frac{2\mu_k}{2\mu_k-1}\right)
		\end{multline*}
		with $X_{n,j,k}:=\frac{n\alpha_k-j}{n^\frac{1}{2\mu_k}}$.
	\end{lemma}
	
	\begin{proof}
		We will apply Taylor's theorem to bound the term on the left hand side of the inequality. We observe that there exist two positive constants $C,c$ such that
		\begin{equation}\label{lemDevJAlpIn}
			\forall (n,j)\in \Dc_k, \forall x\in \left[n,\frac{j}{\alpha_k}\right], \quad \left|\Psi_{\frac{l+k_2}{2\mu_k}}^{(d)}(x)\right|\leq \frac{C}{n^{\frac{l+k_2}{2\mu_k}+d}}.
		\end{equation}
		Thus, the inequality \eqref{lemDevJAlpIn} and Taylor's theorem imply the existence of two positive constants $C,c$ such that for all $(n,j)\in \Dc_k$
		\begin{multline*}
			\left|{H_{2\mu_k}^{\beta_k}}^{(m+k_1)}\left(X_{n,j,k}\right) \frac{(n\alpha_k-j)^{k_1}}{n^\frac{k_1-k_2}{2\mu_k}}\left(\Psi_{\frac{l+k_2}{2\mu_k}}\left(\frac{j}{\alpha_k}\right)-\sum_{k_3=0}^{d-1}\frac{\Psi_{\frac{l+k_2}{2\mu_k}}^{(k_3)}\left(n\right)}{k_3!}\left(\frac{j}{\alpha_k}-n\right)^{k_3}\right)\right| \\ \leq\frac{C}{n^{\frac{l}{2\mu_k}+d\left(1-\frac{1}{2\mu_k}\right)}}\left|X_{n,j,k}\right|^{k_1+d}\exp\left(-c\left|X_{n,j,k}\right|^\frac{2\mu_k}{2\mu_k-1}\right).
		\end{multline*}
		Since the function $x\mapsto x^{k_1+d}\exp\left(-\frac{c}{2}x^\frac{2\mu_k}{2\mu_k-1}\right)$ is bounded, our choice for $d$ allows us to conclude.
	\end{proof}
	
	Lemmas \ref{lemDevH} and \ref{lemDevJAlp} allow us to conclude the proof of Lemma \ref{lemDevAsymp}. 
	
	\subsection{Step 3: Construction of the polynomials $\Pcc^k_\sigma$ satisfying Lemma \ref{lemPrinc} and Theorem \ref{thPrinc}}\label{subsecLemImpTh}
	
	Now that Lemmas \ref{lemEstInterm} and \ref{lemDevAsymp} are proved, we will construct the polynomials $\Pcc^k_\sigma$ in $\C[X,Y]$ which will verify Lemma \ref{lemPrinc} and Theorem \ref{thPrinc}. We start by introducing some notations. We fix $k\in \lc1, \ppp,K\rc$ and $s_k\in\N$. For $l\in \lc0,\ppp,s_k-1\rc$, we define the coefficients $\Acc^k_{s_k,l,m}\in \C$ for $m\in \lc (2\mu_k+1)l, \ppp,(2\mu_k+s_k-1)l+s_k-1\rc$ such that
	\begin{equation}\label{defC}
		\forall \tau\in B_{\varepsilon_\star}(\ut_k), \quad P_{s_k,k}(\tau) \frac{R_{s_k,k}(\tau)^l}{l!} = \sum_{m= (2\mu_k+1)l}^{(2\mu_k+s_k-1)l+s_k-1}\Acc^k_{s_k,l,m}(\tau-\ut_k)^m.
	\end{equation}
	where the polynomial functions $P_{s_k, k}$ and $R_{s_k, k}$ are defined in Lemma \ref{lem_varpi}. Using Lemma \ref{lemEstInterm} and equality \eqref{egHjalp}, we prove that there exist two positive constants $C,c$ such that for all $(n,j)\in\Dc_k$
	\begin{equation}\label{step3Est1}
		\left|\Gcc_j^n-\uz_k^n\uk_k^j\sum_{l=0}^{s_k-1}\sum_{m=(2\mu_k+1)l}^{(2\mu_k+s_k-1)l+s_k-1}\frac{\Acc^k_{s_k,l,m}\alpha_k^{m+l}|\alpha_k|}{\left(\frac{j}{\alpha_k}\right)^\frac{m-2\mu_kl+1}{2\mu_k}}{H_{2\mu_k}^{\beta_k}}^{(m)}\left(Y_{n,j,k}\right)\right|\leq \frac{C}{n^\frac{s_k+1}{2\mu_k}}\exp\left(-c\left(\frac{|n\alpha_k-j|}{n^\frac{1}{2\mu_k}}\right)^\frac{2\mu_k}{2\mu_k-1}\right)
	\end{equation}
	where $Y_{n,j,k}:=\frac{n\alpha_k-j}{\left(\frac{j}{\alpha_k}\right)^\frac{1}{2\mu_k}}$.
	
	We now want to apply Lemma \ref{lemDevAsymp}, so we need to define an integer $d\in \N$ such that
	$$d \geq \frac{s_k+1}{2\mu_k-1}.$$
	We will consider that $d=s_k+1$ so that when we will do computations of the polynomials $\Pcc^k_\sigma$ in Section \ref{secComput}, we will not have to distinguish the value of $d$ depending on the value of $\mu_k$. Then, for $l\in \lc0,\ppp,s_k-1\rc$, $m\in \lc (2\mu_k+1)l, \ppp,(2\mu_k+s_k-1)l+s_k-1\rc$ and $k_1,k_3\in \lc0,\ppp,s_k\rc$, we define the coefficients
	\begin{align}
		\begin{split}
		\Ccc^k_{s_k,l,m,k_1,k_3}&:= \Acc^k_{s_k,l,m}\alpha_k^{m+l}|\alpha_k|\Bcc^k_{m-2\mu_kl+1,k_1,k_3} \\
		&=  \frac{\Acc^k_{s_k,l,m}\alpha_k^{m+l-k_3}|\alpha_k|}{k_1!k_3!}\sum_{k_2=0}^{k_1}\binom{k_1}{k_2}(-1)^{k_1-k_2}\left(\prod_{k_4=0}^{k_3-1} \frac{m-2\mu_kl+1+k_2}{2\mu_k}+k_4\right) 
		\end{split}\label{defCcc}
	\end{align}
	where the coefficients $\Bcc^k_{m-2\mu_kl+1,k_1,k_3}$ are defined in Lemma \ref{lemDevAsymp}. Combining the result of Lemma \ref{lemDevAsymp} with the estimates \eqref{step3Est1}, we prove the existence of two positive constants $C,c$ such that for all $(n,j)\in\Dc_k$
	\begin{multline}\label{step3Est2}
		\left|\Gcc_j^n-\uz_k^n\uk_k^j\sum_{l=0}^{s_k-1}\sum_{m=(2\mu_k+1)l}^{(2\mu_k+s_k-1)l+s_k-1}\sum_{k_1=0}^{s_k}\sum_{k_3=0}^{s_k}\frac{\Ccc^k_{s_k,l,m,k_1,k_3}}{n^\frac{m-2\mu_kl+k_3(2\mu_k-1)+1}{2\mu_k}}{X_{n,j,k}}^{k_1+k_3}{H_{2\mu_k}^{\beta_k}}^{(m+k_1)}\left(X_{n,j,k}\right)\right|\\ \leq \frac{C}{n^\frac{s_k+1}{2\mu_k}}\exp\left(-c\left|X_{n,j,k}\right|^\frac{2\mu_k}{2\mu_k-1}\right)
	\end{multline}
	with $X_{n,j,k}:=\frac{n\alpha_k-j}{n^\frac{1}{2\mu_k}}$. For $\sigma \in \lc 1 ,\ppp,s_k\rc$, we define the polynomial
	\begin{equation}\label{Pcc}
		\Pcc^k_\sigma(X,Y):= \sum_{l=0}^{s_k-1}\sum_{m=(2\mu_k+1)l}^{(2\mu_k+s_k-1)l+s_k-1}\sum_{k_1=0}^{s_k}\sum_{k_3=0}^{s_k} \ind_{m-2\mu_kl+k_3(2\mu_k-1)+1=\sigma}\Ccc^k_{s_k,l,m,k_1,k_3}X^{k_1+k_3}Y^{m+k_1}\in \C[X,Y].
	\end{equation}
	Using the estimates on the derivatives of $H_{2\mu}^\beta$ (Lemma \ref{ineg_H}) to take care of the terms where $m-2\mu_kl+k_3(2\mu_k-1)+1\geq s_k+1$, the inequality \eqref{step3Est2} implies that the polynomials $\Pcc^k_\sigma$ verify the estimates \eqref{inPrinc} of Lemma \ref{lemPrinc}. Lemma \ref{lemPrinc} is proved and Theorem \ref{thPrinc} in the case where the elements $\alpha_k$ are supposed to be distinct ensues from Lemma \ref{lemPrinc} and inequalities \eqref{inHorsD} and \eqref{inHorsDk}.
	
	\section{Closing arguments on Theorem \ref{thPrinc} and proof of Corollary \ref{cor_prin}}\label{secConcluTh}
	
	\subsection{Proof of Theorem \ref{thPrinc} when the elements $\alpha_k$ can be equal}\label{sec_alpha_eg}
	
	As we said in the beginning on Section \ref{sec_GT}, we supposed in the proof that the elements $\alpha_k$ were distinct from one another. In the case where the $\alpha_k$ can be equal, there are some changes that need to be done but the calculations remain similar. Most modifications will happen on the part of the proof contained in Section \ref{sec_GS/GT}. 
	
	First, just as in Section \ref{subsec_est_far}, we would define $\sd_k$, $\ud_k$ and $\Dc_k$ in the same manner but with the added condition that if $\alpha_k=\alpha_l$, then $\sd_k=\sd_l$ and $\ud_k=\ud_l$.
	
	If we consider $k_0\in\lc1,\ppp,K\rc$, we define 
	$$\Jc_{k_0}:= \lc k\in\lc1,\ppp,K\rc, \quad \alpha_k=\alpha_{k_0}\rc.$$
	We observe that for $k\in\Jc_{k_0}$, we have $\Dc_k= \Dc_{k_0}$ because of our new condition.
	
	Lemmas \ref{prop_loin_G} and \ref{prop_loin_H} remain true. The inequality \eqref{inHorsD} thus remains true, however inequality \eqref{inHorsDk} now becomes that for $k_0\in\lc1,\ppp,K\rc$, there exist two constants $C,c>0$ such that for all $(n,j)\in\Dc_{k_0}$
	\begin{multline*}
		\left|\Gcc_j^n -\sum_{k=1}^K\sum_{\sigma=1}^{s_k}\frac{\uz_k^n\uk_k^j}{n^{\frac{\sigma}{2\mu_k}}} \left(\Pcc^k_\sigma\left(X_{n,j,k},\frac{d}{dx}\right)H_{2\mu_k}^{\beta_k}\right)\left(X_{n,j,k}\right)\right|  \leq \sum_{\underset{k\notin \Jc_{k_0}}{k=1}}^K\frac{C}{n^\frac{s_k+1}{2\mu_k}}\exp\left(-c\left|X_{n,j,k}\right|^\frac{2\mu_k}{2\mu_k-1}\right) \\ +\left|\Gcc_j^n -\sum_{k\in \Jc_{k_0}}\sum_{\sigma=1}^{s_{k}}\frac{\uz_{k}^n\uk_{k}^j}{n^{\frac{\sigma}{2\mu_{k}}}} \left(\Pcc^{k}_{\sigma}\left(X_{n,j,k},\frac{d}{dx}\right)H_{2\mu_{k}}^{\beta_{k}}\right)\left(X_{n,j,k}\right)\right|.
	\end{multline*}	
	Therefore, to prove Theorem \ref{thPrinc}, we now have to prove the following lemma which is a modification of Lemma \ref{lemPrinc}.
	\begin{lemma}[Modified Lemma \ref{lemPrinc}]\label{lemPrincNew}
		For all $k_0\in\lc1,\ppp,K\rc$ and $(s_k)_{k\in\Jc_{k_0}}\in \N^{\Jc_{k_0}}$, there exist a family of polynomials $(\Pcc^k_\sigma)_{\sigma\in\lc 1,\ppp,s_k\rc}$ in $\C[X,Y]$ for each $k\in\Jc_{k_0}$ and two positive constants $C,c$ such that for $(n,j)\in\Dc_{k_0}$
		$$\left|\Gcc_j^n -\sum_{k\in\Jc_{k_0}}\sum_{\sigma=1}^{s_k}\frac{\uz_k^n\uk_k^j}{n^{\frac{\sigma}{2\mu_k}}} \left(\Pcc^k_\sigma\left(X_{n,j,k},\frac{d}{dx}\right)H_{2\mu_k}^{\beta_k}\right)\left(X_{n,j,k}\right)\right|\leq \sum_{k\in \Jc_{k_0}} \frac{C}{n^\frac{s_k+1}{2\mu_k}}\exp\left(-c\left|X_{n,j,k}\right|^\frac{2\mu_k}{2\mu_k-1}\right) $$
		with $X_{n,j,k}:=\frac{n\alpha_k-j}{n^\frac{1}{2\mu_k}}$. 
	\end{lemma}
	
	Just as in the case where the elements $\alpha_k$ were supposed distinct, if Lemma \ref{lemPrincNew} is verified, then the families of polynomials $(\Pcc^k_\sigma)_{k,\sigma}$ constructed in Lemma \ref{lemPrincNew} will also verify the estimates \eqref{devPrinc} of Theorem \ref{thPrinc}. Since the equality \eqref{egHjalp} and Lemma \ref{lemDevAsymp} remain true, to prove Lemma \ref{lemPrincNew}, we only have to prove the following Lemma which is a modification of Lemma \ref{lemEstInterm}.
	\begin{lemma}[Modified Lemma \ref{lemEstInterm}]\label{lemEstIntermNew}
		For all $k_0\in\lc1,\ppp,K\rc$ and for all $(s_k)_{k\in\Jc_{k_0}}\in {\N^*}^{\Jc_{k_0}}$, there exist two positive constants $C,c$ such that for all $(n,j)\in \Dc_{k_0}$
		\begin{multline*}
			\left|\Gcc_j^n -\sum_{k\in\Jc_{k_0}}\frac{\uz_k^n\uk_k^j}{2\pi} \int_{-\infty}^{+\infty}P_{s_k,k}(it+\ut_k)\left(\sum_{l=0}^{s_k-1}\frac{(j R_{s_k,k}(it+\ut_k))^l}{l!}\right)\exp\left(it\left(n-\frac{j}{\alpha_k}\right) -\frac{j}{\alpha_k}\frac{\beta_k}{\alpha_k^{2\mu_k}}t^{2\mu_k}\right) dt\right|\\ \leq \sum_{k\in \Jc_{k_0}}\frac{C}{n^\frac{s_k+1}{2\mu_k}}\exp\left(-c\left(\frac{|n\alpha_k-j|}{n^\frac{1}{2\mu_k}}\right)^\frac{2\mu_k}{2\mu_k-1}\right)
		\end{multline*} 
		where the polynomial functions $P_{s_k,k}$ and $R_{s_k,k}$ have explicit expression defined in Lemma \ref{lem_varpi}.
	\end{lemma}
	
	Therefore, there just remains to prove Lemma \ref{lemEstIntermNew} and Theorem \ref{thPrinc} will ensue. We recall that, to prove Lemma \ref{lemEstInterm} in the case where the elements $\alpha_k$ were distinct from one another, we found an expression of the elements $\Gcc_j^n$ as an integral along the path $\Gamma_k$ 
	$$\forall n\in\N^*,\forall j\in\Z, \quad \Gcc_j^n=\frac{1}{2i\pi} \int_{\Gamma_k} e^{n\tau}\Gg_j(\tau)d\tau$$
	and used the triangular inequality to find the inequality \eqref{egDecoup} that we recall here
	$$\left|\Gcc_j^n -\frac{1}{2i\pi} \int_{\Gamma_{k,in}}P_{s,k}(\tau)\left(\sum_{l=0}^{s-1}\frac{(j R_{s,k}(\tau))^l}{l!}\right)e^{n\tau}e^{j\varphi_k(\tau)}d\tau\right| \leq \frac{1}{2\pi}\sum_{l=1}^8 E_l.$$
	We then bounded all the terms $E_i$ to find an estimate on
	 $$\left|\Gcc_j^n -\frac{1}{2i\pi} \int_{\Gamma_{k,in}}P_{s,k}(\tau)\left(\sum_{l=0}^{s-1}\frac{(j R_{s,k}(\tau))^l}{l!}\right)e^{n\tau}e^{j\varphi_k(\tau)}d\tau\right|.$$ 
	 In the case where the elements $\alpha_k$ are no longer supposed to be distinct, the reasoning is the same but with a better suited choice of path to express the elements $\Gcc_j^n$. We fix $k_0\in\lc1,\ppp,K\rc$ and introduce the path $\widetilde{\Gamma}_{k_0}$ which is the ray $\lc -\eta +it, t\in[-\pi,\pi]\rc$ deformed into the path $\Gamma_{k,in}$ inside the balls $B_\varepsilon(\ut_k)$ for $k\in\Jc_{k_0}$ (see Figure \ref{chem_3}). Using Cauchy's formula and taking into account the "$2i\pi$-periodicity" of $\Gg_j(\tau)$, we have that 
	 $$\forall n\in\N^*,\forall j\in\Z, \quad \Gcc_j^n=\frac{1}{2i\pi} \int_{\widetilde{\Gamma}_{k_0}} e^{n\tau}\Gg_j(\tau)d\tau.$$
	 
	 \begin{figure}
	 	\begin{center}
	 		\begin{tikzpicture}[scale=1]
	 			\draw[->] (-2,0) -- (3,0) node[right] {$\Re(\tau)$};
	 			\draw[->] (0,-3.4) -- (0,3.4) node[above] {$\Im(\tau)$};
	 			\draw[dashed] (-2,pi) -- (3,pi);
	 			\draw[dashed] (-2,-pi) -- (3,-pi);
	 			\draw (0,pi) node {$\bullet$} node[below right] {$i\pi$};
	 			\draw (0,-pi) node {$\bullet$} node[above right] {$-i\pi$};

	 			\draw (0,-2) node {$\times$} circle (0.75);
	 			\draw (0.75,-2) node[above right] {$B_\varepsilon(\ut_{k_1})$};
	 			
	 			\draw (0,-0.2) node {$\times$} circle (0.75);
	 			\draw (0.75,-0.2) node[below right] {$B_\varepsilon(\ut_l)$};
	 			
	 			\draw (0,1.9) node {$\times$} circle (0.75);
	 			\draw (0.75,1.9) node[above right] {$B_\varepsilon(\ut_{k_2})$};
	 			
	 			\draw[red,thick] (-0.5,-pi) -- (-0.5,{-2-sqrt(0.75^2-0.5^2)});
	 			\draw[red,thick] (-0.5,{-2+sqrt(0.75^2-0.5^2)}) -- (-0.5,{-0.2-sqrt(0.75^2-0.5^2)});
	 			\draw[red,thick] (-0.5,{-0.2+sqrt(0.75^2-0.5^2)}) -- (-0.5,{1.9-sqrt(0.75^2-0.5^2)});
	 			\draw[red,thick] (-0.5,{1.9+sqrt(0.75^2-0.5^2)}) -- (-0.5,pi);
	 			
	 			\draw[thick,blue] plot [samples = 100, domain={-sqrt(0.2)}:{sqrt(0.2)}] ({0.3-4*abs(\x)^2},\x-2) ;
	 			\draw[blue,thick] (-0.5,{-2-sqrt(0.75^2-0.5^2)}) -- (-0.5,{-2-sqrt(0.2)});
	 			\draw[blue,thick] (-0.5,{-2+sqrt(0.2)}) -- (-0.5,{-2+sqrt(0.75^2-0.5^2)});
	 			\draw[blue,thick,dotted] (-0.5,{-2-sqrt(0.2)}) -- (-0.5,{-2+sqrt(0.2)});

	 			\draw[thick,blue] plot [samples = 100, domain={-sqrt(0.15)}:{sqrt(0.15)}] ({-0.2-2*abs(\x)^2},\x+1.9) ;
	 			\draw[blue,thick] (-0.5,{1.9-sqrt(0.75^2-0.5^2)}) -- (-0.5,{1.9-sqrt(0.15)});
	 			\draw[blue,thick] (-0.5,{1.9+sqrt(0.15)}) -- (-0.5,{1.9+sqrt(0.75^2-0.5^2)});
	 			\draw[blue,thick,dotted] (-0.5,{1.9-sqrt(0.15)}) -- (-0.5,{1.9+sqrt(0.15)});
	 			
	 			\draw[red,thick] (-0.5,{-0.2-sqrt(0.75^2-0.5^2)}) -- (-0.5,{-0.2+sqrt(0.75^2-0.5^2)});
	 			
	 			\draw[blue] (-0.75, 1.9) node[left] {$\Gamma_{k_2,in}$};
	 			\draw[blue] (-0.75, -2) node[left] {$\Gamma_{k_1,in}$};
	 			\draw[red] (-0.75, 0.8) node[left] {$\widetilde{\Gamma}_{k_0,out}$};
	 		\end{tikzpicture}
	 		\caption{A representation of the path $\widetilde{\Gamma}_{k_0}$. Inside the balls $B_\varepsilon(\ut_k)$ where $k$ belongs to $\Jc_{k_0}$, it follows the path $\Gamma_{k,in}$ composed of $\Gamma_{k,res}$ and $\Gamma_{k,p}$. For $l\in\lc1,\ppp,K\rc$, if there is no $k\in\Jc_{k_0}$ such that $\ut_k=\ut_l$, then the path $\widetilde{\Gamma}_{k_0}$ inside $B_\varepsilon(\ut_l)$ just corresponds to the ray $\lc-\eta+it, t\in[-\pi,\pi]\rc$.}
	 		\label{chem_3}
	 	\end{center}
	 \end{figure}
	 
	We end up with an inequality similar to \eqref{egDecoup}.
	\begin{equation}\label{egDecoupNew}
		\left|\Gcc_j^n -\frac{1}{2i\pi} \sum_{k\in\Jc_{k_0}}\int_{\Gamma_{k,in}}P_{s_k,k}(\tau)\left(\sum_{l=0}^{s_k-1}\frac{(j R_{s_k,k}(\tau))^l}{l!}\right)e^{n\tau}e^{j\varphi_k(\tau)}d\tau\right| \leq \frac{1}{2\pi}\left(E_{out}+\sum_{k\in \Jc_{k_0}}\sum_{l=2}^8 E_{l,k}\right)
	\end{equation}
	where $E_{l,k}$ has the same definition as $E_l$ in \eqref{egDecoup} but depends on the $k\in\Jc_{k_0}$ we consider. The term $E_{out}$ is similar to $E_1$ in \eqref{egDecoup} and is equal to 
	$$E_{out} =  \left|\int_{\widetilde{\Gamma}_{k_0,out}} e^{n\tau} \Gg_j(\tau) d\tau\right|,$$
	where $\widetilde{\Gamma}_{k_0,out}$ corresponds to the part of $\widetilde{\Gamma}_{k_0}$ outside the balls $B_\varepsilon(\ut_k)$ for $k\in\Jc_{k_0}$ (see the red path on Figure \ref{chem_3}). Reasoning in the same manner as in the case where the elements $\alpha_k$ are different from one another, we get estimates on the different terms. The minor modifications are left to the reader.  Notice that Lemmas \ref{lemDécal} and \ref{lemExt} are still verified, Lemma \ref{lemEstIntermNew} ensues. Therefore, Theorem \ref{thPrinc} in the case where the elements $\alpha_k$ can be equal is proved for the same polynomials $\Pcc^k_\sigma$ given in Section \ref{subsecLemImpTh}. 
	
	\subsection{Proof of Corollary \ref{cor_prin}}\label{secCor}
	
	We are now going to prove Corollary \ref{cor_prin} that we recall here:
	\begin{corollary*}[Corollary \ref{cor_prin}]
		Let $a\in \ell^1(\Z)$ which verifies Hypotheses \ref{H1} and \ref{H2_bis}. If there exists some integer $J\in\Z$ such that the sequence $b:=(a_{j+J})_{j\in\Z}$ verifies Hypotheses \ref{H2} and \ref{H3}, then for all $s_1,\ppp,s_K\in\N$ there exist a family of polynomials $(\Pcc^k_\sigma)_{\sigma\in\lc1, \ppp,s_k\rc}$ in $\C[X,Y]$ for each $k\in\lc1, \ppp,K\rc$ and two positive constants $C,c$ such that for all $n\in\N^*$ and $j\in \Z$
		$$\left|\Gcc_j^n-\sum_{k=1}^K\sum_{\sigma=1}^{s_k}\frac{\uz_k^n\uk_k^j}{n^\frac{\sigma}{2\mu_{k}}}\left(\Pcc^k_\sigma\left(X_{n,j,k},\frac{d}{dx}\right)H_{2\mu_k}^{\beta_k}\right)\left(X_{n,j,k}\right)\right|\leq\sum_{k=1}^K \frac{C}{n^\frac{s_k+1}{2\mu_k}}\exp\left(-c|X_{n,j,k}|^\frac{2\mu_k}{2\mu_k-1}\right)$$
		with $X_{n,j,k}=\frac{n\alpha_k-j}{n^\frac{1}{2\mu_k}}$.
	\end{corollary*}
	We consider that $a$ satisfies the hypotheses of Corollary \ref{cor_prin}. As we said just before we introduced the corollary, we observe that if we define $\widetilde{F}$ the symbol associated with $b$, then we have that 
	$$\forall \kappa\in\S^1, \quad \widetilde{F}(\kappa)=\kappa^{-J}F(\kappa).$$
	and we have for $k\in\lc1,\ppp,K\rc$
	\begin{equation}\label{G}
		\widetilde{F}(\uk_ke^{i\xi})\underset{\xi\rightarrow0}= \uk_k^{-J}\uz_k\exp(-i(\alpha_k+J) \xi - \beta_k \xi^{2\mu_k} + o(|\xi|^{2\mu_k})).
	\end{equation}
	We fix $s_1, \ppp,s_K\in\N$. Applying Theorem \ref{thPrinc} for the sequence $b$, there exist a family of polynomials $(\Pcc^k_\sigma)_{\sigma\in\lc1,\ppp,s_k\rc}$ in $\C[X,Y]$ for each $k\in\lc1,\ppp,K\rc$ and two positive constants $C,c$ such that for all $n\in\N^*$ and $j\in \Z$
	\begin{multline*}
		\left|\left(\Lcc_b^n\delta\right)_j-\sum_{k=1}^K\sum_{\sigma=1}^{s_k}\frac{\uz_k^n\uk_k^j}{n^\frac{\sigma}{2\mu_{k}}}\left(\Pcc^k_\sigma\left(\frac{n(\alpha_k+J)-j}{n^\frac{1}{2\mu_k}},\frac{d}{dx}\right)H_{2\mu_k}^{\beta_k}\right)\left(\frac{n(\alpha_k+J)-j}{n^\frac{1}{2\mu_k}}\right)\right| \\ \leq\sum_{k=1}^K \frac{C}{n^\frac{s_k+1}{2\mu_k}}\exp\left(-c\left(\frac{|n(\alpha_k+J)-j|}{n^\frac{1}{2\mu_k}}\right)^\frac{2\mu_k}{2\mu_k-1}\right).
	\end{multline*}
	By observing that 
	$$\forall n\in\N^*,\forall j\in\Z, \quad (\Lcc_b^n\delta)_j = (\Lcc_a^n\delta)_{j-nJ}= \Gcc_{j-nJ}^n,$$
	we conclude the proof of Corollary \ref{cor_prin}.	
	
	\section{Computations of the polynomials $\Pcc^k_\sigma$} \label{secComput}
	
	Now that Theorem \ref{thPrinc} is proved, we want to compute the polynomials $\Pcc^k_\sigma$ defined with \eqref{Pcc} in the proof of Theorem \ref{thPrinc}. We separate this section in three parts:
	\begin{itemize}
		\item The coefficients of the polynomials $\Pcc^k_\sigma$ depend on the elements $\Acc^k_{s_k,l,m}$ defined as \eqref{defC}. Based on the definition of the polynomials $P_{s_k,k}$ and $Q_{s_k,k}$ defined in Lemma \ref{lem_varpi}, the elements $\Acc^k_{s_k,l,m}$ are expressed using derivatives of $\varpi_k$ at $\ut_k$. In Section \ref{subsecVarpi}, we present a reliable way to compute the value $\varpi_k^{(n)}(\ut_k)$.
		\item In Section \ref{subsecP}, we compute the polynomials $\Pcc^k_\sigma$ for $\sigma=1,2$. We compare those results with the asymptotic expansion determined in \cite[Theorem 1.2]{R-S}.
		\item In Section \ref{subsecNum}, we compute numerically the polynomials $\Pcc^k_\sigma$ and verify the sharpness of the estimates \eqref{devPrinc} in Theorem \ref{thPrinc} for two specific examples of sequences $a$: 
		
		$\star$ A case where the sequence $a$ has real non negative coefficients. 
		
		$\star$ The sequence $a$ associated to the O3 scheme for the transport equation.
	\end{itemize}
	
	\subsection{Computing the derivatives of $\varpi_k$ at $\ut_k$}\label{subsecVarpi}
	
	The coefficients $\Acc^k_{s_k,l,m}$ defined in \eqref{defC} are expressed using the derivatives of $\varpi_k$ at $\ut_k$. We now present a reliable way to compute $\varpi_k^{(n)}(\ut_k)$. For $\tau\in B_{\varepsilon_\star}(\ut_k)$, $e^{\varpi_k(\tau)}=\kappa_k(e^\tau)$ is an eigenvalue of $\M(e^\tau)$. Lemma \ref{spec_spl} implies that
	\begin{equation}\label{Fvarpi}
		\forall \tau\in B_{\varepsilon_\star}(\ut_k), \quad F(e^{\varpi_k(\tau)})=e^\tau.
	\end{equation}

	For all $n\in \N$, we define the moment function
	\begin{equation}\label{defMn}
		\begin{array}{cccc}
			M_n:& \C^* & \rightarrow & \C\\
			& \kappa & \mapsto & \sum_{j\in \Z}j^na_j\kappa^j
		\end{array}.
	\end{equation}
	We observe that we have the equality $M_0=F$ and
	$$\forall n\in\N, \forall \kappa\in\C^*, \quad M_{n+1}(\kappa)=\kappa\frac{d M_n}{d\kappa}(\kappa),$$
	thus
	\begin{equation}\label{computDer}
		\forall n\in\N, \forall \tau\in B_{\varepsilon_\star}(\ut_k), \quad \frac{d}{d\tau}\left(M_n\left(e^{\varpi_k(\tau)}\right)\right) =\varpi_k^\prime(\tau)M_{n+1}\left(e^{\varpi_k(\tau)}\right).
	\end{equation}

	We will differentiate the equality \eqref{Fvarpi} and use the equality \eqref{computDer} to find an expression of $\varpi_k^{(n)}(\ut_k)$. To do so, we introduce the Bell polynomials (see \cite{Comtet}, Chapter 3.3) defined for $n\in \N$ and $j\in\lc1, \ppp, n\rc$ as 
	$$B_{n,j}(X_1,\ppp,X_{n+1-j}):= \sum \frac{n!}{l_1!\ppp l_{n+1-j}!}\left(\frac{X_1}{1!}\right)^{l_1}\ppp\left(\frac{X_{n+1-j}}{(n+1-j)!}\right)^{l_{n+1-j}}$$
	where the sum is taken over the integers $l_1,\ppp,l_{n+1-j}\in\N$ such that
	\begin{align*}
		j &= l_1+l_2+\ppp+l_{n+1-j},\\
		n& = l_1+2l_2+\ppp+(n+1-j)l_{n+1-j}.
	\end{align*}
	The Bell polynomials $B_{n,j}$ verify the following equalities:
	\begin{equation}
		\forall n \in \N^*, \forall j\in\lc1,\ppp,n\rc, \quad  B_{n,j} = \sum_{i=1}^{n+1-j}\binom{n-1}{i-1}X_i B_{n-i,j-1}, \label{bell1}
	\end{equation}
	\begin{equation}
		\forall n \in \N^*, \forall j\in\lc1,\ppp,n\rc, \forall i \in \lc1,\ppp,n+1-j\rc, \quad  \frac{\partial B_{n,j}}{\partial X_i} = \binom{n}{i} B_{n-i,j-1}.\label{bell2}
	\end{equation}
	
	We can now prove the following lemma which allows us to express recursively the derivatives of $\varpi_k$ at $\ut_k$ with the moments $M_n(\uk_k)$.
	\begin{lemma}\label{expVarpi}
		For all $k\in \lc1,\ppp,K\rc$, we have 
		\begin{align*}
			&\varpi_k^\prime(\ut_k)=\frac{\uz_k}{M_1(\uk_k)},\\
			\forall n\geq 2, \quad & \varpi_k^{(n)}(\ut_k)=\frac{1}{M_1(\uk_k)}\left(\uz_k-\sum_{j=2}^n M_j(\uk_k)B_{n,j}\left(\varpi_k^\prime(\ut_k),\ppp,\varpi_k^{(n+1-j)}(\ut_k)\right)\right).
		\end{align*}
	\end{lemma}
	
	\begin{proof}
		Using the equalities \eqref{computDer}, \eqref{bell1} and \eqref{bell2}, we can prove recursively the following equality for all $n\in\N^*$ and $\tau\in B_{\varepsilon_\star}(\ut_k)$ which looks like Faà di Bruno's formula :
		\begin{equation}\label{derMn}
			\frac{d^n}{d\tau^n}\left(M_0(e^{\varpi_k(\tau)})\right) = \sum_{j=1}^n M_j(e^{\varpi_k(\tau)})B_{n,j}\left(\varpi_k^\prime(\tau),\ppp,\varpi_k^{(n+1-j)}(\tau)\right).
		\end{equation}
		Using the equalities \eqref{derMn}, \eqref{Fvarpi} and $M_0=F$, we conclude the proof of Lemma \ref{expVarpi}.
	\end{proof}
	
	\subsection{Computation of $\Pcc^k_\sigma$ for $\sigma=1,2$}\label{subsecP}
	
	In this section, we will compute the polynomials $\Pcc^k_\sigma$ for $\sigma=1,2$. The goal is to compare the asymptotic expansion \eqref{devPrinc} with the result of \cite[Theorem 1.2]{R-S} and with the local limit theorem (see \cite[Chapter VII, Theorem 13]{Petrov}). We consider $k\in \lc 1,\ppp,K\rc$ and $s_k\in \N^*$. 
	
	$\bullet$ We start to compute the polynomials $\Pcc^k_{1}$. We have using \eqref{Pcc}
	\begin{align*}
		\Pcc^k_1  &= \sum_{l=0}^{s_k-1}\sum_{m=(2\mu_k+1)l}^{(2\mu_k+s_k-1)l+s_k-1}\sum_{k_1=0}^{s_k}\sum_{k_3=0}^{s_k} \ind_{m-2\mu_kl+k_3(2\mu_k-1)+1=1}\Ccc^k_{s_k,l,m,k_1,k_3}X^{k_1+k_3}Y^{m+k_1}\\
		& = \sum_{k_1=0}^{s_k} \Ccc^k_{s_k,0,0,k_1,0}X^{k_1}Y^{k_1}.
	\end{align*}
	Furthermore, for $k_1\in\lc0,\ppp,s_k\rc$, we have using the definition \eqref{defCcc} of $\Ccc^k_{s_k,l,m,k_1,k_3}$ that
	$$\Ccc^k_{s_k,0,0,k_1,0}=\lc\begin{array}{cc}0 & \quad \text{ if }k_1\geq 1,\\ \Acc^k_{s_k,0,0}|\alpha_k| & \quad \text{ if }k_1=0.\end{array}\right. $$
	Furthermore, using the equality \eqref{defC} and the asymptotic expansion \eqref{eqVarpi}, we have 
	$$\Acc^k_{s_k,0,0} = -\mathrm{sgn}(\alpha_k)\varpi_k^\prime(\ut_k) = \frac{1}{|\alpha_k|}.$$
	We then have
	$$\Ccc^k_{s_k,0,0,k_1,0} = \Acc^k_{s_k,0,0}|\alpha_k| = 1.$$
	Therefore, we have proved that
	\begin{equation}\label{Pcc1}
		\forall s_k\in\N\backslash\lc0\rc, \quad \Pcc^k_1=1.
	\end{equation}
	Theorem \ref{thPrinc} implies that there exist two positive constants $C,c$ such that
	\begin{equation}\label{devPrincS1}
		\forall (n,j)\in\N^*\times \Z, \quad \left|\Gcc_j^n - \sum_{k=1}^K\frac{\uz_k^n\uk_k^j}{n^\frac{1}{2\mu_k}}H_{2\mu_k}^{\beta_k}\left(\frac{n\alpha_k-j}{n^\frac{1}{2\mu_k}}\right)\right|\leq \sum_{k=1}^K\frac{C}{n^\frac{1}{\mu_k}}\exp\left(-c\left(\frac{|n\alpha-j|}{n^\frac{1}{2\mu_k}}\right)^\frac{2\mu_k}{2\mu_k-1}\right).
	\end{equation}
	The estimate \eqref{devPrincS1} deduced from Theorem \ref{thPrinc} gives us the same leading term for the asymptotic behavior of $\Gcc_j^n$ as expected from \cite[Theorem 1.2]{R-S}.
	
	$\bullet$ We now compute the polynomials $\Pcc^k_{2}$. We have using \eqref{Pcc}
	\begin{align*}
		\Pcc^k_2  &= \sum_{l=0}^{s_k-1}\sum_{m=(2\mu_k+1)l}^{(2\mu_k+s_k-1)l+s_k-1}\sum_{k_1=0}^{s_k}\sum_{k_3=0}^{s_k} \ind_{m-2\mu_kl+k_3(2\mu_k-1)+1=2}\Ccc^k_{s_k,l,m,k_1,k_3}X^{k_1+k_3}Y^{m+k_1}\\
		& = \sum_{k_1=0}^{s_k} \ind_{\mu_k=1} \Ccc^k_{s_k,0,0,k_1,1}X^{1+k_1}Y^{k_1}+\Ccc^k_{s_k,0,1,k_1,0}X^{k_1}Y^{1+k_1}+\Ccc^k_{s_k,1,2\mu_k+1,k_1,0}X^{k_1}Y^{2\mu_k+1+k_1}.
	\end{align*}

	Furthermore, for $k_1\in\lc0,\ppp,s_k\rc$, we have using the definition \eqref{defCcc} of $\Ccc^k_{s_k,l,m,k_1,k_3}$ that
	$$\Ccc^k_{s_k,0,1,k_1,0}=\lc\begin{array}{cc}0 & \quad \text{ if }k_1\geq 1,\\ \Acc^k_{s_k,0,1}\alpha_k|\alpha_k| & \quad \text{ if }k_1=0,\end{array}\right. $$
	$$\Ccc^k_{s_k,1,2\mu_k+1,k_1,0}=\lc\begin{array}{cc}0 & \quad \text{ if }k_1\geq 1,\\ \Acc^k_{s_k,1,2\mu_k+1}\alpha_k^{2\mu_k+2}|\alpha_k| & \quad \text{ if }k_1=0,\end{array}\right. $$
	\begin{align*}
		\Ccc^k_{s_k,0,0,k_1,1} & = \Acc^k_{s_k,0,0} \alpha_k^{-1}|\alpha_k| \sum_{k_2=0}^{k_1}\binom{k_1}{k_2}(-1)^{k_1-k_2} \frac{k_2+1}{2\mu_k}\\
		& = \lc\begin{array}{cc}0 & \quad \text{ if }k_1\geq 2,\\ \frac{\Acc^k_{s_k,0,0}\mathrm{sgn}(\alpha_k)}{2\mu_k} & \quad \text{ if }k_1=0,1.\end{array}\right. 
	\end{align*}
	Also, using the equality \eqref{defC} and Lemma \ref{expVarpi}, we have 
	$$\Acc^k_{s_k,0,0} = -\mathrm{sgn}(\alpha_k)\varpi_k^\prime(\ut_k),$$
	$$\Acc^k_{s_k,0,1} = -\mathrm{sgn}(\alpha_k)\varpi_k^{(2)}(\ut_k),$$
	$$\Acc^k_{s_k,1,2\mu_k+1} = -\mathrm{sgn}(\alpha_k)\varpi_k^\prime(\ut_k)\frac{\varpi_k^{(2\mu_k+1)}(\ut_k)}{(2\mu_k+1)!}.$$
	Thus, for all $s_k\in\N\backslash\lc0\rc$, 
	\begin{equation}\label{Pcc2}
		 \Pcc^k_2=\ind_{\mu_k=1}\left(-\frac{\varpi_k^\prime(\ut_k)}{2\mu_k}\right)\left(X+X^2Y\right) -\alpha_k^2\varpi_k^{(2)}(\ut_k)Y - \alpha_k^{2\mu_k+3}\varpi_k^\prime(\ut_k)\frac{\varpi_k^{(2\mu_k+1)}(\ut_k)}{(2\mu_k+1)!}Y^{2\mu_k+1}.
	\end{equation}
	Theorem \ref{thPrinc} thus implies that there exist two positive constants $C,c$ such that for all $(n,j)\in\N^*\times \Z$
	\begin{equation}\label{devPrincS2}
		\left|\Gcc_j^n - \uz_k^n\uk_k^j\sum_{k=1}^K\frac{1}{n^\frac{1}{2\mu_k}}H_{2\mu_k}^{\beta_k}\left(X_{n,j,k}\right) + \frac{1}{n^\frac{1}{\mu_k}}\Pcc^k_2\left(X_{n,j,k},\frac{d}{dx}\right)H_{2\mu_k}^{\beta_k}(X_{n,j,k})\right|\leq \sum_{k=1}^K\frac{C}{n^\frac{3}{2\mu_k}}\exp\left(-c\left|X_{n,j,k}\right|^\frac{2\mu_k}{2\mu_k-1}\right)
	\end{equation}
	with $X_{n,j,k}:=\frac{n\alpha_k-j}{n^\frac{1}{2\mu_k}}$.
	
	$\star$ When $\mu_k\geq 2$, the asymptotic expansion \eqref{eqVarpi} implies that $\varpi_k^{(2)}(\ut_k)=0$. Thus, the equality \eqref{Pcc2} becomes
	\begin{equation}\label{Pcc2mu>1}
		\Pcc^k_2=- \alpha_k^{2\mu_k+3}\varpi_k^\prime(\ut_k)\frac{\varpi_k^{(2\mu_k+1)}(\ut_k)}{(2\mu_k+1)!}Y^{2\mu_k+1}.
	\end{equation}
	
	$\star$ We now look at the case $\mu_k=1$. The equality \eqref{Pcc2} becomes
	$$ \Pcc^k_2=\left(-\frac{\varpi_k^\prime(\ut_k)}{2}\right)\left(X+X^2Y\right) -\alpha_k^2\varpi_k^{(2)}(\ut_k)Y - \alpha_k^{5}\varpi_k^\prime(\ut_k)\frac{\varpi_k^{(3)}(\ut_k)}{6}Y^{3}.$$
	As we said in the introduction of the paper, the polynomials satisfying Theorem \ref{thPrinc} are not unique. We will now propose a more convenient choice of polynomials to replace $\Pcc_{k,s_k,2}$. Using Lemma \ref{lemH}, if we define the polynomial
	$$\Qc^k_2(X,Y)= \left(-2\beta_k\varpi_k^\prime(\ut_k)-\alpha_k^2\varpi_k^{(2)}(\ut_k)\right) Y  + \left(-2\beta_k^2\varpi_k^\prime(\ut_k)-\alpha_k^{5}\varpi_k^\prime(\ut_k)\frac{\varpi_k^{(3)}(\ut_k)}{6}\right) Y^3\in\C[X,Y]$$
	we have 
	\begin{equation}\label{egPccQc}
		\Pcc^k_2\left(.,\frac{d}{dx}\right)H_{2}^{\beta_k} = \Qc^k_2\left(.,\frac{d}{dx}\right)H_{2}^{\beta_k}.
	\end{equation}
	We can then replace $\Pcc^k_2$ with $\Qc^k_2$ in the estimate \eqref{devPrincS2} when $\mu_k=1$. This allows us to express the second term of the asymptotic expansion using a linear combination of derivatives of $H_{2}^{\beta_k}$ since $\Qc^k_2(X,Y)$ does not have any terms where $X$ intervenes. We notice that the asymptotic expansion \eqref{eqVarpi} implies that
	\begin{equation}\label{egBeta}
		\frac{\varpi_k^{(2)}(\ut_k)}{2}=\frac{\beta_k}{\alpha_k^3}.
	\end{equation}
	Using Lemma \ref{expVarpi} and equality \eqref{egBeta}, we can prove that actually 
	\begin{equation}\label{egQc}
		\Qc^k_2(X,Y)=-\frac{1}{6\uz_k^2}\left(\uz_k^2M_3(\uk_k)-3\uz_kM_2(\uk_k)M_1(\uk_k)+2M_1(\uk_k)^3\right)Y^3.
	\end{equation}
	We will see in Section \ref{subsecNumProba} that, in the probabilistic case we presented in the introduction of the paper that motivated our result, this expression of $\Qc^k_2$ gives exactly the second term of the asymptotic expansion \eqref{cas_proba} when we apply the local limit theorem (which is fortunate).

	\subsection{Numerical examples}\label{subsecNum}
	
	In this section, we consider some examples of elements $a\in\ell^1(\Z)$ which satisfy the conditions of Theorem \ref{thPrinc} and see how sharp the estimations we found are.
	
	\subsubsection{Probability distribution : real non negative sequences}\label{subsecNumProba}
	
	First, we consider the case where $a$ has real non negative coefficients. If we introduce the sequence $b=(a_{-j})_{j\in\Z}$, then $b$ is the probability distribution of some random variable $X$ supported on $\Z$. We observe that $L_b=\Lcc_a$, so, recalling that $b^n=b\ast \ppp\ast b$, we have
	$$\forall n\in\N^*,\forall j\in\Z, \quad b^n_j=\Gcc_j^n.$$
	We will settle on $a\in\ell^1(\Z)$ such that $a_j=0$ for $j\neq -1,0,1$ and 
	$$a_{-1}=2/3, a_0=1/6, a_1=1/6.$$
	This sequence verifies Hypothesis \ref{H1}. In this case, we have $r=p=1$. Also, $F(1)=1$ and 
	$$\forall \kappa\in\S^1\backslash\lc1\rc, \quad |F(\kappa)|<1.$$
	The function $F$ satisfies that
	$$F(e^{i\xi})\underset{\xi\rightarrow 0}=\exp(-i\alpha\xi -\beta\xi^2+o(\xi^2))$$
	where $\alpha=\E(X)=\frac{1}{2}$ and $\beta=\frac{V(X)}{2}=\frac{7}{24}$. We have $\mu=1$ in this case and Hypothesis \ref{H2} is satisfied with $K=1$, $\uk_1=1$ and $\uz_1=1$. It also directly satisfies Hypothesis \ref{H3} since $K=1$. Since $K=1$, we lose the subscript $k$ in most the notations that follow. The sequence $a$ verifies Hypotheses \ref{H1}, \ref{H2} and \ref{H3}, so we can apply Theorem \ref{thPrinc}. As an example, we will apply Theorem \ref{thPrinc} for $s=2$ and use the calculations of Section \ref{subsecP} to determine the terms of the asymptotic expansion:
	
	$\bullet$ Using the equality \eqref{Pcc1} on $\Pcc^k_1$, the leading order term of the asymptotic expansion given by Theorem \ref{thPrinc} is
	\begin{align*}
		\frac{1}{\sqrt{n}}H_{2}^\beta\left(\frac{n\alpha-j}{\sqrt{n}}\right) &=\frac{1}{\sqrt{4\pi\beta n}}\exp\left(-\frac{|j-n\alpha|^2}{4\beta n}\right)\\
		&=\frac{1}{\sqrt{2\pi V(X) n}}\exp\left(-\frac{|j-n\E(X)|^2}{2V(X) n}\right).
	\end{align*}

	$\bullet$ We notice that using the moments function $M_n$ defined with \eqref{defMn}, we have
	$$\forall n\in \N, \quad M_n(1)=(-1)^n\E(X^n).$$
	Using the equalities \eqref{egPccQc} and \eqref{egQc} that respectively links the polynomials $\Pcc^k_2$ and $\Qc^k_2$ and allows us to compute the polynomial $\Qc^k_2$, the second order term of the asymptotic expansion given by Theorem \ref{thPrinc} is 
	\begin{align*}
		&\frac{1}{n}\left(-\frac{1}{6}(M_3(1)-3M_2(1)M_1(1)+2M_1(1)^3)\right)\left(H_{2}^\beta\right)^{(3)}\left(X_{n,j}\right)\\
		 &=\frac{\E((X-\E(X))^3)}{6(2\beta)^2n}\left(H_{2}^\frac{1}{2}\right)^{(3)}\left(\frac{X_{n,j}}{\sqrt{2\beta}}\right)\\
		 &=\frac{q_1\left(\frac{X_{n,j}}{\sqrt{V(X)}}\right)}{n}
	\end{align*}
	where $X_{n,j}:=\frac{n\E(X)-j}{\sqrt{n}}$ and the function $q_1:\R\rightarrow \R$ is defined as
	$$\forall x\in\R,\quad q_1(x):=-\frac{\E((X-\E(X))^3)}{6\sqrt{2\pi}V(X)^2}(x^3-3x)e^{-\frac{x^2}{2}}. $$
	
	Theorem \ref{thPrinc} then states that there exist two constants $C,c>0$ such that
	\begin{equation}\label{ine_proba}
		\forall n\in\N^*,\forall j\in\Z,\quad \left|\mathrm{Err}(n,j)\right|\leq  \frac{C}{n^\frac{3}{2}}\exp\left(-c\left|X_{n,j}\right|^2\right),
	\end{equation}
	with $X_{n,j}=\frac{n\E(X)-j}{\sqrt{n}}$ and 
	$$\mathrm{Err}(n,j):=\Gcc_j^n-\frac{1}{\sqrt{2\pi V(X) n}}\exp\left(-\frac{|X_{n,j}|^2}{2V(X)}\right)-\frac{q_1\left(\frac{X_{n,j}}{\sqrt{V(X)}}\right)}{n}$$.
	
\begin{figure}
		\begin{center}
			\includegraphics[width=160mm]{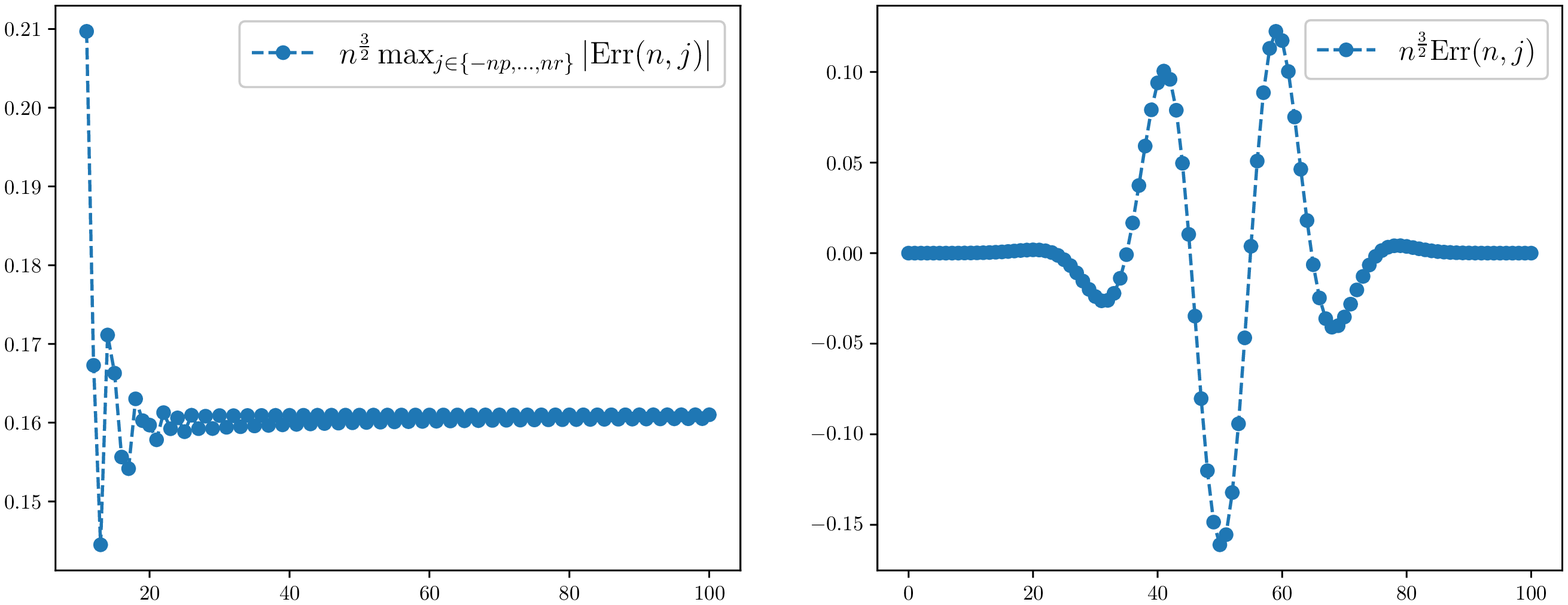}
			\caption{On the left : A representation of $n^\frac{3}{2}\max_{j\in\lc-nr,\ppp,np\rc} |\mathrm{Err}(n,j)|$ depending on $n$. As expected knowing that $-r<\alpha<p$, we see that the function is bounded and even seems to converge. On the right : We fixed $n=100$ and represented $j\in\Z\mapsto n^\frac{3}{2}\mathrm{Err}(n,j)$.}
			\label{im_cas_proba}
		\end{center}
	\end{figure}

	The estimate \eqref{ine_proba} is exactly the asymptotic expansion of the elements $b^n_j=\Gcc_j^n$ we expected via the local limit theorem (see \cite[Chapter VII, Theorem 13]{Petrov} for more details). This behavior is represented on Figure \ref{im_cas_proba} where we even see that the remainder $n^\frac{3}{2}Err(n,j)$ seems to scale like $f\left(\frac{n\alpha-j}{\sqrt{n}}\right)$. This would correspond to the next term in the asymptotic expansion of $\Gcc_j^n$.
	
	\subsubsection{The O3 scheme for the transport equation}
	
	We will now consider an example linked to finite difference schemes. We consider the transport equation 
	$$\partial_t u + a \partial_x u =0, \quad (t,x)\in\R_+\times \R$$
	with Cauchy data at $t=0$. The O3 scheme is an explicit third order accurate finite difference approximation of the previous transport equation. We refer to \cite{Despres} for a detailed analysis of this scheme. It corresponds to the numerical scheme \eqref{an_num} for $a\in\ell^1(\Z)$ such that $a_j=0$ for $j\notin \lc-2,-1,0,1\rc$ and 
	$$a_{-2}=-\frac{\lambda a (1-(\lambda a)^2)}{6} , \quad a_{-1} = \frac{\lambda a (1+\lambda a)(2-\lambda a)}{2},\quad		a_0=\frac{ (1-(\lambda a)^2)(2-\lambda a)}{2} ,\quad a_1 =- \frac{\lambda a (1-\lambda a)(2-\lambda a)}{6},$$
	with $\lambda =\frac{\Delta t}{\Delta x}>0$. The parameter $\lambda a$ is the Courant number. We have in this case that $r=2$ and $p=1$. For $\lambda a\in]-1,1[\backslash\lc0\rc$, we have that $F(1)=1$ and
	$$\forall \kappa\in\S^1\backslash\lc1\rc, \quad |F(\kappa)|<1.$$
	Also, there exists $\beta\in\R_+^*$ such that 
	$$F(e^{i\xi})\underset{\xi\rightarrow 0}=\exp(-i\lambda a \xi - \beta \xi^4+o(\xi^4)).$$
	We have $\mu=2$ in this case and Hypothesis \ref{H2} is satisfied with $K=1$, $\uk_1=1$ and $\uz_1=1$. Since $K=1$, we lose the subscript $k$ in most the notations that follow. The sequence $a$ verifies hypotheses \ref{H1}, \ref{H2} and \ref{H3}, so we can apply Theorem \ref{thPrinc}. As an example, we will apply Theorem \ref{thPrinc} for $s=3$ and $\lambda a=\frac{1}{2}$. 
	
	$\bullet$ Using the equality \eqref{Pcc1}, we have
	$$\Pcc_1=1.$$
	
	$\bullet$ Using the equality \eqref{Pcc2mu>1} and Lemma \ref{expVarpi} to compute $\varpi^{(3)}(1)$, we have
	$$\Pcc_2=0.$$
	
	$\bullet$ Using the equality \eqref{Pcc} to express the polynomial $\Pcc_3$ and Lemma \ref{expVarpi} to compute the coefficients $\Acc_{l,m}$, we numerically compute the polynomials $\Pcc_{3}$:
	$$\Pcc_3=p_3Y^6$$
	where $p_3\approx-1,953125.10^{-3}$.
	
	Theorem \ref{thPrinc} then states that there exist two constants $C,c>0$ such that
	\begin{equation}\label{ine_O3}
		\forall n\in\N^*,\forall j\in\Z,\quad \left|\mathrm{Err}(n,j)\right|\leq  \frac{C}{n}\exp\left(-c\left|X_{n,j}\right|^\frac{4}{3}\right),
	\end{equation}
	with $X_{n,j}=\frac{n\alpha-j}{n^\frac{1}{4}}$ and 
	$$\mathrm{Err}(n,j):=\Gcc_j^n-\sum_{\sigma=1}^3\frac{1}{n^\frac{\sigma}{4}}\Pcc_{\sigma}\left(X_{n,j},\frac{d}{dx}\right)H_{4}^\beta(X_{n,j})$$.
	
	\begin{figure}
		\begin{center}
			\includegraphics[width=160mm]{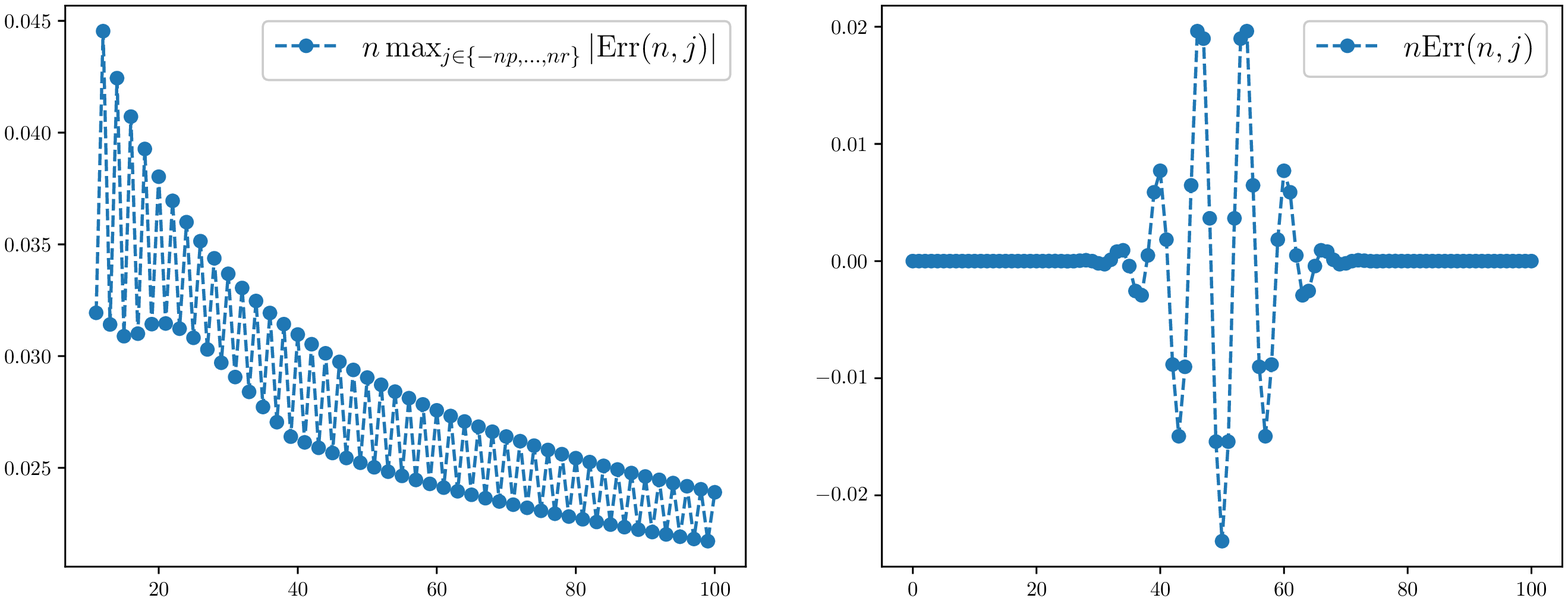}
			\caption{For these figures, we chose $\lambda a=1/2$. On the left : A representation of $n\max_{j\in\lc-nr,\ppp,np\rc} |\mathrm{Err}(n,j)|$ depending on $n$. As expected, the function seems to be bounded. On the right : We fixed $n=100$ and represented $j\in\Z\mapsto n\mathrm{Err}(n,j)$. We observe the exponential decay in $j$. Also, we can see a particular shape of curve that arises that would correspond to the next term in the asymptotic expansion of $\Gcc_j^n$.}
			\label{im_O3}
		\end{center}
	\end{figure}
	
	This behavior is represented on Figure \ref{im_O3} where we even see that the remainder $nErr(n,j)$ seems to scale like $f\left(\frac{n\alpha-j}{\sqrt{n}}\right)$. Hence, the estimate \eqref{ine_O3} seems to be sharp.

	\section{Appendix: Proof of auxiliary results}\label{sec_appendix}
	
	\subsection{Proof of the Lemma \ref{arp}}
	
	We recall here the statement of Lemma \ref{arp}. 
	
	\begin{lemma*}[Lemma \ref{arp}]
		For $a\in\ell^1(\Z)$ which verifies Hypotheses \ref{H1} and \ref{H2}, we have that $a_{-r}$ and $a_p$ belong to $\D$.
	\end{lemma*}
	
	\begin{proof}
		We introduce the polynomial function $g$ defined by
		$$\forall \kappa\in\C, \quad g(\kappa):= \sum_{l=-r}^{p} a_l\kappa^{l+r}.$$
		For all $\kappa\in\S^1$, Hypothesis \ref{H2} implies that 
		$$|g(\kappa)|=\left|\kappa^rF(\kappa)\right|=|F(\kappa)|\leq 1.$$
		Observing that $g$ is not a constant function, the maximum principle for holomorphic functions \cite{Rudin} allows us to conclude that 
		$$|a_{-r}|=|g(0)|<1.$$
		The same kind of argument allows us to conclude for the coefficient $a_p$.
	\end{proof}
	
	\subsection{Proof of the Lemma \ref{ineg_H}}
	
	We recall here the statement of Lemma \ref{ineg_H}. 
	
	\begin{lemma*}[Lemma \ref{ineg_H}]
		For $\mu\in\N^*$, $\beta\in \C$ with positive real part and $m\in \N$, there exist two constants $C,c>0$ such that
		$$\forall x\in \R, \quad  \left|{H_{2\mu}^\beta}^{(m)}(x)\right|\leq C\exp\left(-c|x|^\frac{2\mu}{2\mu-1}\right).$$
	\end{lemma*}

	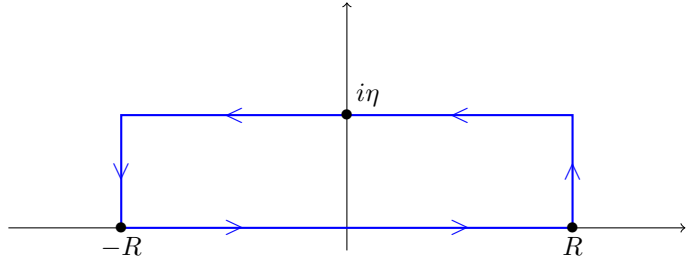
\begin{figure}
		\begin{center}
			\begin{tikzpicture}[scale=3]
				\draw (-1.5,0) -- (-1,0);
				\draw[->] (1,0) -- (1.5,0);
				\draw[->] (0,-0.1) -- (0,1);
				\draw[thick,color=blue] (-1,0) -- (1,0) node[near start, sloped] {$>$} node[near end, sloped] {$>$} -- (1,0.5)node[midway, sloped] {$>$} -- (-1,0.5) node[near start, sloped] {$<$} node[near end, sloped] {$<$} -- cycle node[midway, sloped] {$>$};
				\draw (-1,0) node {$\bullet$};
				\draw (-1,0) node[below] {$-R$};
				\draw (1,0) node {$\bullet$};
				\draw (1,0) node[below] {$R$};
				\draw (0,0.5) node {$\bullet$};
				\draw (0,0.5) node[above right] {$i\eta$};
			\end{tikzpicture}
			\caption{Integrating path for the proof of Lemma \ref{ineg_H}.}
			\label{int_path}
		\end{center}
	\end{figure}
	
	\begin{proof}
		We fix $\eta\in\R$ that we will choose more precisely later. Integrating the function $z\mapsto (iz)^m\exp(izx-\beta z^{2\mu})$ on the rectangle depicted in the Figure \ref{int_path} using the Cauchy formula and passing to the limit $R\rightarrow +\infty$, we obtain
		$$\forall \eta\in\R, \quad {H_{2\mu}^\beta}^{(m)}(x) = \frac{1}{2\pi}\int_\R (i(t+i\eta))^me^{i(t+i\eta)x}e^{-\beta (t+i\eta)^{2\mu}}dt.$$
		Thus,
		$$\left|{H_{2\mu}^\beta}^{(m)}(x)\right| \leq \frac{e^{-\eta x}}{2\pi}\int_\R (t^2+\eta^2)^\frac{m}{2}\exp\left(-\Re\left(\beta (t+i\eta)^{2\mu}\right)\right)dt.$$
		Using Young's inequality, we can show that there exists a constant $c>0$ such that 
		$$\forall t\in\R, \quad\Re\left(\beta (t+i\eta)^{2\mu}\right)\geq \frac{\Re(\beta)}{2}t^{2\mu}-c \eta^{2\mu}.$$
		and thus there exists $C>0$ independent from $x$ and $\eta$ such that 
		$$\left|{H_{2\mu}^\beta}^{(m)}(x)\right| \leq C(1+|\eta|^m) e^{-\eta x+c\eta^{2\mu}}.$$
		Optimizing $e^{-\eta x+c\eta^{2\mu}}$ with respect to $\eta$ yields the desired result.
	\end{proof}

	\subsection{Proof of the Lemma \ref{est_Beps}}
	
	We recall here the statement of Lemma \ref{est_Beps}. 
	
	\begin{lemma*}[Lemma \ref{est_Beps}, Inequalities in $B_{\varepsilon_\star}(\ut_k)$]
		There exists $C>0$ such that for all $\tau \in B_{\varepsilon_\star}(\ut_k)$ and $(n,j)\in\Dc_k$, we have 
		$$ \left|e^{n\tau} \left( e^{j\varpi_k(\tau)} -e^{jQ_{s,k}(\tau)}\right)\right| \leq C n|\tau-\ut_k|^{2\mu_k+s} \exp(n \Re(\tau-\ut_k)+j (\Re(\varpi_k(\tau)) + \mathrm{sgn}(\alpha_k)|\xi_{s,k}(\tau)(\tau-\ut_k)^{2\mu_k+s}|)) $$
		and
		$$\left|e^{n\tau+j\varphi_k(\tau)} \left( e^{jR_{s,k}(\tau)} -\sum_{l=0}^{s-1}\frac{(jR_{s,k}(\tau))^l}{l!}\right)\right| \leq C \left(n|\tau-\ut_k|^{2\mu_k+1}\right)^{s} \exp(n \Re(\tau-\ut_k)+j (\Re(\varphi_k(\tau)) + \mathrm{sgn}(\alpha_k)|R_{s,k}(\tau)|)). $$
	\end{lemma*}
	
	\begin{proof}
		We begin with the first inequality. We define the holomorphic function $\Sc$ such that $$\forall z\in \C, \quad \Sc(z) = \lc\begin{array}{cc} 1 & \text{ if }z=0, \\ \frac{\sinh(z)}{z} & \text{ else.} \end{array}\right.$$
		
		We consider $(n,j)\in\Dc_k$ and $\tau \in B_{\varepsilon_\star}(\ut_k)$. We have 
		\begin{multline*}
			\left|e^{n\tau} \left( e^{j\varpi_k(\tau)} - e^{jQ_{s,k}(\tau)}\right)\right| 
			= |j||\xi_{s,k}(\tau)||\tau-\ut_k|^{2\mu_k+s}\left|\Sc\left(j\frac{\xi_{s,k}(\tau)(\tau-\ut_k)^{2\mu_k+s}}{2}\right)\right|\\ \exp\left(n\Re(\tau) + j \left(\Re(\varpi_{k}(\tau))-\Re\left( \frac{\xi_{s,k}(\tau)(\tau-\ut_k)^{2\mu_k+s}}{2}\right)\right)\right).
		\end{multline*}
		
		We observe that the function $z\in\C\mapsto |\Sc(z)|\exp(-|z|)$ is bounded. Therefore, because the function $\xi_k$ can be bounded on $B_{\varepsilon_\star}(\ut_k)$ and $\Re(\tau)= \Re(\tau-\ut_k)$,  
		\begin{multline*}
			\left|e^{n\tau} \left( e^{j\varpi_k(\tau)} - \uk_k^je^{j\varphi_k(\tau)}\right)\right| \\ \lesssim n|\tau-\ut_k|^{2\mu_k+s} \exp\left(n\Re(\tau-\ut_k) + j \left(\Re(\varpi_k(\tau))-\Re\left( \frac{\xi_{s,k}(\tau)(\tau-\ut_k)^{2\mu_k+s}}{2}\right)\right)+|j|\frac{|\xi_{s,k}(\tau)(\tau-\ut_k)^{2\mu_k+s}|}{2}\right).
		\end{multline*}
		Since we have
		\begin{multline*}
			j \left(\Re(\varpi_k(\tau))-\Re\left( \frac{\xi_{s,k}(\tau)(\tau-\ut_k)^{2\mu_k+s}}{2}\right)\right)+|j|\frac{|\xi_{s,k}(\tau)(\tau-\ut_k)^{2\mu_k+s}|}{2} \\ \leq j\left(\Re(\varpi_k(\tau))+\mathrm{sgn}(\alpha_k)|\xi_{s,k}(\tau)(\tau-\ut_k)^{2\mu_k+s}|\right).
		\end{multline*}
		
		The proof of the second inequality is similar. We define the holomorphic function $\Psi_s$ such that $$\forall z\in \C, \quad \Psi_s(z) = \frac{1}{z^s}\left(e^z-\sum_{l=0}^{s-1}\frac{z^l}{l!}\right).$$
		We then have
		$$e^{n\tau+j\varphi_k(\tau)} \left( e^{jR_{s,k}(\tau)} -\sum_{l=0}^{s-1}\frac{(jR_{s,k}(\tau))^l}{l!}\right)=(jR_{s,k}(\tau))^s\Psi_s(jR_{s,k}(\tau))e^{n\tau+j\varphi_k(\tau)} .$$
		We observe that the function $z\in\C\mapsto |\Psi_s(z)|\exp(-|z|)$ is bounded. Therefore,
		$$\left|e^{n\tau+j\varphi_k(\tau)} \left( e^{jR_{s,k}(\tau)} -\sum_{l=0}^{s-1}\frac{(jR_{s,k}(\tau))^l}{l!}\right)\right|\lesssim (n|\tau-\ut_k|^{2\mu_k+1})^{s}\exp(n\Re(\tau-\ut_k) +j\Re(\varphi_k(\tau))+|jR_{s,k}(\tau)|).$$
		We can then conclude the proof of the second inequality.
	\end{proof}

	\subsection{Proof of the Lemma \ref{ine_taup}}
	
	We recall here the statement of Lemma \ref{ine_taup}.
	
	\begin{lemma*}[Lemma \ref{ine_taup}, Inequalities on $\Gamma_{k,p}$]
		For $(n,j)\in\N^*\times \Z$ such that $\mathrm{sgn}(j)=\mathrm{sgn}(\alpha_k)$ and $\tau\in \Gamma_{k,p}$, we have
		
		$\bullet$ Case A: $\rho_k\left(\frac{\zeta_k}{\gamma_k}\right) \in \left[-\frac{\eta}{2},\varepsilon_{k,0}\right]$
		\begin{align*}
			n\Re(\tau-\ut_k)+j(\Re(\varpi_k(\tau)) + \mathrm{sgn}(\alpha_k)|\xi_{s,k}(\tau)(\tau-\ut_k)^{2\mu_k+s}|) &\leq -nc_\star\Im(\tau-\ut_k)^{2\mu_k} - \frac{n}{\alpha_k} (2\mu_k-1)\gamma_k\left(\frac{|\zeta_k|}{\gamma_k}\right)^\frac{2\mu_k}{2\mu_k-1},\\
			n\Re(\tau-\ut_k)+j\Re(\varpi_k(\tau)) &\leq -nc_\star\Im(\tau-\ut_k)^{2\mu_k} - \frac{n}{\alpha_k} 	(2\mu_k-1)\gamma_k\left(\frac{|\zeta_k|}{\gamma_k}\right)^\frac{2\mu_k}{2\mu_k-1},\\
			n\Re(\tau-\ut_k)+j(\Re(\varphi_k(\tau)) + \mathrm{sgn}(\alpha_k)|R_{s,k}(\tau)|) &\leq -nc_\star\Im(\tau-\ut_k)^{2\mu_k} - \frac{n}{\alpha_k} 	(2\mu_k-1)\gamma_k\left(\frac{|\zeta_k|}{\gamma_k}\right)^\frac{2\mu_k}{2\mu_k-1}.
		\end{align*}	
		
		$\bullet$ Case B: $\rho_k\left(\frac{\zeta_k}{\gamma_k}\right) >\varepsilon_{k,0}$
		\begin{align*}
			n\Re(\tau-\ut_k)+j(\Re(\varpi_k(\tau)) + \mathrm{sgn}(\alpha_k)|\xi_{s,k}(\tau)(\tau-\ut_k)^{2\mu_k+s}|)&\leq -\frac{n}{\alpha_k}(2\mu_k-1)A_R\ud_k\varepsilon_{k,0}^{2\mu_k},\\
			n\Re(\tau-\ut_k)+j\Re(\varpi_k(\tau)) &\leq -\frac{n}{\alpha_k}(2\mu_k-1)A_R\ud_k\varepsilon_{k,0}^{2\mu_k},\\
			n\Re(\tau-\ut_k)+j(\Re(\varphi_k(\tau)) + \mathrm{sgn}(\alpha_k)|R_{s,k}(\tau)|) & \leq -\frac{n}{\alpha_k}(2\mu_k-1)A_R\ud_k\varepsilon_{k,0}^{2\mu_k}.
		\end{align*}
		
		$\bullet$ Case C: $\rho_k\left(\frac{\zeta_k}{\gamma_k}\right) <-\frac{\eta}{2}$
		\begin{align*}
			n\Re(\tau-\ut_k)+j(\Re(\varpi_k(\tau)) + \mathrm{sgn}(\alpha_k)|\xi_{s,k}(\tau)(\tau-\ut_k)^{2\mu_k+s}|) &\leq -\frac{n}{\alpha_k}(2\mu_k-1)A_R\ud_k\left(\frac{\eta}{2}\right)^{2\mu_k},\\
			n\Re(\tau-\ut_k)+j\Re(\varpi_k(\tau)) &\leq -\frac{n}{\alpha_k}(2\mu_k-1)A_R\ud_k\left(\frac{\eta}{2}\right)^{2\mu_k},\\
			n\Re(\tau-\ut_k)+j(\Re(\varphi_k(\tau)) + \mathrm{sgn}(\alpha_k)|R_{s,k}(\tau)|) & \leq-\frac{n}{\alpha_k}(2\mu_k-1)A_R\ud_k\left(\frac{\eta}{2}\right)^{2\mu_k}.
		\end{align*}
	\end{lemma*}
	
	\begin{proof}
		In every case, the second inequality is a direct consequence of the first one. Furthermore, the proof of the first and third inequalities are very similar. For the first one, we will use inequality \eqref{estVarpi} and the third one will rely on inequality \eqref{estR}. Thus, we will focus in each case on the first inequality.
		
		We consider $(n,j)\in\N^*\times \Z$ such that $\mathrm{sgn}(j)=\mathrm{sgn}(\alpha_k)$ and $\tau \in \Gamma_{k,p}$. Using first the inequality \eqref{estVarpi}, the fact that $\tau\in\Gamma_{k,p}$ and finally the inequality \eqref{ine_Re}, we have
		\begin{align*}
			n\Re(\tau-\ut_k)+j(\Re(\varpi_k(\tau))+ \mathrm{sgn}(\alpha_k)|\xi_{s,k}(\tau)(\tau-\ut_k)^{2\mu_k+s}|)&
			\leq   n\Re(\tau-\ut_k)- \frac{j}{\alpha_k} \Psi_k(\tau_p) \\
			&\leq -nc_\star \Im(\tau-\ut_k)^{2\mu_k} +\frac{n}{\alpha_k}\left(\gamma_k \tau_p^{2\mu_k} - 2\mu_k \zeta_k\tau_p \right).
		\end{align*}
		
		$\bullet$ First, we consider the case A. Then, we have $\tau_p= \rho_k\left(\frac{\zeta_k}{\gamma_k}\right)$. Therefore, 
		\begin{equation}
			\gamma_k\tau_p^{2\mu_k} - 2\mu_k\zeta_k\tau_p = -(2\mu_k-1)\gamma_k \left(\frac{|\zeta_k|}{\gamma_k}\right)^\frac{2\mu_k}{2\mu_k-1}\leq0. \label{cas1}
		\end{equation}
		
		$\bullet$ We consider the case B. Because $\tau_p= \varepsilon_{k,0}$, we have 
		$$n\Re(\tau-\ut_k)+j(\Re(\varpi_k(\tau))+ \mathrm{sgn}(\alpha_k)|\xi_{s,k}(\tau)(\tau-\ut_k)^{2\mu_k+s}|) \leq\frac{n}{\alpha_k}\left(\gamma_k \varepsilon_{k,0}^{2\mu_k} - 2\mu_k \zeta_k\varepsilon_{k,0} \right) .$$
		
		We recall that $\rho_k\left(\frac{\zeta_k}{\gamma_k}\right)>\varepsilon_{k,0}$ and that $\rho_k\left(\frac{\zeta_k}{\gamma_k}\right)$ is the only real root of $-\zeta_k+\gamma_k x^{2\mu_k-1}=0$. Therefore, $-\zeta_k\leq -\gamma_k\varepsilon_{k,0}^{2\mu_k-1}$ and
		\begin{equation}
			\gamma_k\tau_p^{2\mu_k} - 2\mu_k\zeta_k\tau_p \leq -(2\mu_k-1)\gamma_k \varepsilon_{k,0}^{2\mu_k}\leq0.\label{cas2}
		\end{equation}
		
		Using \eqref{ineg_gamma} to bound $\gamma_k$, we deduce the inequality \eqref{ine_varpi_cas2}.
		
		$\bullet$ Finally, we place ourselves in case C. We have that $\tau_p= -\frac{\eta}{2}$, so
		$$n\Re(\tau-\ut_k)+j(\Re(\varpi_k(\tau))+\mathrm{sgn}(\alpha_k) |\xi_{s,k}(\tau)(\tau-\ut_k)^{2\mu_k+s}|) \leq\frac{n}{\alpha_k}\left(\gamma_k \left(\frac{\eta}{2}\right)^{2\mu_k} + 2\mu_k \zeta_k\frac{\eta}{2}\right) .$$
		
		We recall that $\rho_k\left(\frac{\zeta_k}{\gamma_k}\right)<-\frac{\eta}{2}$ and that $\rho_k\left(\frac{\zeta_k}{\gamma_k}\right)$ is the only real root of $-\zeta_k+\gamma_k x^{2\mu_k-1}=0$. Then,
		$\zeta_k \leq -\gamma_k\left(\frac{\eta}{2}\right)^{2\mu_k-1}$ and
		\begin{equation}
			\gamma_k\tau_p^{2\mu_k} - 2\mu_k\zeta_k\tau_p \leq -(2\mu_k-1)\gamma_k \left(\frac{\eta}{2}\right)^{2\mu_k}\leq0.\label{cas3}
		\end{equation}
		
		Using \eqref{ineg_gamma} to bound $\gamma_k$, we deduce the inequality \eqref{ine_varpi_cas3}.
	\end{proof}

	\subsection{Proof of the Lemma \ref{ine_res}}
	
	We recall here the statement of Lemma \ref{ine_res}. 
	
	\begin{lemma*}[Lemma \ref{ine_res}]
		For $(n,j)\in\N^*\times \Z$ such that $\mathrm{sgn}(j)=\mathrm{sgn}(\alpha_k)$ and $\tau\in \Gamma_{k,res}$, we have in all cases
		\begin{align*}
			n\Re(\tau-\ut_k)+j(\Re(\varpi_k(\tau))+ \mathrm{sgn}(\alpha_k) \left|\xi_{s,k}(\tau)(\tau-\ut_k)^{2\mu_k+s}\right|) &\leq -n\frac{\eta}{2}, \\
			n\Re(\tau-\ut_k)+j(\Re(\varpi_k(\tau)) &\leq -n\frac{\eta}{2}, \\
			n\Re(\tau-\ut_k)+j(\Re(\varphi_k(\tau)) + \mathrm{sgn}(\alpha_k)|R_{s,k}(\tau)|)&\leq -n\frac{\eta}{2}. 
		\end{align*}
	\end{lemma*}

	\begin{proof}
		For the same reasons as for the proof of Lemma \ref{ine_taup}, we will only focus on the first inequality. We consider $(n,j)\in\N^*\times \Z$ such that $\mathrm{sgn}(j)=\mathrm{sgn}(\alpha_k)$ and $\tau \in \Gamma_{k,res}$. Using the inequality \eqref{estVarpi} and the facts that $\Im(\tau-\ut_k)^{2\mu_k}\geq\ell_{k,p}^{2\mu_k}$ and  $-\eta+i\ell_{k,p}+\ut_k\in \Gamma_{k,p}$, we have
		\begin{align*}
			n\Re(\tau-\ut_k)+j(\Re(\varpi_k(\tau))+ \mathrm{sgn}(\alpha_k)\left|\xi_{s,k}(\tau)(\tau-\ut_k)^{2\mu_k+s}\right|)&\leq -n\eta -\frac{j}{\alpha_k}\left( -\eta -A_R\eta^{2\mu_k} + A_I\Im(\tau-\ut_k)^{2\mu_k}\right)\\
			&\leq -n\eta -\frac{j}{\alpha_k} \Psi_k(\tau_p).
		\end{align*}
		
		We know that $\eta+\tau_p\geq \frac{\eta}{2}$, so
		$$-n\eta- \frac{j}{\alpha_k}(\tau_p - A_R\tau_p^{2\mu_k}) = -n(\eta+\tau_p)+\frac{n}{\alpha_k}\left(\gamma_k\tau_p^{2\mu_k}-2\mu_k\zeta_k\tau_p\right) \leq -n\frac{\eta}{2} + \frac{n}{\alpha_k}\left(\gamma_k\tau_p^{2\mu_k} - 2\mu_k\zeta_k\tau_p\right).$$
		
		We proved at the end of the proof of Lemma \ref{ine_taup} that, in the three cases A, B and C, $\gamma_k\tau_p^{2\mu_k} - 2\mu_k\zeta_k\tau_p$ are non positive (see \eqref{cas1}, \eqref{cas2} and \eqref{cas3}). This concludes the proof.
	\end{proof}

	\textbf{Acknowledgments : } The author would like to thank Jean-François Coulombel and Grégory Faye for their many useful advice and suggestions as well as their attentive reading of the paper. He also would like to thank the referees for their numerous comments and their suggestion to search for an asymptotic expansion up to any order and not only up to order $1$.
	
	\bibliographystyle{alpha}
	\bibliography{Coeuret_2022_locLimTh}{}

\end{document}